\documentclass[12pt,a4paper]{amsart}
\theoremstyle{plain}

\usepackage[colorlinks]{hyperref}
\usepackage{enumerate, amssymb}

\advance\hoffset-5mm \advance\textwidth40mm

\def\bdi{\begin{diagram}}
\def\edi{\end{diagram}}

\theoremstyle{plain}

\newtheorem{thm}{Theorem}[section]
\newtheorem{cor}[thm]{Corollary}
\newtheorem{lem}[thm]{Lemma}
\newtheorem{prop}[thm]{Proposition}
\theoremstyle{definition}
\newtheorem{defi}[thm]{Definition}
\newtheorem{defis}[thm]{Definitions}
\newtheorem{conj}[thm]{Conjecture}
\newtheorem{conv}[thm]{Convention}
\newtheorem{nota}[thm]{Notation}
\newtheorem{rem}[thm]{Remark}
\newtheorem{rems}[thm]{Remarks}
\newtheorem{exa}[thm]{Example}
\newtheorem{exas}[thm]{Examples}
\newtheorem{prob}[thm]{Problem}
\newtheorem{probs}[thm]{Problems}
\newtheorem{ques}[thm]{Question}
\newtheorem{sit}[thm]{}

\newcommand{\AVF}{ \operatorname{{\rm AVF}}}
\newcommand{\HVF}{ \operatorname{{\rm HVF}}}
\newcommand{\Span}{ \operatorname{{\rm Span}}}

\newcommand{\Lie}{ \operatorname{{\rm Lie}}}

\newcommand{\Hol}{ \operatorname{{\rm Hol}}}
\newcommand{\Ker}{ \operatorname{{\rm Ker}}}

\newcommand{\ED}{ \operatorname{{\rm ED}}}

\newcommand{\Aut}{ \operatorname{{\rm Aut}}}

\newcommand{\GL}{ \operatorname{{\rm GL}}}
\newcommand{\SL}{ \operatorname{{\rm SL}}}

\newcommand{\Lin}{ \operatorname{{\rm Lin}}}

\def\deg{\mathop{\rm deg}}

\def\reg{{\mathop{\rm reg}}}

\def\codim{\mathop{\rm codim}}

\def\lto{\longrightarrow}
\def\hto{\hookrightarrow}

\renewcommand{\epsilon}{\varepsilon}

\def\and{\quad\mbox{and}\quad}

\newcommand{\C}{\ensuremath{\mathbb{C}}}

\newcommand{\R}{\ensuremath{\mathbb{R}}}

\newcommand{\Z}{\ensuremath{\mathbb{Z}}}
\newcommand{\N}{\ensuremath{\mathbb{N}}}
\newcommand{\G}{\ensuremath{\mathbb{G}}}

\newcommand{\bk}{{\ensuremath{\rm \bf k}}}

\newcommand{\hD}{{\hat D}}

\newcommand{\bX}{{\bar X}}
\newcommand{\hX}{{\hat X}}

\newcommand{\hZ}{{\hat Z}}

\newcommand{\tD}{{\tilde D}}

\newcommand{\tG}{{\tilde G}}
\newcommand{\tH}{{\tilde H}}
\newcommand{\tI}{{\tilde I}}

\newcommand{\tV}{{\tilde V}}
\newcommand{\tX}{{\tilde X}}

\newcommand{\brH}{{\breve H}}
\newcommand{\brM}{{\breve M}}
\newcommand{\brY}{{\breve Y}}
\newcommand{\brZ}{{\breve Z}}
\newcommand{\fm}{{\mathfrak m}}

\newcommand{\bF}{{\bar F}}

\newcommand{\bQ}{{\bar Q}}
\newcommand{\bY}{{\bar Y}}
\newcommand{\bZ}{{\bar Z}}

\newcommand{\cW}{{\ensuremath{\mathcal{W}}}}

\def\fm{{\mathfrak m}}

\newcommand{\cL}{{\ensuremath{\mathcal{L}}}}

\newcommand{\cF}{{\ensuremath{\mathcal{F}}}}
\newcommand{\cG}{{\ensuremath{\mathcal{G}}}}
\newcommand{\cS}{{\ensuremath{\mathcal{S}}}}
\newcommand{\cE}{{\ensuremath{\mathcal{E}}}}
\newcommand{\cA}{{\ensuremath{\mathcal{A}}}}

\newcommand{\cO}{{\ensuremath{\mathcal{O}}}}
\newcommand{\cC}{{\ensuremath{\mathcal{C}}}}

\newcommand{\cH}{{\ensuremath{\mathcal{H}}}}

\newcommand{\cP}{{\ensuremath{\mathcal{P}}}}
\newcommand{\cQ}{{\ensuremath{\mathcal{Q}}}}
\newcommand{\cR}{{\ensuremath{\mathcal{R}}}}
\newcommand{\cN}{{\ensuremath{\mathcal{N}}}}
\newcommand{\cT}{{\ensuremath{\mathcal{T}}}}
\newcommand{\cX}{{\ensuremath{\mathcal{X}}}}
\newcommand{\cY}{{\ensuremath{\mathcal{Y}}}}

\newcommand{\cZ}{{\ensuremath{\mathcal{Z}}}}
\newcommand{\p}{\partial}
\newcommand{\de}{\delta}
\newcommand{\id}{{\rm id}}

\renewcommand{\rho}{\varrho}

\newcommand{\brho}{\bar\varrho}

\def\bals#1\eals{\begin{align*}#1\end{align*}}
\def\bal#1\eal{\begin{align}#1\end{align}}

\def\SAut{\mathop{\rm SAut}}

\def\A{{\mathbb A}}

\def\CC{{\mathbb C}}

\def\PP{{\mathbb P}}

\renewcommand{\phi}{\varphi}

\newcommand{\bnum}{\begin{enumerate}}
\newcommand{\enum}{\end{enumerate}}
\renewcommand{\emptyset}{\varnothing}

\newcommand{\brem}{\begin{rem}}
\newcommand{\brems}{\begin{rems}}
\newcommand{\erem}{\end{rem}}
\newcommand{\erems}{\end{rems}}
\newcommand{\bprob}{\begin{prob}}
\newcommand{\eprob}{\end{prob}}
\newcommand{\bprobs}{\begin{probs}}
\newcommand{\eprobs}{\end{probs}}
\newcommand{\bques}{\begin{ques}}
\newcommand{\eques}{\end{ques}}
\newcommand{\bexa}{\begin{exa}}
\newcommand{\bexas}{\begin{exas}}
\newcommand{\eexa}{\end{exa}}
\newcommand{\eexas}{\end{exas}}
\newcommand{\bdefi}{\begin{defi}}
\newcommand{\edefi}{\end{defi}}
\newcommand{\bdefis}{\begin{defis}}
\newcommand{\edefis}{\end{defis}}
\newcommand{\bcor}{\begin{cor}}
\newcommand{\ecor}{\end{cor}}
\newcommand{\blem}{\begin{lem}}
\newcommand{\elem}{\end{lem}}
\newcommand{\bconv}{\begin{conv}}
\newcommand{\econv}{\end{conv}}
\newcommand{\bconj}{\begin{conj}}
\newcommand{\econj}{\end{conj}}
\newcommand{\bprop}{\begin{prop}}
\newcommand{\eprop}{\end{prop}}
\newcommand{\bthm}{\begin{thm}}
\newcommand{\ethm}{\end{thm}}
\newcommand{\bnota}{\begin{nota}}
\newcommand{\enota}{\end{nota}}
\newcommand{\bsit}{\begin{sit}}
\newcommand{\esit}{\end{sit}}
\newcommand{\be}{\begin{equation}}
\newcommand{\ee}{\end{equation}}
\newcommand{\bproof}{\begin{proof}}
\newcommand{\eproof}{\end{proof}}
\def\ba{\begin{array}}
\def\ea{\end{array}}

\begin{document}
\title[Extensions of isomorphisms of subvarieties in flexible varieties]{Extensions of isomorphisms of subvarieties in flexible varieties}

\author{Shulim Kaliman}
\address{
University of Miami, Coral Gables, FL 33124, USA}
\email{kaliman@math.miami.edu}

\date{\today}
\maketitle

\begin{abstract} Let $X$ be an algebraic variety isomorphic  to the complement of a closed subvariety of
dimension at most $n-3$ in $\A^n_\bk$. We find some conditions under which an isomorphism of two closed subvarieties of $X$ can be extended
to an automorphism of $X$. We also study the similar problem for subvarieties of affine quadrics and $\SL(n,\bk)$.
\end{abstract}

\thanks{
{\renewcommand{\thefootnote}{} \footnotetext{ 2010
\textit{Mathematics Subject Classification:}
14R10,\,14R20, 32M17, 32M25.\mbox{\hspace{11pt}}}}

{\footnotesize \tableofcontents}

\section*{Introduction} 

Let $X$ be a smooth 
quasi-affine variety  over an algebraically closed field $\bk$ of characteristic zero
and $Y_1$ and $Y_2$ be closed subvarieties of $X$.
We study the following {\bf extension problem}:

{\em  Under what restrictions on $Y_i$ and $TY_i$ an isomorphism $Y_1\to Y_2$ extends to an automorphism of $X$?}

This question makes sense when $X$ itself possesses a large automorphism
 group $\Aut (X)$ which leads to the notion of a flexible variety \cite{AFKKZ}.
Recall that it is a quasi-affine algebraic variety of  dimension at least 2 on which the group $\SAut (X)$ (generated by elements of all one-parameter unipotent subgroups of $\Aut (X)$) acts $m$-transitively  for every $m>0$.
The simplest  example of a flexible variety is $X=\A^n_\bk$ and the extension problem was studied extensively for such an $X$.

The starting point (and an inspiration)  of that research was, of course, the Abhyankar-Moh-Suzuki theorem \cite{AbMo}, \cite{Su} which
states that given two plane curves isomorphic to a line one can be transferred to the other by an automorphism of $\A^2_\bk$.
Then, disproving an Abhyankar's conjecture,  Jelonek \cite{Je} established 
that if one requires
that $4\dim Y_1 +2 \leq n$ then one gets a positive answer to the extension problem in $\A^n_\bk$
for the case of smooth $Y_i$ \footnote{ Abhyankar's conjecture was also disproved slightly earlier by
 Craighero \cite{Cr} but he did not consider the extension problem in the full generality.}. In the non-smooth case we have to take into consideration $\dim TY_i$
and the more general result established by the author \cite{Ka91} and Vasudevan Srinivas \cite{VS} states the following:

\bthm\label{in.t0} Let $\varphi : Y_1 \to Y_2$ be an isomorphism of two closed subvarieties of $\A^n_\bk$ such that $n \geq \ED (Y_1)+1$ where $\ED (Y_1)=\max (2\dim Y_1+1, \dim TY_1)$.
Then $\varphi$ extends to an automorphism of $\A_\bk^n$.
\ethm

The first result in the case of a flexible variety different from $\A_\bk^n$ is due to  Van Santen (formerly Stampfli) \cite{St} who proved that  given two curves isomorphic to a line
in an algebraic variety $X$ isomorphic to the  special linear group $\SL(n,\bk)$ over $\bk$
one can be transferred to the other by an automorphism of $X$ provided that $n\geq 3$
(in particular, any such  curve is an orbit of a $\G_a$-action). This theorem was generalized later
in the paper of Van Santen and Feller \cite{FS} where they showed that the same is true if one considers two curves isomorphic to a line in a connected linear algebraic group modulo some exceptions. Then Van Santen jointly with J. Blanc showed (among other facts)  that there are closed surfaces isomorphic
to $\A_\bk^2$ in $\SL (2,\bk)$ which cannot be transferred to each other by an 
automorphism of $\SL (2,\bk)$ (as an algebraic variety).

In the present paper beside $\SL (n,\bk)$ we study smooth quadrics in $\A_\bk^n$
and the case of $X$ equal to the complement to a codimension at least 2 subvariety in $\A_\bk^n$
(we call it the Gromov-Winkelmann case since these authors established the flexibility of such an $X$ \cite{Wi}, \cite{Gr1}).  The main results of our paper are the following.

\bthm\label{in.t1} Let $Z$, $Y_1$, and $Y_2$ be closed subvarieties of $\A_\bk^n$  such that $Y_1 \cap Z=Y_2 \cap Z =\emptyset$,  $\dim Z\leq n-3$
and $\ED (Y_1)\leq n-2$. Let $\varphi : Y_1 \to Y_2$ be an isomorphism and $X =\A_\bk^n \setminus Z$.
Suppose also that either

{\rm (a)}  $\dim Z + \dim Y_1 \leq n-3$, or

{\rm (b)} $\dim Y_1 =1$ and  $\dim Z = n-3$.

Then there exists an automorphism $\gamma \in \SAut (X)$ for which $\gamma|_{Y_1}=\varphi$.\footnote{If $\ED (Z) \leq n-1$
then there is no need for this theorem. Indeed, one can consider the isomorphism $\psi : Y_1 \cup Z \to Y_2 \cup Z$ such
that $\psi|_{Y_1} =\varphi$ and $\psi|_Z =\id_Z$. Then Theorem \ref{in.t0} implies that $\psi$ extends to an automorphism
of $\A_\bk^n$ as soon as $\ED (Y_1) \leq n-1$.}
\ethm


\bthm\label{in.t2} 
Let  $m\geq 6$ and $X$ be a hypersurface in $\A_\bk^m$ that is a nonzero fiber of a non-degenerate quadratic form. Suppose that
$\varphi : Y_1\to Y_2$ is an isomorphism of two closed subvarieties of $X$.
 Let   $\ED (Y_i) +\dim Y_i \leq m-2$.
Then $\varphi$ extends to an automorphism of $X$ which belongs to $\SAut (X)$.
\ethm

\bthm\label{in.t3} Let $X=\SL (n,\C)$ and
$\varphi : Y_1\to Y_2$ be an isomorphism of two closed subvarieties of $X$  such that either

{\rm (i)} $\ED (Y_i)+\dim Y_i \leq n-2$, $H_i (Y_1)=0$ for $i \geq 3$ 
and $H_2 (Y_1)$ is a free abelian group; or

{\rm (ii)} $\dim Y_1$ is a curve and $\ED (Y_i)\leq n-2$, or;

{\rm (iii)} $Y_1$ is a once-punctured curve and $\ED (Y_1) \leq 2n-3$. 

Then there exists a holomorphic automorphism $\beta$ of $X$ such that  $\beta|_{Y_1} =\varphi$.
\ethm

\bthm\label{in.t4} Let $\varphi : Y_1\to Y_2$ be an isomorphism of two closed
 subvarieties of $X\simeq \SL (n,\bk )$ with $n \geq 3$ such that $Y_i$ is isomorphic to $\A_\bk^k$. Suppose that 
 either $k\leq \frac{n}{3}-1$ or $k=1$.
Then there exists $ \alpha \in \SAut (X)$ such that $\alpha|_{Y_1} =\varphi$.

\ethm

The paper is organized as follows. In the first six sections we develop some technique which is valid for a wide class of flexible varieties.
More precisely,  in Section 1 we consider a morphism $\kappa : X \to P$ 
of smooth irreducible varieties 
such that $X$ is equipped with an action of a group $G \subset \Aut (X)$ which
preserves every fiber $\kappa^{-1}(p), \, p \in P$ and acts transitively on it.
 We describe some conditions under which one can 
find an algebraic family $\cA\subset G$ of automorphisms of $X$ such that 
for two subvarieties $Y$ and $Z$ of $X$ and a general 
element $\alpha \in \cA$ the varieties $Y$ and $\alpha (Z)$ are transversal (this is a relative version of the transversality theorem from \cite{AFKKZ} that dealt with the case when $P$ was
a singleton). In Section 2 we remind
some facts about flexible varieties and in Section 3 we prove a relative version of a theorem from \cite{AFKKZ} which yields automorphisms of a given flexible
variety with prescribed jets at a finite number of points.  

Sections 4-6 are crucial.
Namely,  the proof of Theorem \ref{in.t0} is heavily based on the fact that for a general linear projection
$\theta : \A_\bk^n \to \A_\bk^{n-1}$ the variety $\theta (Y_i)$ is closed in $\A_\bk^{n-1}$
and the restriction of $\theta$ yields an isomorphism $Y_i \to \theta (Y_i)$.  Note that
$\theta$ can be viewed as the composition of a fixed projection $\theta_0$ and a general linear automorphism
of $\A_\bk^n$.  
Hence in Section 4 we imitate this idea for  a morphism $\rho : X \to Q$ over $P$ as above
but with each fiber of $\kappa$ being a $G$-flexible variety. We show
that when the morphism $\rho$ and the variety $Q$ are smooth and
 $Z$ is a closed subvariety of $X$ with $\ED (Z) \leq \dim Q$ 
 then for a general element $\alpha$ of some algebraic family $\cA \subset G$ the morphism
$\rho|_{\alpha (Z)} : \alpha (Z) \to Q$ is an injection and, furthermore, the induced morphism
$T\alpha (Z) \to TQ$ of the Zariski tangent bundles is also an injection.
However,  a priori this map is not proper and in Section 5 we describe some conditions
under which the morphism $\rho|_{\alpha (Z)}$ is also proper and, therefore, 
$\rho\circ \alpha (Z)$ is closed in $Q$ and $\rho|_{\alpha (Z)} : \alpha (Z) \to \rho\circ \alpha (Z)$
is an isomorphism. These facts are 
already sufficient for the proof of Theorem \ref{in.t1} in Section 7.  Section 6 is devoted
to the independently interesting case when $\rho$
is a partial quotient morphism of a $\G_a$-action on $X$. In this situation we cannot
claim that $\rho$ and $Q$ are smooth and cannot guarantee the properness of $\rho|_{\alpha (Z)}$.
However, we establish that for any given finite subset $S \subset Z$ and a general $\alpha \in \cA$ one can find
a neighborhood $V'$ of $\rho (\alpha (S))$ in $\overline{ \rho (\alpha (Z))}\subset Q$ such that for
$V=\rho^{-1} (V') \cap \alpha (Z)$ the restriction of $\rho|_{\alpha (Z)}$ yields an isomorphism $V \to V'$.
Applications of this result go beyond the present paper (e.g., see \cite{KaKuTr}).

Section 8 contains Theorem \ref{in.t2} and in Sections 9 and 10 we prove some technical results
that enable us to obtain Theorems \ref{in.t3} and \ref{in.t4} in Section 11
where we also present an example of a topological
obstruction for the extension problem in the case of general flexible varieties.

{\em Acknowledgements.} The author is deeply indebted to his referee for catching mistakes in the original version of this manuscript
and for an unusually thorough review which was a great help to the author.

\section{Algebraically generated groups of automorphisms}

Let $X$ be an irreducible algebraic variety and $\Aut (X)$ be the group of its algebraic automorphisms. 
Recall the following terminology introduced by Ramanujam \cite{Ra1}.

\bdefi\label{agga.d1}  (1) Given an irreducible algebraic variety $\cA$ and
a map $\varphi:\cA\to\Aut(X)$ we say that $(\cA,\phi)$
is an {\em algebraic family of automorphisms on $X$} if the induced map
$\cA\times X\to X$, $(\alpha,x)\mapsto \varphi(\alpha).x$ is a morphism.

(2) In the case when $\cA$ is a connected algebraic group and the induced map 
$\cA\times X\to X$ is not only a morphism but also an action of $\cA$ on $X$ we call this family a connected algebraic subgroup of $\Aut (X)$.
\edefi

\bdefi\label{agga.d1r}
Following \cite[Definition 1.1]{AFKKZ} we call a subgroup $G$ of $\Aut (X)$ algebraically generated if it is generated as an abstract group by a family 
$\cG$ of connected algebraic subgroups of $\Aut (X)$.
\edefi



We have the following important fact \cite[Theorem 1.15]{AFKKZ} (which is the analogue of  the Kleiman transversality theorem \cite{Kl}
for algebraically generated groups).

\bthm\label{agga.t1} ({\rm Transversality Theorem}) Let a subgroup $G\subseteq \Aut(X)$ be
algebraically generated by a system $\cG$ of connected algebraic
subgroups closed under conjugation in $G$. Suppose that $G$ acts
with an open orbit $O\subseteq X$.

Then there exist subgroups $H_1,\ldots, H_m\in \cG$ such that for
any locally closed reduced subschemes $Y$ and $Z$ in $O$ one can
find a Zariski dense open subset $U=U(Y,Z)\subseteq H_1\times
\ldots \times H_m$ such that every element $(h_1,\ldots, h_m)\in
U$ satisfies the following: \bnum[(a)]
\item  {\em The translate $(h_1\cdot\ldots\cdot h_m).Z_\reg$
meets $Y_\reg$
transversally. } 
\item $\dim (Y\cap (h_1\cdot\ldots\cdot h_m).Z)\le
\dim Y+\dim Z-\dim X$. \footnote{We put the dimension of empty sets equal to $-\infty$.}\\
In particular $Y\cap (h_1\cdot\ldots\cdot h_m).Z=\emptyset$
if $\dim Y+\dim Z<\dim X$.
\enum
\ethm

We need to generalize \cite[Theorem 1.15]{AFKKZ} further.

\bthm\label{agga.t2} {\rm (Collective Transversality Theorem)}  Let $X$ and $P$ be smooth irreducible
algebraic varieties and $\kappa : X \to P$ be a smooth morphism  (in particular $X\times_PX$ is smooth
and $\dim X\times_PX =2 \dim X - \dim P$).
Let a group $G\subseteq \Aut(X)$ be
algebraically generated by a system $\cG$ of connected algebraic
subgroups closed under conjugation in $G$. 
Suppose that the $G$-action transforms every fiber $\kappa^{-1} (p)$ into
itself and, furthermore, the restriction of the $G$-action to $\kappa^{-1} (p)$ is transitive for every $ p\in P$.

Then there exist subgroups $H_1,\ldots, H_m\in \cG$  such that for
any locally closed reduced subschemes $Y$ and $Z$ in $X$ one can
find a Zariski dense open subset $U=U(Y,Z)\subseteq H_1\times
\ldots \times H_m$ so that every element $(h_1,\ldots, h_m)\in
U$ satisfies the following:  

{\rm (i)} 
 $\dim (Y\cap (h_1\cdot\ldots\cdot h_m).Z)\le \dim (Y\times_PZ) + \dim P -\dim X$. 
 
 In particular, when  $\dim Y\times_PZ \leq \dim Y+ \dim Z -\dim P$ one has 

{\rm (ii)} $\dim (Y\cap (h_1\cdot\ldots\cdot h_m).Z)\le \dim Y+\dim Z -\dim X$. 

Furthermore, suppose that the inequality $\dim Y\times_PZ \leq \dim Y+ \dim Z -\dim P$ holds
and also that $Z$, $Y \times_PZ$, and  $Y\times_PX$ are smooth.
Then

{\rm (iii)}  $(h_1\cdot\ldots\cdot h_m).Z$ meets $Y$ transversally.

\ethm

The proof of Theorem \ref{agga.t2} is an adjustment of the proof of  \cite[Theorem 1.15]{AFKKZ}. 
Hence following the latter we establish first three facts which in the case of a singleton $P$ are nothing but Propositions 1.5, 1.8, and 1.16 in  \cite{AFKKZ}.

\bprop\label{1.2} {\rm (Analogue of \cite[Proposition 1.5]{AFKKZ})}  Let the assumptions of Theorem \ref{agga.t2} hold.
There are (not necessarily distinct) subgroups
$H_1,\ldots, H_m\in \cG$ such that for every $p \in P$ and each $x \in \kappa^{-1} (p)$ one has
\be\label{1.2.a}  \kappa^{-1} (p)=
H_1.(H_2. \cdots .(H_m.x).
\ee
\eprop

\bproof
Let us
introduce the partial order on the set of sequences in $\cG$ such that for $\cH=(H_1,\ldots, H_m)$ and $\cH'= (H'_1,\ldots H'_s)$ one has
$$
\cH\succcurlyeq \cH'
\Longleftrightarrow
\exists i_1<\ldots <i_s:\quad (H_1',\ldots , H_s')=
(H_{i_1},\ldots, H_{i_s})\,.
$$

Assuming first that $\kappa$ has a section $\lambda : P \to X$ with $S=\lambda (P)$
we consider $\cH. S=\bigcup_{x\in S}\cH.x$ where $\cH.x= H_1.(H_2. \cdots .(H_m.x)$. Then
such a set $\cH.S$ is constructible (since it is the image of the algebraic variety $H_1\times \ldots \times H_m \times S$ under a morphism). In particular, $X_\cH :=X \setminus \cH. S$ is a constructible set.
Furthermore, the following property holds:

\begin{center} if $\cH  \succcurlyeq \cH'$, then $\cH'.S \subset \cH. S$ and, therefore, $X_\cH \subset X_{\cH'}$.\end{center}

Because of transitivity for every $y \in \kappa^{-1}(p)$ we can find a sequence $\cH=(H_1,\ldots, H_m)$  and $g=h_1\cdot\ldots\cdot h_m$ (where $h_i \in H_i$) for which
$y=g.x$ where $x = S\cap \kappa^{-1}(p)$.  Hence $y \in \tilde \cH. S$ for any sequence $\tilde \cH$ of the form  $\tilde \cH=( \cH_1 , \cH)$.  In particular, choosing $y$ in any  given irreducible component 
$C$ of $X_{\cH_1}$ we guarantee $C$  is not contained in $X_{\tilde \cH}$. That is, $\dim X_{\tilde \cH}\cap C < \dim C$ since $X_{\tilde \cH}\subset X_{\cH_1}$.

Thus, enlarging $ \cH$ we can reduce the dimension of   $X_{ \cH}$ and continuing this process we can make $X_{ \cH}=\emptyset$. In particular, for every
$p \in P$ and  $x = S\cap \kappa^{-1}(p)$ one has $ \cH .x =\kappa^{-1}(p)$.  Therefore, for every $y \in \kappa^{-1}(p)$ there exists $g=h_1\cdot\ldots\cdot h_m$ as before (with $\cH$ now
being independent from $y$) for which $y=g.x$. Let $\cH^t =(H_m, H_{m-1}, \ldots , H_1)$ and $\tilde \cH =(\cH, \cH^t)$, i.e., $x \in \cH^t.y$.
Then one has  $\tilde \cH .y =\kappa^{-1}(p)$ for every $y \in \kappa^{-1}(p)$, i.e we get the desired
conclusion in the presence of a section $\lambda$.

In the general case consider an \'etale neighborhood $W$ of a point $p \in P$ and suppose that $V\subset P$ is the image of $W$ under the natural morphism. Since $\kappa$ is smooth one can suppose that
for an appropriate choice of $W$  the natural projection $\tau: X\times_P W \to W$ has a section. Consider the induced $G$-action on $X\times_P W$.
By the previous argument there is a sequence $\cH$ such that for every $w\in W$ and $z \in \tau^{-1}(w)$ one has $\cH.z=\tau^{-1}(w)$. Applying the natural projection
$ X\times_P W \to X$ we see that for every $p \in V$ and every $y \in \kappa^{-1}(p)$ we have $\cH .y =\kappa^{-1}(p)$. Choosing now a finite number of \'etale neighborhoods that cover $X$,
we can enlarge $\cH$ so that it works for each of these neighborhoods. This implies the desired conclusion. 

\eproof

\bprop\label{1.5} {\rm (Analogue of \cite[Proposition 1.8]{AFKKZ})}   Let the assumptions of Theorem \ref{agga.t2} hold. 
Assume that the generating family $\cG$ of
connected algebraic subgroups is closed under conjugation in $G$,
i.e., $gHg^{-1}\in \cG$ for all $g\in G$ and $H\in \cG$. Then
there is a sequence $\cH=(H_1,\ldots, H_m) $ in  $\cG$ such that
for all $p \in P$ and $x\in \kappa^{-1}(p)$ the tangent space $T_x \kappa^{-1}(p)$
is spanned by the tangent spaces
$$
T_x(H_1.x),\ldots ,T_x(H_m.x)
$$
to the orbits $H_1.x,\ldots ,H_m.x$ at $x$.
\eprop

\bproof
Let $\cH=(H_1,\ldots, H_m) $  be a sequence in $\cG$ satisfying the conclusion of  Proposition \ref{1.2}. Consider the map
$$\Phi_\cH:H_1\times\ldots \times H_m\times X\to X\times_P X \, \, {\rm given \, \, by} \, \, ((h_1,\ldots, h_m), x)\to ((h_1\cdot \ldots \cdot h_m).x,x).$$ 
The fiber $\tau^{-1}(x)$ of the second projection $\tau : X\times_PX \to X$ over $x \in X$ is naturally isomorphic to $\kappa^{-1}(p)$ where $p=\kappa (x)$. 
Hence, by Proposition \ref{1.2} $\Phi_\cH$ is a surjective map while the assumptions of Theorem \ref{agga.t2} imply that $\tau$ is a smooth morphism.
Let us consider the map of relative tangent bundles
$$
d\Phi_\cH: T(H_1\times \ldots \times H_m\times X/X)\to
\Phi_\cH^*(T((X\times_PX)/X))\,
$$
and its restriction to $\{ (e,\ldots,e)\}\times X\cong X$ (where $e$ is the identity in the group $G$),
$$
d\Phi_\cH: T_eH_1\times \ldots \times T_eH_m
\times X\to \Phi_\cH^*(T((X\times_PX)/X))\,.
$$
The set $U_\cH$ of points in $X$ where this map
is surjective is, of course, open. By \cite[Proposition 1.8]{AFKKZ} for every $p \in P$ and $x\in \kappa^{-1}(p)$ the tangent space $T_x \kappa^{-1}(p)$ is spanned by
the tangent spaces $T_x(H.x)$, where $H\in \cG$. 
Hence $\bigcup_\cH U_\cH$ coincides with $X$. Since an increasing union of open subsets stabilizes,
we obtain that $X=U_\cH$ for $\cH$ sufficiently large (with respect to the partial order introduced in Proposition \ref{1.2}). This yields the desired conclusion.

\eproof

\bprop\label{1.16}{\rm (Analogue of \cite[Proposition 1.16]{AFKKZ})}
 Let the assumption of Theorem \ref{agga.t2} hold.
Then there is a sequence $H_1,\ldots, H_m$ in $\cG$
so that for a suitable open dense subset $U\subseteq
H_{m}\times\ldots \times H_{1}$, the map
\be\label{mapop}
\Phi_{m}: H_m\times \ldots\times H_1\times X \lto X\times_P X
\quad\mbox{with} \quad (h_m,\ldots,h_1,x)\mapsto
((h_m\cdot\ldots\cdot h_1).x ,x)
\ee
is surjective and smooth on $U\times X$.

\eprop

\bproof

By Proposition \ref{1.2}  there are subgroups
$H_1,\ldots, H_m\subseteq G$  in $\cG$ such that
$\Phi_{m}$ is surjective.  Let $U_m\subset H_m\times \ldots \times H_1\times X$ be the set of points
where $\Phi_m$ is smooth. Then $U_m$ is non-empty by  \cite[Chapter III, Corollary 10.7]{Har} and it is open
by [SGA1, Exp. II, Prop. 1.1].      
Consider the
complement $A_m=(H_m\times \ldots\times H_1\times X)\backslash
U_m$. Let us study the effect of increasing the number $m$ of factors in the product $H_m\times \ldots\times H_1$.

 Suppose that $H$ is an element
of $\cG$ and $\Phi'$  plays the same role for the sequence $H_1,\ldots, H_m, H$ as $\Phi_m$ for  the sequence $H_1,\ldots, H_m$.
Note that for a fixed $h \in H$ the restriction of $\Phi'$ yields the morphism $\{ h \} \times H_m\times \ldots\times H_1\times X \lto X\times_P X$
that  is the composition of $\Phi_m$ and the automorphism $\varphi$ of
$X \times_P X$ given by $(x,y)\to (h.x, y)$ which implies smoothness of $\Phi'$ on $H\times U_m$.      
Thus, $U_{m+1}\supseteq H_{m+1}\times U_m$ and
$A_{m+1}\subseteq H_{m+1}\times A_m$. Increasing the
number of factors by $H_{m+1},\ldots, H_{m+k}$ in a suitable way,
we can achieve that
\be\label{auxia} \dim A_{m+k}< \dim
(H_{m+k}\times\ldots\times H_{m+1}\times A_m) \,. \ee
Indeed, if
$(h_m,\ldots,h_1,x)\in A_m$ and $y=(h_m\cdot\ldots\cdot h_1).x$
then by  Proposition
\ref{1.5} for suitable $H_{m+k},\ldots, H_{m+1}$ the map
$$
H_{m+k}\times\ldots\times H_{m+1}\times X\lto X\times_P X, \, (h_{m+k}, \ldots, h_{m+1},z) \to ((h_{m+k} \cdot \ldots \cdot h_{m+1}). z,z)
$$
is smooth in all points $(e,\ldots, e, y)$ where $e$ is the identity of the group $G$. In particular, $\Phi_{m+k}$ is smooth in all points
$(e,\ldots, e, h_m,\ldots, h_1,x)$ with $x\in X$, i.e.,
$$
\{ (e,\ldots, e)\} \times A_m\cap  A_{m+k}=\emptyset.
$$
Now \eqref{auxia} follows.

Thus increasing the number of factors suitably we can achieve
that
$$\dim A_m<\dim(H_m\times\ldots\times H_1)\,.$$
That is, the image of $A_m$ under the projection
$$
\pi: H_{m}\times\ldots\times H_{1}\times X\lto
H_{m}\times\ldots \times H_{1}
$$
is nowhere dense. Hence there is an open dense subset
$U\subseteq H_{m}\times\ldots \times H_{1}$
such that $\Phi_m: U\times X\to X\times_P X$ is smooth.

\eproof
 
\brem\label{agga.r1r} Let $\Phi_m$ and $H_1, \ldots, H_m$ be as in Proposition \ref{1.16} and let 
$H$ be an element $\cG$. Suppose that $\Phi'$ (resp. $\Phi''$)  plays the same role for the sequence $H_1,\ldots, H_m, H$ (resp. $H, H_1,\ldots, H_m$) as $\Phi_m$  
for  the sequence $H_1,\ldots, H_m$. We showed already in the proof of Proposition \ref{1.16} that
that $\Phi'$ is smooth on $H\times U_m$. Similarly, $\Phi''$ is smooth on the preimage $U''$ of $U_m$ under the natural 
projection $H_m \times \ldots \times H_1 \times H \times X \to H_m \times \ldots \times H_1 \times X$.
Indeed, the smoothness of $\Phi_m|_{U_m}$ is equivalent to the fact that
it factors locally through an \'etale morphism $U_m \to \A_\bk^n \times (X \times_PX)$ over $X \times_PX$ or, in other words, that $\Phi_m|_{U_m}$ admits 
a local \'etale section $s : X\times_PX\to U_m$
(where $s(x,y)= (\hat s (x,y), y)$ with $\hat s(x,y) \in H_m\times \ldots\times H_1$). Denoting by  $\varphi$ the automorphism of
$X \times_P X$ given by $(x,y)\to (h.x, y)$ one observes that
$(x,y) \to (\hat s \circ \varphi (x,y),h, y)\in  H_m\times \ldots\times H_1 \times h \times X$
is a local section of $\Phi''|_{U''}$, i.e., $\Phi''$ is smooth on $U''$.

\erem

\bproof[Proof of Theorem \ref{agga.t2}] 
By  Proposition \ref{1.16}  there are subgroups $H_1,\ldots, H_m$ in $\cG$ such
that $\Phi_m:U \times X\to X\times_P X$ is smooth for some open
subset $U\subseteq H_m\times\ldots \times H_1$. 
Let $\cC=\Phi_m^{*}(Y\times_P Z)\cap (U\times X)$.

Consider first the case in (i) when $Y\times_P Z$ is smooth.
Then $\cC$ is smooth.
By \cite[Chapter III, Corollary 10.7]{Har} the general fibers of the
projection $\pi_\cC : \cC\to U$  
are smooth as well. 
Suppose that $\pi_\cC$ is dominant (otherwise the general fibers of $\pi_\cC$ are empty). Shrinking $U$ we may now assume that all fibers of $\pi_\cC $  are 
smooth.  Then the dimension of every fiber of $\Phi_m$ is $\dim U -\dim X + \dim P$.
Thus $\dim \cC = \dim U + \dim P -\dim X + \dim Y\times_P Z$ and the dimension of every fiber $\pi_\cC^{*}(h)$  of $\pi_\cC$ is 
\be\label{agga.eq3} \dim Y\times_P Z +\dim P - \dim X . \ee
Observe now that  for a point $h=(h_m,\ldots, h_1)\in U$ the fiber
$\pi_\cC^{*}(h)$ maps bijectively via $\{ h \} \times X \to X, \, (h,x) \to (h_m \cdot \ldots \cdot h_1).x$     onto $Y\cap
(h_1\cdot\ldots\cdot h_m).Z$ which yields (i) (and therefore (ii)) in the case of smooth $Y\times_P Z$.




In the general case stratifying $Z$ and $Y$ we can find Zariski dense open subsets $Z_0 \subset Z$ and $Y_0\subset Y$ such that
$Y_0\times_P Z_0$ is smooth.  
By Formula (\ref{agga.eq3}) we see that  for a general $h \in U$ the dimension
of $Y_0 \cap (h_1\cdot\ldots\cdot h_m).Z_0$ is at most $\dim Y\times_PZ +\dim P - \dim X$. Let $Z'=Z\setminus Z_0$ and $Y'=Y \setminus Y_0$ and 
consider, say, the pair $(Y_0, Z')$. We can suppose that  $Y_0\times_P Z'$ is smooth (otherwise stratify further). Then the same argument with Formula (\ref{agga.eq3})
implies that  the dimension of $Y_0 \cap (h_1\cdot\ldots\cdot h_m).Z'$ is at most $\dim Y_0\times_P Z' +\dim P - \dim X
\leq \dim Y\times_P Z +\dim P - \dim X$. Repeating this procedure for the pairs
$(Y',Z_0)$ and $(Y',Z')$ we get (i) and (ii) in full generality.

For (iii) consider  $\cZ=U\times Z$ and $\cY=\Phi_m^{*}(Y\times_P X)\cap (U\times X)$, i.e., $\cC$ (as a scheme) is the
intersection of $\cZ$ and $\cY$.  As before, shrinking $U$ we can suppose that all fibers of the
natural projections $\pi_\cY : \cY\to U$ and  $\pi_\cZ : \cZ \to U$   
are smooth. 
 Observe also that
 $\pi_\cZ^{*}(h)= \{ h\} \times Z\subset \{ h\} \times X$ and  $\pi_\cY^{*}(h)=\{ h\} \times (h_1^{-1} \cdot \ldots \cdot h_m^{-1}).Y\subset \{ h\} \times X$ \footnote{In particular,
the last equality implies that $Y$ is smooth under the assumption of (ii).}. It remains to note that if these two smooth
subvarieties of  $\{ h \}\times X$ do not meet transversely and the dimension of  their intersection $\pi_\cC^{*}(h)$ is $\dim Y +\dim Z - \dim X$ then $\pi_\cC^{*}(h)$ cannot be smooth (as a scheme).
Indeed, the absence of transversality implies that at some closed point of $x \in \pi_\cC^* (h)$  
the dimension of the intersection of the tangent spaces $T_x\pi_\cY^{*}(h)$ and $T_x\pi_\cZ^{*}(h)$  
is greater than $\dim Y +\dim Z - \dim X$. Note that this intersection coincides with the intersection of
kernels of the differentials of all functions from the defining ideals of $\pi_\cZ^{*}(h)$ and $\pi_\cY^{*}(h)$,
and, therefore, from the defining ideal of $\pi_\cC^{*}(h)$. However, for a smooth $\pi_\cC^{*}(h)$
this dimension must be equal to  $\dim Y +\dim Z - \dim X$.  
Hence the smoothness of $\pi_\cC^{*}(h)$
established before yields (iii) which concludes the proof.

\eproof

\brem\label{agga.r1}
(1) Suppose that $\kappa (Y)$ is dense in $P$ and all fibers of $\kappa|_Y$ are of the same dimension (say, $\kappa|_Y$ is flat). 
That is, the dimension of each of these fibers is $\dim Y -\dim P$. Then  $\dim Y\times_{P} Z = \dim Y + \dim Z- \dim P$
and we are under the assumption of (ii) in Theorem \ref{agga.t2}, i.e., the dimension of $Y\cap (h_1\cdot\ldots\cdot h_m).Z$
is at most $\dim Y +\dim Z - \dim X$.


(2) Let us emphasize the following fact the first part of which follows from the argument in the proof of Theorem \ref{agga.t2}.
\erem

\bprop\label{agga.p2}

If a sequence $H_1,\ldots, H_m \in \cG$ satisfies Proposition \ref{1.16}  then it satisfies also
Theorem \ref{agga.t2}. Furthermore, for any element $H$ of $\cG$
the sequence  $H_1,\ldots, H_m, H$ (resp. $H, H_1,\ldots, H_m$) satisfies Theorem \ref{agga.t2} as well.
\eprop

\bproof The second statement is true because by Remark \ref{agga.r1r}
the sequence  $H_1,\ldots, H_m, H$ (resp. $H, H_1,\ldots, H_m$) satisfies Proposition \ref{1.16}.

\eproof

\brem\label{agga.r2} 
(1) Let us consider an application of Proposition \ref{agga.p2}.
Suppose that $X$ and $P$ are smooth irreducible
algebraic varieties and $\kappa : X \to P$ be a morphism  which is not in general smooth or even dominant.
Let a group $G\subseteq \Aut(X)$ be
algebraically generated by a system $\cG$ of connected algebraic
subgroups closed under conjugation in $G$. 
Suppose that the $G$-action transforms every fiber $\kappa^{-1} (p)$ into
itself and, furthermore, the restriction of the $G$-action to $\kappa^{-1} (p)$ is transitive for every $ p\in P$.

By the Generic Smoothness theorem (e.g., see \cite[Chapter III, Corollary 10.7]{Har}) we can present $\kappa (X)$ as a disjoint union $\bigcup_{k=1}^nP_k$ of
smooth varieties such that for $X_k=\kappa^{-1}(P_k)$ the morphism $\kappa|_{X_k} : X_k \to P_k$ is smooth.
For $Y$ and $Z$ as in Theorem \ref{agga.t2} let $Y_k=X_k \cap Y$ and $Z_k=X_k\cap Z$.
Then by Theorem \ref{agga.t2}  there exist subgroups $H_1^k,\ldots, H_{m_k}^k\in \cG$  such that 
one can find a Zariski dense open subset $U_k=U(Y_k,Z_k)\subseteq H_1^k\times
\ldots \times H_{m_k}^k$ so that  for every element $(h_1^k,\ldots, h_{m_k}^k)\in
U_k$ we have the inequality
$$\dim (Y_k\cap (h_1^k\cdot\ldots\cdot h_{m_k}^k).Z_k)\le \dim (Y_k\times_{P_k}Z_k) + \dim P_k -\dim X_k \, .$$
Consider now a general element $\alpha$ in
$$H_1^1\times
\ldots \times H_{m_1}^1 \times H_1^2 \times \ldots \times   H_1^n\times \ldots \times H_{m_n}^n.$$
Then Proposition \ref{agga.p2} implies that one has
$$\dim (Y\cap \alpha .Z)\le \max \{ \dim (Y_k\times_{P_k}Z_k) + \dim P_k -\dim X_k| k=1, \ldots n\}  \, .$$

(2) Let us consider now the case when in (1)
all nonempty fibers of the morphism $\kappa|_Y : Y \to P$ are of the same dimension.
By Remark \ref{agga.r1} (1) we have  $\dim (Y_k\times_{P_k}Z_k) =  \dim Y_k + \dim Z_k - \dim P_k$ for every $k=1, \ldots , n$.
Hence 
$$\dim (Y\cap \alpha .Z)\le \max \{ \dim Y_k + \dim Z_k -\dim X_k| k=1, \ldots n\} \, .$$
Furthermore,  if all nonempty fibers of $\kappa|_X : X \to P$ are of the same dimension then 
$\dim Y_k -\dim X_k \leq \dim Y -\dim X$ and we have the following.
\erem

\bprop\label{agga.p3}  Let $X$ and $P$ be smooth irreducible
algebraic varieties and $\kappa : X \to P$ be a flat morphism.
Let a group $G\subseteq \Aut(X)$ be
algebraically generated by a system $\cG$ of connected algebraic subgroups closed under conjugation in $G$. 
Suppose that the $G$-action transforms every fiber $\kappa^{-1} (p)$ into
itself and, furthermore, the restriction of the $G$-action to $\kappa^{-1} (p)$ is transitive for every $ p\in P$.
Let $Y$ and $Z$ be locally closed reduced subschemes of $X$ such that 
all nonempty fibers of the morphis  $\kappa|_Y : Y \to P$ are of the same dimension.
Then there exist subgroups $H_1,\ldots, H_{m}\in \cG$  such that such that for any $Y$ and $Z$ as above
one can find a Zariski dense open subset $U=U(Y,Z)\subseteq H_1\times
\ldots \times H_{m}$ so that for every element $(h_1,\ldots, h_{m})\in
U_k$ we have the inequality
$\dim (Y\cap (h_1\cdot\ldots\cdot h_m).Z)\le \dim Y+\dim Z -\dim X$. 

\eprop
}

\section{Flexible varieties}

\bdefi\label{fm.d1}  (1) A derivation $\sigma$ on the ring $A$ of regular functions on a quasi-affine algebraic variety $X$ is called locally nilpotent
if for every $0\ne a \in A$ there exists a natural $n$ for which $\sigma^n (a)=0$. For the smallest $n$ with this property one defines
the degree of $a$ with respect to $\sigma$ as $\deg_\sigma a=n-1$. This derivation can be viewed as a vector field on $X$ which
we also call locally nilpotent. The  flow of this vector field is an algebraic $\G_a$-action on $X$, i.e., the action of the group $(\bk, +)$ 
which can be viewed as a one-parameter unipotent group $U$ in the group $\Aut (X)$ of all algebraic automorphisms of $X$.
In fact, every $\G_a$-action is generated by a locally nilpotent vector field (e.g, see \cite{Fre}).

(2) A smooth quasi-affine algebraic variety $X$ of dimension at least 2
is called flexible if for every $x \in X$ the tangent space $T_xX$ is spanned by the tangent vectors
to the orbits of one-parameter unipotent subgroups of $\Aut (X)$ through $x$.

(3) The subgroup $\SAut (X)$ of $\Aut X$ generated by all one-parameter unipotent subgroups is called special.
\edefi

We have the following \cite[Theorem 01]{AFKKZ} and \cite[Theorem 2.12]{FKZ-GW}.

\bthm\label{fm.t1} For every irreducible smooth quasi-affine algebraic variety $X$ of dimension at least 2 the following are equivalent

{\rm (i)} the special subgroup $\SAut (X)$ acts transitively on $X$;

{\rm (ii)}  the  special subgroup $\SAut (X)$ acts infinitely transitively on $X$ (i.e., for every natural $m$
the action is $m$-transitive \footnote{Recall that a group $G$ acts $m$-transitively on a space $Y$ if for any two $m$-tuples $(y_1, \ldots , y_m)$ and $(y_1', \ldots , y_m')$
of distinct points in $Y$ there is an element $\alpha \in G$ such that $\alpha (y_i) =y_i'$ for every $i=1, \ldots , m$.});

{\rm (iii)} $X$ is flexible.

\ethm

\bdefi\label{fm.d3} (1) For every locally nilpotent vector fields $\sigma$ and each function $f \in \Ker \sigma$ from its kernel the field
$f\sigma$ is called a replica of $\sigma$. Recall that such a replica is automatically locally nilpotent.

(2)  Let $\cN$ be a set of complete algebraic vector fields on $X$. We say that a subgroup $G_\cN \subset  \Aut (X)$ is generated by $\cN$ if $G_\cN$ is generated 
by the elements of all one-parameter  groups that are flows of  complete vector fields from $\cN$.

(3) A collection of locally nilpotent vector fields $\cN$ is called saturated if $\cN$ is closed under conjugation by elements in $G_\cN$
and for every $\sigma \in \cN$ each replica of $\sigma$ is
also contained in $\cN$.

\edefi

\bdefi\label{fm.d3a} Let $\cN$ be a saturated set of locally nilpotent vector fields on $X$ and  $G=G_\cN$ be a subgroup of $\SAut (X)$
which is generated by  $\cN$.
Then $X$ is called $G$-flexible if for any $x \in X$ the vector space $T_xX$ is generated by  the values of locally nilpotent vector fields from $\cN$ at $x$.
\edefi

\brem\label{fm.r1r} A priori the notion of $G$-flexibility depends on $\cN$ while $\cN$ is not determined uniquely by the group $G=G_\cN$.
However, the following generalization of Theorem \ref{fm.t1} ( see \cite[Theorem 2.12]{FKZ-GW}) shows that
this notion depends on $G$ only.

\erem

\bthm\label{fm.t1a}
Let $X$ be an irreducible smooth quasi-affine algebraic variety of dimension at least 2
and $G$ be a subgroup generated by a saturated set $\cN$ of locally nilpotent vector fields on $X$. 
Then the following are equivalent

{\rm (i)} $G$ acts transitively on $X$;

{\rm (ii)}  $G$ acts infinitely transitively on $X$;

{\rm (iii)} $X$ is $G$-flexible.

\ethm

The next fact is a straightforward consequence of  Theorem \ref{fm.t1a}.

\bprop\label{fm.p1} Let $X$ be a $G$-flexible variety where $G$ is as in Definition \ref{fm.d3a} and $\Delta$ be the diagonal in $X\times X$.
Then the natural action of $G$  on the variety $X\times X\setminus \Delta$ is transitive.\footnote{Precaution:
this natural action of $G$ on $X\times X\setminus \Delta$ is not infinitely transitive.}
\eprop

The following result will be very useful later in this paper.

\bthm\label{fm.t2} {\rm (\cite[Theorem 4.14 and Remark 4.16]{AFKKZ})} Let $x_1, \ldots , x_{m}$ be distinct points in a $G$-flexible variety $X$ of $\dim X=n$ where $G$
is generated by a saturated set $\cN$ of locally nilpotent vector fields on $X$.
Then there exists an automorphism $\alpha \in G \subset \SAut (X)$ such that it fixes the
 points $x_1, \ldots , x_{m}$ and for every $i$ the linear map ${\rm d}\alpha |_{T_{x_i}X}$ coincides
with a prescribed element $\beta_i$ of $\SL (n,\bk)$. 
Furthermore,  let $k \in \N$ and $\gamma_i$ be a $k$-jet of an isomorphism
between two \'etale neighborhoods of $x_i$ in $X$
which preserves  $x_i$ and a local volume form at $x_i$.  Then $\alpha$ can be chosen so that for every $i=1, \ldots  , m$  the $k$-jet  of
$\alpha$ at $x_i$ coincides with $\gamma_i$.

\ethm

By the Rosenlicht Theorem (e.g., see \cite[Theorem
2.3]{PV}) for $X$, $A$, and $U$ as in Definition \ref{fm.d1} (1)
one can find a finite set of $U$-invariant functions $a_1,\ldots,a_m\in A$,
which separate general $U$-orbits in $X$. They generate a morphism $\rho : X \to Q$ into an affine
algebraic variety $Q$ (in particular, $\dim Q = \dim X -1$ because general $U$-orbits are one-dimensional). Note that this set of invariant functions can be chosen so that $Q$ is normal (since $X$ is normal).

\bdefi\label{fm.d3} Such a morphism $\rho : X \to Q$ into a normal $Q$ is called a partial quotient.
In the case when $a_1, \ldots , a_m$ generate the subring $A^U$ of $U$-invariant elements of $A$
such a morphism is called the categorical quotient.\footnote{However, in general $A^U$ is not finitely generated by the 
Nagata's example. That is, why, following \cite{FKZ-GW} we prefer to formulate some results for partial quotients.}
\edefi

\bprop\label{fm.p2a} Let $G\subset \SAut (X)$ be generated by a saturated set $\cN$ of locally nilpotent vector fields on a 
smooth quasi-affine  algebraic variety $X$ (i.e., $X$ is $G$-flexible).
Suppose that $\cN_0$ is a nontrivial subset of $\cN$ which is also saturated and closed under conjugation by elements of $G$.
Then $X$ is $G_0$-flexible where the group $G_0$ is generated by  $\cN_0$.
\eprop

\bproof By Theorem \ref{fm.t2} for any point $x\in X$ where $\sigma\in \cN_0$ does not vanish we can find $\alpha_1, \ldots , \alpha_n \in G$ where $n=\dim X$ such that $\alpha_i (x)=x$ and the values
of the fields $\alpha_{1*}(\sigma), \ldots , \alpha_{n*} (\sigma)$ at $x$ generate $T_xX$.
 That is, $T_xX$ is generated by the values of the fields from $\cN_0$.
 Since $G$ acts transitively on $X$ we can guarantee the same for every $x \in X$
 using conjugation by elements in $G$. Thus,
$X$ is $G_0$-flexible which is the desired conclusion. 

\eproof

\bdefi\label{adp.d2} Let $\delta_1$ and $\delta_2$ be a pair of locally nilpotent vector fields on $X$.
Suppose that $\Ker \delta_1$ and $\Ker \delta_2$ are finitely generated algebras (in particular,
$\delta_i$ admits the categorical quotient morphism $\rho : X \to Q_i$).
We say that $(\delta_1, \delta_2)$ is a compatible pair if

 (i) the vector space ${\rm Span}
(\Ker \delta_1 \cdot \Ker \delta_2)$ generated by elements from
$\Ker \delta_1 \cdot \Ker \delta_2$ contains a nonzero ideal in
$\bk [X]$ (so called associated ideal of the pair) and

(ii) some element $a \in \Ker \delta_2$ is of degree
1 with respect to $\delta_1$, i.e., $\delta_1 (a) \in \Ker \delta_1
\setminus \{ 0 \}$.

  \edefi

\brem\label{adp.r1} (1) In \cite{KaKu08} the vector field $\delta_2$ was also allowed to be semi-simple, but we do not consider this case here.

(2) For every locally nilpotent vector field $\delta$ on $X$ its nonzero replica $f\delta$ has the same kernel
and for each $a \in \bk [X]$ the degree of $a$ with respect to $\delta$ coincides with its degree with respect to $f\delta$. 
Hence, for a compatible pair $(\delta_1, \delta_2)$
and $f_i \in \Ker \delta_i \setminus \{ 0\}, \, i=1,2$ the pair $(f_1\delta_1, f_2 \delta_2)$ is also compatible.
Furthermore, the conjugation by any element of $\Aut (X)$ transfers $(\delta_1, \delta_2)$ into another compatible pair.

(3) The assumption that $\Ker \delta_1$ and $\Ker \delta_2$ are finitely generated was unfortunately missed in \cite{KaKu08}.

(4) It is worth mentioning that in (ii) the field $a\delta_1$ is complete and, furthermore, (ii) holds for any commuting pair of nontrivial non-equivalent\footnote{Recall that two 
locally nilpotent derivations $\delta$ and $\sigma$ on $\bk [X]$ are equivalent if $\Ker \delta = \Ker \sigma$.} locally nilpotent derivations.

\erem

\bnota\label{fm.n1} For every affine algebraic variety $X$ we denote by $\AVF (X)$ the space of all algebraic vector fields on $X$.
If $K$ is an ideal in $\bk [X]$ then $\AVF_K (X)$ is the subspace of $\AVF (X)$ generated by all fields of
the form $f \sigma$ where $f \in K$ and $\sigma \in \AVF (X)$. Given a saturated set $\cN$ of locally nilpotent vector
fields on $X$ and a closed subvariety $Z$ in $X$ we denote by $\Lie_{\rm alg}^\cN (X, Z)$ the Lie algebra generated by the complete vector fields
vanishing on $Z$ that are of the form $b\delta$ where $\delta \in \cN$ and $b \in \bk [X]$ has degree $\deg_\delta b\leq1$.  
If $x \in X$ then $\mu_x$ is its vanishing maximal ideal in $\bk [X]$ and for every $\bk [X]$-module $M$
its localization at $\mu_x$ will be denoted by $(M)_{\mu_x}$.

\enota

\blem\label{fm.l1} Let $X$ and $Z$ be as in Notation \ref{fm.n1} and $I\subset \bk^{[n]}$ be the vanishing ideal of $Z$.
Let $M$ be a $\bk [X]$-submodule of $\AVF (X)$ such that for every $x \in X \setminus Z$ the localization 
$(M)_{\mu_x}$ coincides with the localization of $\brM =\AVF(X)$ at $\mu_x$.
Then there exists $k>0$ such that $M$ contains $\AVF_{I^k} (X)$.
\elem

\bproof Since $X$ is affine the $\bk [X]$-module $\brM$ is finitely generated, i.e., we have $\brM=\sum_{i=1}^m \brM_i$ where $\brM_i=\bk [X]\sigma_i$ for
some vector field $\sigma_i \in \AVF [X]$. Let $M_i=M\cap \brM_i$. Note that $\brM_i$ has a natural structure of the ring $\bk [X]$
and $M_i$ can be viewed is an ideal $K_i$ in $\bk [X]$. Since  the operations of localization and intersection commute
we have $(M_i)_{\mu_x}=(\brM_i)_{\mu_x}$ for every $x \notin Z$.
Hence, the zero locus of $K_i$ is contained in $Z$.  By Nullstellensatz there exists $k_i$ for which $I^{k_i} \subset K_i$ and, therefore,
$I^{k_i} \brM_i \subset M_i$. Letting $k=\max\{k_i|i=1, \ldots, m\}$ we get the desired conclusion.
\eproof

\bthm\label{adp.t1} Let $X$, $Z$ and $\cN$ be as in Notation \ref{fm.n1} and $I\subset \bk^{[n]}$ be the vanishing ideal of $Z$.
Suppose that every vector field in $\cN$ is tangent to $Z$ (i.e., the flows of these vector fields map $Z$ onto itself).
Let $G \subset \SAut (X)$ be the group generated by $\cN$
and let $X\setminus Z$ be $G$-flexible. 
Suppose that $X$ admits a pair of compatible locally nilpotent vector fields $\delta_1$ and $\delta_2\in \cN$.
Then for some $k\geq 1$ the algebra
$\Lie_{\rm alg}^\cN (X, Z)$ contains the space $ \AVF_{I^k} (X)$. 

\ethm

\bproof Let $\rho_i : X \to Q_i:={\rm Spec}\, \Ker \delta_i$ be the quotient morphism associated with $\delta_i$, let $Z_i$ be the closure of $\rho_i (Z)$ in $Q_i$, and let $a$ be as in Definition \ref{adp.d2}.
Choose a nonzero function $h_i \in \Ker \delta_i \simeq \bk [Q_i]$ that vanishes on $Z_i$ and note that for every $f_i \in \Ker \delta_i, \, i=1,2$ the fields
$f_1h_1\delta_1, af_2h_2\delta_2, af_1h_1\delta_1$, and $f_2h_2\delta_2$ are contained in $\Lie_{\rm alg}^\cN (X, Z)$. 
Hence
\be\label{adp.eq1} [f_1h_1\delta_1, af_2h_2\delta_2]-[af_1h_1\delta_1,f_2h_2\delta_2]=f_1f_2 h_1h_2 \delta_1(a) \delta_2\ee
also belongs to  $\Lie_{\rm alg}^\cN (X, Z)$.   By condition (i) in Definition \ref{adp.d2} $\Span (\Ker \delta_1 \cdot \Ker \delta_2)$ 
contains a nontrivial ideal $J$. Hence, $J_2:=h_1h_2J\subset J \cap I$ is also
nontrivial ideal in $\bk [X]$. By Formula (\ref{adp.eq1}) one has
 $J_2\delta_2 \subset \Lie_{\rm alg}^\cN (X, Z)$. 

Let $\tilde \cN$ be the set $\{ \delta''\}$ of locally nilpotent vector fields $\delta''\in \cN$ such that for some
$\delta' \in \cN$ the pair $(\delta', \delta'')$ is compatible. By Remark \ref{adp.r1} (2) this set is saturated and closed under conjugation by elements of $G$.
By Proposition \ref{fm.p2a} $\tilde \cN$ generates a subgroup $\tG$ of $\SAut (X\setminus Z)$ such that $X\setminus Z$ is $\tG$-flexible.
Let $\sigma_1, \ldots , \sigma_m \in \tilde \cN$, $\sigma_2=\delta_2$
and $J_i$ plays the same role for $\sigma_i$ as $J_2$ above for $\delta_2$, i.e., $J_i \sigma_i \subset \Lie_{\rm alg}^\cN (X, Z)$.
Applying \cite[Proposition 1.8]{AFKKZ} (with $H_i$ being the one-parameter unipotent group associated with $\sigma_i$) 
we can suppose that for every $x \in X$ the values of the fields $\sigma_1, \ldots , \sigma_m$
generate $T_xX$. Set $L_1=J_1 \cdot \ldots \cdot J_m$, i.e., $L_1 \sigma_i \subset \Lie_{\rm alg}^\cN (X, Z)$ for every $i$ and $L_1$ is contained in 
$I$. 

Put $M_1=\AVF_{L_1} (X)$.
Then $M=\sum_{i=1}^m L_1\sigma_i$ is a $\bk [X]$-submodule of  $M_1$ such that 
for every point $x$ in the complement (in $X$) to the zero locus $W_1$ of $L_1$  and every nonzero $v \in T_xX$ there exists a vector field $\sigma \in M$ whose value at $x$ is $v$.
Hence, $M/(\mu_x M) = M_1/(\mu_x M_1)$ when $x \in X\setminus W_1$.
In particular, the localizations $(M)_{\mu_x}$ and $(M_1)_{\mu_x}$ 
satisfy the relation $(M_1)_{\mu_x}=\mu (M_1)_{\mu_x}+(M)_{\mu_x}$ and the Nakayama lemma \cite[Corollary 2.7]{AM} implies that $(M)_{\mu_x}=(M_1)_{\mu_x}$.
Hence by \cite[Proposition 3.9]{AM} $M = M_1$, i.e., $ \AVF_{L_1}(X)\subset \Lie_{\rm alg}^\cN (X, Z)$.

Using conjugations by elements of $\tG$, we can transform $L_1$ into a sequence
of ideals $L_1, \ldots , L_k \subset I$ 
such that for every $x \in X \setminus Z$ there exists $i$ for which $x$ is not in the zero locus of $L_i$. 
Consider the $\bk [X]$-module $N= \sum_{i=1}^k\AVF_{L_i} (X)$. By construction $(N)_{\mu_x}$ coincides with
the localization of $\AVF (X)$ at $\mu_x$ for every $x \in X \setminus Z$. Hence, by Lemma \ref{fm.l1} 
$N$ contains $\AVF_{I^k}(X)$ for some $k >0$. Since $ \AVF_{L_i}(X)\subset \Lie_{\rm alg}^\cN (X, Z)$ for every $i$,
we see that $N\subset \Lie_{\rm alg}^\cN (X, Z)$ which yields the desired conclusion.
\eproof

\brem\label{adp.r2} (1) Let the assumptions of Theorem \ref{adp.t1} hold with the following execption:
we do not assume that the fields from $\cN$ are tangent to $Z$ and that $X\setminus Z$ are $G$-flexible, but
we suppose that $X$ is $G$-flexible. Then the conclusion of Theorem \ref{adp.t1} remains valid.
Indeed, consider the saturated subset $\cN_Z$ of $\cN$ that consists of all fields that vanish on $Z$
and let $G_Z\subset \SAut (X)$ be the group generated by $\cN_Z$. Then $X\setminus Z$ is $G_Z$-fleixible
by \cite{FKZ-GW}. Hence, replacing $\cN$ and $G$ by $\cN_Z$ and $G_Z$ respectively we get the assumptions
and, therefore, the conclusion of Theorem \ref{adp.t1}.

(2) If $\cN$ is the set of all locally nilpotent vector fields on $X$ 
and for some $k\geq 1$ the algebra $\Lie_{\rm alg}^\cN (X, Z)$ 
contains the space of all algebraic vector fields on $X$ that 
vanish on $Z$ with multiplicity at least $k$ then 
we say that the pair $(X,Z)$ has the algebraic density property. In this
terminology Theorem \ref{adp.t1} is a generalization of \cite[Theorem 4]{KaKu08}
which established the algebraic density property for pairs of the form $(\C^n, Z)$ 
where $Z$ is a closed subvariety  of $\C^n$ with $\dim Z \leq n-2$.
The algebraic density property in the complex case has some remarkable consequences.
In particular,  as in \cite{KaKu08}
we get two interesting facts which will not be used in the sections below.  
\erem

\bthm\label{adp.t2}  {\rm (cf. \cite[Theorem 4.10.6]{For})}   Let $X$ be a complex affine flexible variety and $Z$ be a closed subvariety of $X$
whose codimension is at least 2. 
Suppose that $X$ admits a pair of compatible vector fields. 
Let $\Phi_t : \Omega_0 \to \Omega_t=\Phi_t(\Omega_0) \subset X \setminus Z \, (t\in [0,1])$ be a
${\cC}^1$-isotropy consisting of injective holomorphic maps between Runge domains\footnote{Recall that an open
subset $\Omega$ of a Stein manifold $Y$ is called a Runge domain if every holomorphic function on $\Omega$ can be approximated
in the compact-open topology by holomorphic functions on $Y$.}
with $\Phi_0={\rm Id}_{\Omega_0}$. Suppose also that each $\Omega_t$ is Stein.
Then $\Phi_1$ can be approximated uniformly on compacts of $\Omega_0$ by holomorphic automorphisms of $X$ identical on $Z$.
\ethm

\bproof Consider generators $f_1, \ldots, f_m$ of the vanishing ideal $I\subset \bk [X]$ of $Z$.  These functions have no common zeros on any $\Omega_t$.
By the weak Nullstellensatz for Stein spaces
(e.g., see \cite[Theorem 4.25]{Oni}) there are holomorphic functions $g_1^t, \ldots , g_m^t$ on $\Omega_t$
for which $\sum_{i=1}^m f_ig_i^t =1$ on $\Omega_t$. Since $\Omega_t$ is a Runge domain on every compact $K_t \subset \Omega_t$
these functions $g_1^t, \ldots , g_m^t$ can be uniformly approximated by global holomorphic functions on $X$, i.e.,
we get a holomorphic function $h$ on $X$ vanishing on $Z$ which is as close to 1 on $K_t$ as we wish. Furthermore, 
replacing $h$ by $h^k$  for a given $k\geq 1$ we can suppose that $h$ is contained in  $\tI^k$  where $\tI$ is the vanishing ideal of $Z$
in the algebra $\Hol (X)$ of holomorphic functions  on $X$. Hence every holomorphic vector field $\nu_t$ on $\Omega_t$
can be approximated in the compact-open topology by holomorphic fields from  $\HVF_{\tI^k} (X)$ where $\HVF (X)$ is the space
of all holomorphic vector fields on $X$. Since $X$ is affine every holomorphic function (resp. vector field) on $X$ can be approximated by regular functions
(resp. algebraic vector fields) on $X$ in the compact-open topology
and also $\tI$ is generated by $I$ over $\Hol (X)$ (e.g., see \cite[Theorem 4]{Ka91}). Hence $\nu_t$ can be approximated in
the compact-open topology by elements of $\AVF_{I^k} (X)$ and, therefore, by Theorem \ref{adp.t1} by elements of $\Lie_{\rm alg}^\cN (X, Z)$.
Let $\tilde \nu_t$ be an element of $\Lie_{\rm alg}^\cN (X, Z)$ uniformly close to $\nu_t$ on $K_t$
and let $\psi_s$ and $\tilde \psi_s$ be the flows of $\nu_t$ and $\tilde \nu_t$ respectively with $s$ being the time parameter. 
Suppose that for some $s_0>0$ both of these flows
are defined for every $x \in K_t$. Recall that all generators of $\Lie_{\rm alg}^\cN (X, Z)$ can be chosen as complete
algebraic vector fields that vanish on $Z$. By \cite[Corollary 4.8.4]{For} $\tilde \psi_{s_0}$ can be approximated by compositions
of flows of these generators. Hence we have the following

\noindent (*) \hspace{.5cm} the element $\psi_{s_0}$ of the flow of $\nu_t$ can be uniformly approximated on $K_t$
by global holomorphic automorphisms of $X$ identical on $Z$.

To make use of (*) following the proof of \cite[Theorem 4.9.2]{For} we consider the trace
$\tilde \Omega = \{ (t,z): t\in [0,1], z \in \Omega_t\}$ of the isotropy $\{ \Phi_t\}$ in $\R \times X$
and we treat $\Phi_t$ as the flow of the continuous time dependent vector field 
$$V(t,z)=\dot{\Phi}_t (\Phi_t^{-1} (z))$$ where dot denotes the derivative on $t$.
The field $V$ is continuous on $\tilde \Omega$ and holomorphic on $\Omega_t$ for every fixed $t \in [0,1]$.
Divide the interval $[0,1]$ into $N$ subintervals $[t_k, t_{k+1}]$ of length $1/N$ for a given natural $N$
and consider the locally constant holomorphic vector field $\tilde V(t,z)$ which is equal to $V(t_k,z)$ on every interval $[t_k, t_{k+1}]$.
Let $\phi_t$ be the flow of $\tilde V$.  Similarly to the estimates in the proof of \cite[Theorem 4.8.2]{For} one can check
that as $N \to + \infty$ the flow $\phi_t$ converges
to $\Phi_t$ uniformly on compacts in $\Omega$ for all $t \in [0,1]$.
Since by (*) $\varphi_t$ can be approximated by global holomorphic automorphisms of $X$ identical on $Z$ so does $\Phi_t$.
This yields the desired conclusion.
\eproof

\bcor Let $X$ be a complex affine flexible variety of dimension $n$ and $Z$ be a closed subvariety of $X$
whose codimension is at least 2.  Suppose that $X$ admits a pair of compatible vector fields.
Then every $x \in X \setminus Z$ 
has a neighborhood $\Omega$ in $X\setminus Z$ that is a Fatou-Bieberbach domain
(i.e., $U$ is biholomorphic to $\C^n$).
\ecor

\bproof Since $X\setminus Z$ is flexible it suffices to prove this statement for some point $x_0$ in $X\setminus Z$ only.
Choose a dominant morphism $\varphi : X \to \C^n$ and choose $x_0\in X \setminus Z$ so that for a ball $B_0\subset \C^n$ with 
center at $\varphi (x_0)$ the component $B$ of the preimage $\varphi^{-1}(B_0)$, containing $x_0$, is naturally biholomorphic to $B_0$.
Taking $B$ small enough we can suppose that there are regular functions $g_1, \ldots , g_m \in \C [X]$ such
that each $|g_i|$ does not exceed 1 on $B$ while it is greater than 1 at each point of $\varphi^{-1}(B_0) \setminus B$.
This implies that $B_0$ is $Hol (X)$-convex. Hence it is a Runge domain in $X$ by the Oka-Weil theorem \cite[Theorem 2.2.5]{For}.
For an analytic coordinate system $(z_1, \ldots , z_n)$ on $B$ with the origin at $x_0$ consider the homothety $\Phi : (z_1, \ldots , z_n) \to (z_1/2, \ldots , z_n/2)$.
By Theorem \ref{adp.t2} $\Phi$ can be approximated by a global
holomorphic automorphism $F$ of $X$ identical on $Z$. Since $F(B)\subset B$ by the Brouwer fixed point theorem $F$ has a fixed point in $B$.
Without loss of generality we can suppose that this point is $x_0$ and reducing the size of $B$ we can suppose that $B$ contains no other fixed point but $x_0$. 
The eigenvalues $\lambda_1, \ldots , \lambda_n$ of the map ${\rm d}F$
at $x_0$ must be close to those of ${\rm d} \Phi$ which are $1/2$. In particular, we can suppose that $|\lambda_1|\geq |\lambda_2| \geq \ldots \geq
|\lambda_n|$ and $|\lambda_1|^2 < |\lambda_n|$. In particular, these eigenvalues satisfy the assumptions of \cite[Theorem 9.1]{RoRu}
and we can copy the argument of Rosay and Rudin. Namely, consider the basin $\Omega$ of attraction of $F$, i.e., $x \in \Omega$ if there exists
$N>0$ such that $F^N(x) \in B$. By continuity for every $y$ in a small neighborhood of $x$ one has $F^N(y) \in B$. Thus $\Omega$ is open (in
the standard topology). Treating $B$ as a neighborhood of the origin in $T_{x_0} X\simeq \C^n$ and letting $A={\rm d}F$ we can consider
an injective map $A^{-M}\circ F^M: F^{-N}(B)\to \C^n$ where $M\geq N$. Choose positive constants $\beta> |\lambda_1|$ and 
$\alpha<|\lambda_n|$ such that $\beta^2 < \alpha$. Since the local Jacobi matrix of the map $A^{-1}\circ F$ at $x_0$ is the identity matrix
the computation in \cite[Theorem 9.1]{RoRu}
\footnote{ See the formula immediately after formula (9) in  \cite[Theorem 9.1]{RoRu}. 
Formally, this formula is proven for $X=\C^n$ but it works in our case as well without change.}
implies that for  every compact $K \subset F^{-N}(B)$
there exist some $M_0>N$ and a positive real number $b$ such that for every $x \in K$ one has 
$$||A^{-M}\circ F^M (x)-A^{-M-1}\circ F^{M+1}(x)||=||A^{-M} (F^M(x) -A^{-1}\circ F (F^M(x)))||\leq b(\beta^2 / \alpha)^M$$
where $M\geq M_0$. Hence we have a well-defined holomorphic map $\Phi : \Omega \to \C^n$
where $\Psi =\lim_{M\to \infty} A^{-M}\circ F^M$. Since at every point of  $F^{-N}(B)$ the local Jacobian of $A^{-M}\circ F^M$
does not vanish we have the following alternative: either the local Jacobian of $\Psi$ does not vanish or it is identically zero.
However, the local Jacobian of $\Psi$ at $x_0$ is 1 and we have the former, i.e., $\Psi$ is an open map. Furthermore, it is
injective since otherwise the maps $ A^{-M}\circ F^M$ are not injective for large $M$. Note also that $\Psi =A^{-1}\circ \Psi \circ F$,
i.e., $\Psi$ and $A^{-1} \circ \Psi$ have the same range. Since the linear operator $A^{-1}$ is an expansion it follows that
$\Psi (\Omega) =\C^n$. Thus the basin $\Omega$ of attraction of $F$ is biholomorphic to $\C^n$ and we have the desired conclusion.\eproof

\brem\label{adp.r2r} For $X=\C^n$ the question about  Fatou-Bieberbach domains in the complement of
subvariety of codimension 2 was posed by Siu and answered by Buzzard and Hubbard \cite{BuHu} (see also \cite{For01}).

\erem

\section{Relative version of Theorem \ref{fm.t2}}

Let us prove first the following analogue of \cite[Theorem 3.1]{AFKKZ}.

\bthm\label{rv.t1} Let $\rho : X \to Q$ be a dominant morphism of quasi-affine algebraic varieties, $Q_0$ be
a Zariski open dense subset of $Q$, and $X_0 = \rho^{-1}(Q_0)$.  
Let  every fiber $\rho^{-1}(q), \, q \in Q_0$ be $G$-flexible where $G \subset \SAut (X_0)$
is a subgroup generated by a saturated set $\cN$ of locally nilpotent vector fields on $X_0$
 which are tangent to the fibers of $\rho$.
Suppose that $q_1, \ldots, q_m$ are distinct points in $\rho (X_0)$
and  $\alpha_i \in G|_{\rho^{-1}(q_i)}$.
Then there exists an automorphism $\alpha$ of $X$ over $Q$ 
such that $\alpha|_{\rho^{-1}(q_i)} =\alpha_i$  for every $i=1, \ldots , m$ and $\alpha|_{X_0} \in G$.

\ethm

\bproof  Suppose first that $Q=Q_0$. For a locally nilpotent vector field $\sigma$ denote by $\exp (t\sigma)$ 
the element of the one-parameter group associated with $\sigma$ at time $t\in \bk$.
By definition $\alpha_i$ is of the form 
$$\exp (t_{1,i}\sigma_{1,i}) \circ \cdots \circ \exp (t_{n(i),i}\sigma_{n(i),i})|_{\rho^{-1}(q_i)}$$
where $n(i)$ is a natural number depending on $i$ and $t_{j,i} \in \bk$. Choose regular functions $f_i$ on $Q$ such that $f_i (q_i)=1$ and $f_i(q_j)=0$ for every $j \ne i$.
Using the natural embedding $\bk  [Q]\subset \bk [X]$ we treat $f_i$ as a function on $X$. Then by the assumption $f_i\in  \Ker \sigma$ for every $\sigma \in \cN$,
i.e., the replica $f_i\sigma \in \cN$.  It remains to put  $\alpha =$
$$\exp ({t_{1,1}}f_1\sigma_{1,1}) \circ \cdots \circ \exp ({t_{n(1),1}}f_1\sigma_{n(1),1})  \circ \cdots \circ \exp ({t_{1,m}}f_m\sigma_{1,m}) \circ \cdots \circ \exp ({t_{n(m),m}}f_m\sigma_{n(m),m})$$
and we are done in the case of $Q=Q_0$.

In the general case this $\alpha$ is only an automorphism of $X_0$ and its extension to $X$ may have poles on $X\setminus X_0$.
However, one can choose $f_i$ from the above so that they vanish on $Q\setminus Q_0$ with sufficiently high multiplicity.
That is, we can assume that each $f_i\sigma_{j,i}$ is a regular vector field on $X$ that vanishes on $X\setminus X_0$. Then
the extension of $\alpha$ becomes a regular automorphism on $X$ whose restriction to $X \setminus X_0$ is the identity map.
\eproof

\brem\label{rv.r0} 
(1) Let $Q$ be affine. Then Theorem \ref{rv.t1} remains valid with the same proof when the finite set $q_1, \ldots, q_m$ is replaced by a collection of $m$ disjoint closed suvarieties of $Q$ which are contained in $Q_0$.

(2) One can consider a more general situation when every $\alpha_i$ is the restriction of an element of $G$ to a $k$-infinitesimal neighborhood $V_i$
of $\rho^{-1}(q_i)$.\footnote{For every reduced subvariety $Y$ of $X$ with a defining ideal $I\subset \bk [X]$ one can treat an automorphism of the $k$-infinitesimal neighborhood of $Y$
as an automorphism of the ring $\bk [X]/I^k$.} Then we can still find  $\alpha \in G$
for which $\alpha|_{V_i} =\alpha_i, \, i=1, \ldots , m$. For this it suffices to require that each $f_i$ vanishes with
multiplicity at least $k$ at $q_j$ (where $j \ne i$) and takes the value 1 with multiplicity at least $k$ at $q_i$. 
\erem

\bnota\label{rv.n1} Further in this section $X$ is a smooth algebraic variety of dimension $n$, $G \subset \SAut (X)$ is a subgroup generated by a saturated set $\cN$ of locally nilpotent vector fields on $X$,
$G_z\subset G$ is the isotropy group of a point $z \in X$, $\fm_z$ is the maximal ideal in the local ring $\cO_{X,z}$ at $z$, and $A_m(X,z)=\fm_z/\fm_z^{m+1}$
(in particular, $A_1(X,z)$ coincides with the cotangent space $T_z^*X$).
We consider the set $\Aut (A_m(X,z))$ of
$\bk$-algebra isomorphisms $f:A_m (X,z)\to A_m(X,z)$ satisfying the following condition:
\be\label{rv.eq1b}  \text{the Jacobian } J(f)= 1\mod\fm_z^{m+1}. \ee

Let $u_1, \ldots, u_n\in \mu_z$ be such that they generate the cotangent space.  Then we call the $n$-tuple $(u_1, \ldots , u_n)$ a local coordinate system
at $z$ (indeed, if the ground field $\bk =\C$ then $u_1, \ldots, u_n$ form a local analytic coordinate system in a small neighborhood of $z$ in the standard topology).
In terms of  this local coordinate system  elements of $\Aut (A_m(X,z))$ can be described
as follows. The $\bk$-algebra $A_m$ is contained in the quotient
$A/\fm_A^{m+1}$ of the local ring 
$A=\bk [[u_1,\ldots, u_n]]$  of formal power series with respect to the power of its maximal ideal $\fm_A$. Therefore, we treat any map $f\in \Aut(A_m(X,z))$ as  an $n$-tuple of polynomials
$(F_1,\ldots, F_n)\in (A/\fm_A^{m+1})^n$ in $n$ variables $u_1, \ldots , u_n$ of degree at most $m$ such that they vanish at the origin and  the determinant of the matrix $[\frac{\p F_i}{\p x_j}]_{i,j}$ is 1
modulo terms of degree higher than $m$. In particular, each $F_i$ is the sum of homogeneous $k$-forms where $k$ runs from 1 to $m$. Let $F_i'$ be the $m$-form present is this sum
and $\theta_{z,m}$ be the linear map from $\Aut(A_m(X,z))$ to the space of $n$-tuples of $m$-forms given by  
$$\theta_{z,m} (f)=(F_1', \ldots , F_n').$$ 
Suppose also that  $\lambda (f)$ is the $n$-tuple of linear parts of $f$. In particular, $\lambda (f) \in \SL (n,\bk)$ (because of the assumption on the Jacobian).
Note that $\SL (n,\bk)$ admits different natural actions on the space $\theta_{z,m}(\Aut(A_m(X,z)))$ of $n$-tuples $F(\bar u)$ of $m$ forms  in $n$ variables  (i.e., $\bar u = (u_1, \ldots , u_n)$)
for which we use the following notations
$$\lambda ._l F(\bar u) = \lambda (F(\bar u)), \,\, \,  \lambda ._r F(\bar u)= F(\lambda (\bar u)), \,\,\, {\rm and} \, \, \, \lambda .F(\bar u)=\lambda^{-1}(F(\lambda (\bar u))).$$

\enota

\blem\label{rv.l1}  Let Notation \ref{rv.n1} hold and $\Aut_{m-1}(A_m(X,z))$ be  the subgroup of the group $\Aut (A_m(X,z))$ consisting  of those automorphisms $f$ for which $f\equiv\id\mod\fm_z^{m}$
(i.e.,  $(f -\theta_{z,m} (f)) (\bar u)$ coincides with the $n$-tuple $(u_1, \ldots, u_n)$). Then we have the following.

{\rm (a)} For every $f \in  \Aut_{m-1}(A_m(X,z))$ and  $g \in  \Aut (A_m(X,z))$ one has

\[
    {\begin{array}{c}
    g-\theta_{z,m}(g)= f\circ g-\theta_{z,m}(f\circ g)= g\circ f-\theta_{z,m}(g\circ f),\\
    \theta_{z,m} (g\circ f)= \lambda(g) ._l (\theta_{z,m} (f)) +  \theta_{z,m} (g) \, \, {\rm and} \, \,
    \theta_{z,m} (f\circ g)= \lambda(g) ._r (\theta_{z,m} (f)) +  \theta_{z,m} (g),\\
       g^{-1}\circ f\circ g\in  \Aut_{m-1}(A_m(X,z)) \text{ and } \theta_{z,m} (g^{-1}\circ f\circ g)=\lambda(g).\theta_{z,m} (f).             
    \end{array} } 
\]
In particular, if $g$ is also in  $ \Aut_{m-1}(A_m(X,z))$ then  $\theta_{z,m} (g\circ f)=\theta_{z,m} (f\circ g) =\theta_{z,m} (f) +  \theta_{z,m} (g)$.

{\rm (b)} For $m \geq 2$ the set  $\cF_{z,m}:=\theta_{z,m} (\Aut_{m-1}(A_m(X,z)))$ is  the linear space of $n$-tuples $F(\bar u)$ of $m$-forms in $n$ variables of divergence zero and
the  $\SL (n,\bk)$-action on $\cF_{z,m}$ given by $\lambda . F(\bar u)$ is irreducible.

{\rm (c)} There is a natural homomorphism  $J_{z,m} : G_z \to  \Aut (A_m(X,z))$ such that
in the case when $X$ is $G$-flexible one has $J_{z,1} (G_z) =  \Aut (A_1(X,z))\simeq \SL (n,\bk)=\SL (T_z^*X)$.

{\rm (d)} If $\p$ is a locally
nilpotent vector field on $X$ with a zero of order $m\ge 2$ at $z$
then $\theta_{z,m}( J_{z,m}(\exp(t\p)))=t\theta_{z,m}( J_{z,m} (\exp(\p)))$.

\elem

\bproof Statement (a) is straightforward (see also \cite[Lemma 4.12]{AFKKZ}).  The first clause in statement (b) can be found in \cite[Lemma 4.13]{AFKKZ} and the second on in \cite[ IX.10.2]{Pr}.
Define $J_{m,z}(\alpha)$ as the operation of taking the $m$-jet of $\alpha \in G_z$.  Since  $\alpha \in \SAut (X)$ 
we see that the Jacobian $J(\alpha)\equiv 1$.
Hence $J_{z,m}  (G_z) \subset   \Aut (A_m(X,z))$. The fact that $J_{z,1}(G_z) \simeq \SL (T_z^*X)$
when $X$ is $G$-flexible follows from \cite[Corollary 4.3]{AFKKZ} which concludes (c).

For (d) note (as in  \cite[Lemma 4.12]{AFKKZ}) that  $\exp(\p) \in G_{z,m}$ and it induces the map
$\id + \hat \p \in \Aut_{m-1}(A_m)$ where $\hat \p$ denotes the derivation on $A_m$ induced by $\p$. Hence $\exp(t\p)$ induces $\id +t \hat \p$ which is (d) and we are done.

\eproof

\bnota\label{rv.n2}  In addition to Notation \ref{rv.n1} suppose that $\rho : X \to Q$ is a smooth morphism of smooth 
quasi-affine algebraic varieties
such that every fiber $Y=\rho^{-1}(q), \, q\in Q$ is  of dimension at least 2 and $G$-flexible (i.e., we  are
 under the assumptions of Theorem \ref{rv.t1} and $G$ is generated by a saturated set $\cN$
 of locally nilpotent vector fields such that every $\delta \in\cN$ is tangent to the fibers of $\rho$). In particular,
for every $z \in X$ and $q=\rho (z)$ a local coordinate system can be chosen in the form $(\bar u, \bar v):=(u_1, \ldots, u_k, v_1, \ldots , v_{n-k})$ where $u_i$ and $v_i$ are regular functions on $X$ such that $v_1, \ldots , v_{n-k}$ are the lifts
 of functions on $Q$
that form a local coordinate system at $q \in Q$ while  the restriction $(u_1, \ldots, u_k)$ yields a local coordinate system at $z \in Y$.

Note that in such a coordinate system  for every $\alpha \in G_z$ its image $J_{z,m} (\alpha)\in \Aut (A_m(X,z))$ is of the form 
\be\label{rv.eq1} J_{z,m} (\alpha) =(F_1 (\bar u, \bar v), \ldots, F_k (\bar u, \bar v), v_1, \ldots , v_{n-k})\ee
where $F_i$ is a polynomial of degree at most $m$ and  the determinant of the matrix $[\frac{\p F_i}{\p u_j}]_{i,j}$ is 1 up to terms of degree higher than $m$.
\enota

\blem\label{rv.l1b} Let $\delta \in \cN$ and $\delta_q$ be the restriction of $\delta$ to $Y=\rho^{-1}(q)$ for $q \in Q$.
Then a partial quotient morphism $\tau : X \to P$ (resp. $\tau_q : Y\to P_q$) of $\delta$ (resp. $\delta_q$)
can be chosen so that $\rho$ factors through $\tau$ and $\tau|_{Y}$ factors through $\tau_q$
(i.e., $\rho =\theta \circ \tau$ and $\tau|_{Y} = \kappa_q \circ \tau_q$ for some morphisms $\theta : P \to Q$
and $\kappa_q : P_q \to P$).
\elem

\bproof The quasi-affine variety $Q$ is contained as an open subset in an affine variety $\breve{Q}$. Under the natural embedding
$\bk [\breve{Q}] \hookrightarrow \bk [X]$ generators of $\bk [\breve{Q}]$ can be treated as elements of $\Ker \delta$.
Thus we can choose $\tau$ so that $\bk [P]$ contains these generators which implies that $\rho$ factors through $\tau$,
i.e., $\rho =\theta \circ \tau$.
Similarly, we can choose $\tau_q$ so that $\bk [P_q]$ contains generators of the ring of regular functions
on an affine variety in which $\theta^{-1}(q)$
is an open subset.  This implies that $\tau|_{Y}$ factors through $\tau_q$.
\eproof

\blem\label{rv.l2} Let $Z$ be a finite subset of $Y=\rho^{-1}(q)$ and $z\in Z$. 
Suppose that there exists $\delta \in \cN$ such that  $\delta|_{Y} \ne 0$ and partial quotient morphisms
$\tau$ and $\tau_q$ from Lemma \ref{rv.l1b} can be chosen so that 
\be\label{rv.eq2b} \kappa_q(P_q) \text{ is dense in } \theta^{-1}(q).\ee

Then such a field $\delta$ and  a coordinate system $(\bar u, \bar v)$ as in Notation \ref{rv.n2} can be chosen so that
the following holds:

{\rm (i)} $\delta$ induces a trivial derivation on $A_m(X,w)$ for every $w \in Z \setminus \{ z\}$;

{\rm (ii)} $u_i$ belongs to the kernel $\Ker \delta$ for every $i \geq 2$;

{\rm (ii)} the derivation $\sigma$ on $A_m(X,z)\subset A/\fm_A^{m+1}$ (where $A=\bk [[\bar u, \bar v)]]$) induced by $\delta$ coincides with $\sigma :=\p/\p u_1$.

\elem

\bproof  Treating $\delta$ as a derivation $\delta : \bk [X] \to \bk [X]$ consider its conjugate $\tilde \delta =g^*\circ \delta \circ (g^*)^{-1}: \bk [X] \to \bk [X]$
for $g\in G$. Similarly, let $\tilde \delta_q=g|_Y^*\circ \delta_q \circ (g|_Y^*)^{-1}$.
Note that   $\Ker \tilde \delta =g^*(\Ker \delta)$ and  $\Ker \tilde \delta_q =g|_Y^*(\Ker \delta_q)$. 
Therefore, for $\tilde \tau  =\tau \circ g : X \to P$ and $\tilde \tau_q=\tau_q \circ g|_Y : Y \to P_q$ one has
$\tilde \tau^*(\bk[P])\subset \Ker \tilde \delta$ and $\tilde \tau_q^*(\bk[P_q])\subset \Ker \tilde \delta_q$. Hence these morphisms
commute with the $\G_a$-actions associated with $\tilde \delta$ and $\tilde \delta_q$ respectively. Each fiber of $\tilde \tau$ (resp. $\tilde \tau_q$)
is the image of a fiber of $\tau$ (resp. $\tau_q$) under the action of $g$. This implies that the general fibers of $\tau$ and $\tau_q$ are isomorphic
to lines and these morphisms are partial quotient morphisms
for $\tilde \delta$ and $\tilde \delta_q$ respectively.   Furthermore, one can see that $\rho=\rho\circ g =\theta \circ \tau \circ g=\theta \circ \tilde \tau$ and,
similarly, $\tilde \tau|_Y=\kappa_q\circ \tilde \tau_q$.  In particular, the assumption (\ref{rv.eq2b}) holds for the locally nilpotent vector field $\tilde \delta$.
Since $\cN$ is saturated $\tilde \delta \in \cN$ and we can replace $\delta$ by $\tilde \delta$ while trying to achieve (i)-(iii) for the finite set $Z$. 

Note that the validity
of   (i)-(iii) for the pair $(\tilde \delta, Z)$ is equivalent to the validity of these conditions for the pair $(\delta, g^{-1}(Z))$.
Indeed, (i) is obvious and for the rest let $(u_1, \ldots, u_k, v_1, \ldots , v_{n-k}) $ be a local coordinate system at $z$ such that (i)-(iii) are true for the pair $(\tilde \delta, Z)$.
Let $(\hat u_1, \ldots, \hat u_k, v_1, \ldots , v_{n-k}) $ be a local coordinate system at $g^{-1}(z)$ such that $\hat u_i= u_i\circ g$, i.e.,
the local form of $g$ is $$(\hat u_1, \ldots, \hat u_k, v_1, \ldots , v_{n-k}) \to (\hat u_1\circ g^{-1}, \ldots, \hat u_k\circ g^{-1}, v_1, \ldots , v_{n-k}).$$
Since $u_i \in \Ker \tilde \delta_i$ for $i\geq 2$ and $\tilde \delta =g^*\circ \delta \circ (g^*)^{-1}$ we see that $\hat u_i \in \Ker \delta$, i.e., we have (ii) for the pair $(\delta, g^{-1}(Z))$. Condition (iii) is equivalent to
the fact that modulo $t^2$ the exponent $\exp (t \tilde \delta), \, t \in \bk$ has the following local form at $z$
$$(u_1, \ldots, u_k, v_1, \ldots , v_{n-k}) \to ( u_1+t+tf (\bar u,\bar v), u_2 \ldots, u_k, v_1, \ldots , v_{n-k})$$ where $f$ is a function
vanishing at the origin with multiplicity $m+1$ or higher. Since $\exp (t \tilde \delta)=g^* \circ \exp( t\delta) \circ (g^*)^{-1}$ the local form of $\exp (t \delta)$ is
$$(\hat u_1, \ldots, \hat u_k, v_1, \ldots , v_{n-k}) \to ( \hat u_1+t+tf (\hat u_1, \ldots , \hat u_k, \bar v), \hat u_2 \ldots, \hat u_k, v_1, \ldots , v_{n-k})$$
and we have (iii) for the pair $(\delta, g^{-1}(Z))$.

Therefore, we replace $Z$ by $g(Z)$ while keeping $\delta$ intact. By virtue of infinite transitivity (Theorem \ref{fm.t1a}) we can suppose now that
$Z$ consists of general points $\{ w_i \}$ of $Y$ with $z=w_1$.
Hence         $p_i'=\tau_q (w_i)$ are distinct
general (and, therefore, smooth) points $p_i'$ of $P_q$ such that for some neighborhood $U_i\subset P_q$ of $p_i'$ one has a natural 
$U_i$-isomorphism  $\tau_q^{-1}(U_i) \simeq U_i \times \A_\bk$.
Furthermore, we can suppose also that $\delta$ is nontrivial on $\tau_q^{-1}(p_1')$. 

Let $p_i=\kappa_q(p_i')$.
By (\ref{rv.eq2b}) we can suppose that $\{ p_i \}$ are general points of $\theta^{-1}(q)$, i.e., there is a function $f \in \bk [P]$ that vanishes at each $p_i, \, i\geq 2$ with multiplicity at least $m$ but has  $f(p_1) =1$
also with multiplicity at least $m$.
Then replacing $\delta$ by its  replica $f \delta \in \cN$ we get (i).

Since $p_1$ is a smooth point of $\theta^{-1}(q)$  a local coordinate system at $p_1 \in P$ can be chosen in the form
$(u_2, \ldots , u_k$,  $v_1,\ldots, v_{n-k})$ where each $u_i, \, i\geq 2$ is a regular function on $P$. 
Taking $u_1$ as an appropriate extension to $X$ of a coordinate function  on the  $\G_a$-orbit $\tau_q^{-1}(p_1') \simeq \A_\bk$ we can 
treat  $(u_1, u_2, \ldots , u_k, v_1,\ldots, v_{n-k})$
as a local coordinate system at $z \in X$.  This is the desired coordinate system with  $\delta$  being a desired derivation which concludes the proof.

\eproof

\brem\label{rv.r2b} Condition (\ref{rv.eq2b}) is very mild and it is automatically true when $q$ is a general point of $Q$.
\erem

\blem\label{rv.l3} Let $z \in Z\subset Y$ and $(\bar u, \bar v)$ be as in Lemma \ref{rv.l2}.  Suppose that $G_Z =\bigcap_{w \in Z}G_w$ and $G_{Z,z}^m$ is the subgroup
of $G_Z$ such that it induces the identity map on the $m$-th infinitesimal neighborhood of every point  $w \in Z \setminus \{ z\}$.  Then the image $I_m:=J_{z,m} (G_{Z,z}^m)$ contains 
all automorphisms of $A_m(X,z)$ with the following coordinate form
\be\label{rv.eq2} (\bar u, \bar v):=(u_1, \ldots, u_k, v_1, \ldots , v_{n-k}) \to (u_1 + \ell_1^m(\bar v), \ldots, u_k+\ell_k^m (\bar v), v_1, \ldots , v_{n-k})\ee
where every $\ell_i^m (\bar v)$ is an $m$-form in variables $v_1, \ldots , v_{n-k}$. Furthermore,   
for every $\lambda (\bar u)=(\lambda_1 (\bar u), \ldots, \lambda_k(\bar u))$ in $\SL (T_z^*Y)$ and each $k$-tuple 
$(\ell_1^1(\bar v), \ldots,   \ell_k^1(\bar v))$ of 1-forms the subgroup $I_1$ contains
 the following automorphism of $T_z^*X$
\be\label{rv.eq3} (\bar u, \bar v) \to (\lambda_1 (\bar u) + \ell_1^1(\bar v), \ldots, \lambda_k(\bar u)+\ell_k^1 (\bar v), v_1, \ldots , v_{n-k}).\ee

\elem

\bproof
Let $\sigma$ be as in Lemma \ref{rv.l2}. 
Then for every $h \in \Ker \delta$ the automorphism
$\exp ( h\delta) \in G_{Z,z}^m$ induces the automorphism $\exp (h'\sigma )$ on $A_m(X,z)$ (where $h'$ is the image of $h$ in $A_m(X,z)$)
of the following coordinate form
\be\label{rv.eq4} (\bar u, \bar v) \to (u_1 + h', u_2, \ldots, u_k, v_1, \ldots , v_{n-k}).\ee
In particular, choosing $h$ so that $h'$ is equal to an $m$-form $l_1^m(\bar v)$ we see that the automorphism 
\be\label{rv.eq5} (\bar u, \bar v) \to (u_1 + \ell_1^m (\bar v), u_2, \ldots, u_k, v_1, \ldots , v_{n-k})\ee is contained in $I_m$.
Note  the variety  $(Y \setminus Z) \cup \{ z \}$ is $G_{Z,z}^m$-flexible by \cite{FKZ-GW}.   Therefore,
by Theorem \ref{fm.t2} for every  $\lambda (\bar u)\in \SL (T_z^*Y)$ 
the subgroup $I_1$ contains an automorphism as in Formula (\ref{rv.eq3}) for certain 1-forms $\ell_i^1(\bar v)$.

Using notations of Lemma \ref{rv.l1} we suppose that $f$ and $g\in G_{Z,z}$ are such that
$J_{z,1} (g)$ is given by Formula (\ref{rv.eq3}) and  $J_{z,m} (f)$ is given by Formula (\ref{rv.eq5}). In particular,   $J_{z,m} (f) \in \Aut_{m-1}(A_m(X,z))$
and $\theta_{z,m} (J_{z,m} (f))=(\ell_1^m(\bar v), 0, \ldots, 0)$ while 
$$\theta_{z,1} (J_{z,m} (g))=(\lambda_1 (\bar u) + \ell_1^1(\bar v), \ldots, \lambda_k(\bar u)+\ell_k^1 (\bar v), v_1, \ldots , v_{n-k}).$$
 By Lemma \ref{rv.l1} (a) $J_{z,m} (g^{-1}\circ f \circ g) \in \Aut_{m-1}(A_m(X,z))$
and  $\theta_{z,m} (J_{z,m} (g^{-1}\circ f \circ g))$ is obtained from $\theta_{z,m} (J_{z,m} (f))$ via conjugation by
$\theta_{z,1} (J_{z,m} (g))$. Using such conjugations 
we see that automorphisms of the form
\be\label{rv.eq6} (\bar u, \bar v) \to (u_1, \ldots, u_{i-1}, u_i+\ell_i^m (\bar v), u_{i+1}, \ldots,  u_k, v_1, \ldots , v_{n-k}) \ee
are also contained in $I_m$. Taking the product of automorphism in Formula (\ref{rv.eq6}) with $i$ running from 1 to $k$ we
obtain every automorphism from Formula (\ref{rv.eq2}) as an element of $I_m$.  This implies in turn that 
any automorphism from Formula (\ref{rv.eq3}) is contained in $I_1$ regardless of the choice of  $\ell_i^1(\bar v)$
 which concludes the proof.
\eproof

Now we can formulate the main result of this section.

\bthm\label{rv.t2} Let Notation \ref{rv.n2} hold and  $V (A_m(X,z))$ be  the subset of $\Aut (A_m(X,z))$
which consists of automorphisms as in Formula (\ref{rv.eq1}) satisfying the assumption on the determinant of $[\frac{\p F_i}{\p u_j}]_{i,j}$ in Formula (\ref{rv.eq1b}).
Let $Z$ be a finite subset of $X$ such that for every
$q \in \rho (Z)$ there exists $\delta \in \cN$ (where $\delta$ may depend on $q$) for which condition (\ref{rv.eq2b}) in Lemma \ref{rv.l2} holds.  Suppose that
$G_Z =\bigcap_{z \in Z}G_z$ and $J_{Z,m} : G_Z \to \prod_{z \in Z} V (A_m(X,z))$ is the natural homomorphism.
Then $J_{Z,m}$ is surjective
(in brief, one can choose an element $\alpha \in G$ so that for every $z \in Z$ the $m$-jet of $\alpha$ at $z$ coincides with a prescribed jet from Formula (\ref{rv.eq1})
satisfying  the assumption on the determinant).
\ethm

\bproof By Theorem \ref{rv.t1} and Remark \ref{rv.r0} (2) it suffices to consider the case when $\rho (Z)$ is a singleton $q \in Q$. 
Furthermore, assume that for every $z \in Z$ we can find an element $\alpha_z$ in the subgroup $G_{Z,z}^m \subset G_Z$ from Lemma \ref{rv.l3} 
such that at $z$ this automorphism has a prescribed $m$-jet   from Formula (\ref{rv.eq1})
satisfying  the assumption on the determinant.  Recall that for  every  
$w \in Z \setminus \{ z \}$ the $m$-jet  of $\alpha_z$ at $w$ is the identity map in the $m$-th infinitesimal neighborhood of $w$.
Hence the composition of such automorphisms $\alpha_z$ with $z$ running over $Z$ (in any order) yields an automorphism $\alpha \in G_Z$
such that at every point of $Z$ it has a prescribed
$m$-jet  from Formula (\ref{rv.eq1})
satisfying  the assumption on the determinant. Therefore,  it is enough to consider the case when $Z$ is a singleton $z$ which we do below.

We shall use induction on $m$ with the case of $m=1$ provided by Lemma \ref{rv.l3}.
Assume now that the statement is true for $m-1$. 
Let $h \in V(A_m (X,z))$.  Then $h -\theta_{z,m} (h) =\tilde h  \in V(A_{m-1}(X,z))$. By the assumption
 there exists $\alpha \in G_{Z}$ for which $J_{z,m-1}(\alpha) =\tilde h$. Let $J_{Z,m}(\alpha) =g$ and let $\lambda (g)$ be as in Notation \ref{rv.n1}. 
 Consider the $n$-tuple $\tilde f=\theta_{z,m}(h)-\theta_{z,m}(g)$ of $m$-forms and let $f \in \Aut_{m-1}(A_m(X,z))$ be such that $\theta_{z,m} (f)=\lambda (g)^{-1} ._r \tilde f$.
 Then by Lemma \ref{rv.l1} (a) we have $f\circ g=h$. That is, it suffices to show that $f$ belongs to $J_{z,m}(G_Z)$, or, equivalently $\theta_{z,m} (f)$ is contained in
  $$I:=\theta_{z,m} (J_{z,m}(G_Z) \cap \Aut_{m-1}(A_m(X,z))).$$
  Note that since $h$ and $g$ are in $V(A_{m-1}(X,z))$  the last $n-k$ coordinates of the $n$-tuple $\tilde f$ (and, therefore, $\theta_{z,m}(f)$) are equal to
 zero for $m>1$.        
 Note also that any element of $I$ (and, thus, $\theta_{z,m}(f)$) is of the form $\sum_\mu \mu p_\mu (\bar u)$ where $\mu$ is a monomial in coordinates $\bar v$ 
 (with $(\bar u, \bar v)$ from Notation \ref{rv.n2}) and 
 $p_\mu (\bar u)$ is an $n$-tuple of  homogeneous polynomials in $\bar u$ of degree $l=m-\deg \mu$ such that
 the last $n-k$ coordinates of this $n$-tuple are equal to zero and its divergence is also zero (the latter fact follows from the assumption of 
 the Jacobian, see \cite[Lemma 4.13]{AFKKZ}). Thus applying  
 Lemma \ref{rv.l1} (a) again we see  that it suffices to show that 
$\mu p_\mu (\bar u)$ belongs to $I$. Furthermore, we can suppose that $l>0$ since the case of $l=0$ is taken care of by Lemma \ref{rv.l3}.

Recall that $\delta$ can be chosen so that conditions (i)-(iii) from Lemma \ref{rv.l2} are satisfied. In particular, $u_2$ belongs to the kernel
$\Ker \delta$ as well as every $v_j$ and, therefore, $\mu$.  Since $\cN$ is saturated the replica  $\mu u_2^{ l-1}\delta $ belongs to $\cN$.
Thus, $\exp (\mu u_2^{ l-1}\delta ) \in G_{Z}$, and 
$$e:=J_{z,m} (\exp (\mu u_2^{ l-1}\delta ) )$$
 belongs to $J_{z,m}(G_Z) \cap \Aut_{m-1}(A_m(X,z))$. 
Consider a finte set $\{ \gamma_i \}$ in $G_Z$ and  $g_i =J_{z,m} (\gamma_i)$.
Then  $$e_i:=J_{z,m} (\gamma_i^{-1} \circ \exp (\mu u_2^{ l-1}\delta) \circ \gamma_i ) \in J_{z,m}(G_Z) \cap \Aut_{m-1}(A_m(X,z)) $$ 
and $\theta_{z,m} (e_i)=\lambda(g_i) . \theta_{z,m} (e)$ is the result of the natural action of
 $\lambda (g_i) \in \SL (T_zX)$ from Lemma \ref{rv.l1} (b).  
 Consider the projection $\kappa$ of the space of $n$-tuples of $m$-forms  in $\bar u$ and $\bar v$ to the space of 
 $k$-tuples
 forgetting the last $n-k$ coordinates. Note that $\kappa (\theta_{z,m}(e))$ belongs to $\mu \theta_{z,l}( \Aut_{ l-1} (A_{l}(Y,z)))$.
 Recall also that by Lemma \ref{rv.l3} we can suppose that $\lambda(g_i)$ is given by 
 $$ (\bar u, \bar v) \to (\lambda_1 (\bar u), \ldots, \lambda_k(\bar u), v_1, \ldots , v_{n-k})$$ where $(\lambda_1 (\bar u), \ldots, \lambda_k(\bar u))$ is a prescribed
 element of $\SL(T_zY) \simeq \SL (k, \bk)$. Therefore,  $\kappa ( \theta_{z,m}(e_i))\in \mu \theta_{z,l}( \Aut_{ l-1} (A_{l}(Y,z)))$ and
 by Lemma \ref{rv.l1} (b) we can suppose that these elements  
 generate the vector space $\mu \cF_{z,l}$ where $\cF_{z,l}= \theta_{z,l}( \Aut_{ l-1} (A_{ l}(Y,z))$ is the space of $k$-tuples of $l$-forms in $\bar u$
 whose divergence is zero. 
 Lemma \ref{rv.l1} (d) implies that $t_i\lambda(g_i)  .  \theta_{z,m} (e)$ is also 
 contained in $I$ for every $t_i \in \bk$. 
  Applying Lemma \ref{rv.l1} (a) again we see that  every  linear combination of the elements $\lambda(g_i). \theta_{z,m} (e)$  belongs to $I$
 and, therefore, $\kappa (I)$ contains the space $\mu \theta_{z,l}\cF_{z,l}$. Hence $\mu p (\bar u) \in I$ which 
 concludes the proof.
\eproof

\bcor\label{rv.c1} Let Notation \ref{rv.n2} hold and condition (\ref{rv.eq2b}) in Lemma \ref{rv.l2} be valid for every
$q \in Q$. Suppose that 
$Z$ is a finite subset of $X$, and $m$ is a natural number. 
Let $S_z'$ and $S_z''$ be \'etale sections of $\rho : X \to Q$ through $z \in Z$.
Then there exists an automorphism $\alpha \in G$ such that for every $z \in Z$ one has $\alpha (z)=z$ and $\alpha (S_z')$ is tangent to $S_z''$
at $z$ with multiplicity at least $m$.
\ecor

\section{General projections for flexible varieties. I}




\bnota\label{gp1.d1} 
(1) If $\kappa : X \to P$ is a morphism of algebraic varieties then we denote
by $\Aut (X/P)$ (resp. $\SAut (X/P)$) the subgroup of $\Aut (X)$ (resp. $\SAut (X)$) that preserves each fiber of $\kappa$.

 (2) If $\kappa : X \to P$ is a smooth morphism of smooth varieties we denote by $T(X/P)$ the relative tangent bundle
(which is the kernel of the induced map $TX \to \kappa^* (TY)$). Furthermore, if $Z$ is a subvariety of $X$ then we still
denote by $T(Z/P)$ the intersection of the Zariski tangent bundle $TZ$ with $T(X/P)$. Similarly, for every $z \in Z$ we
let $T_z(Z/P)=T_zZ \cap T(Z/P)$.
\enota

The aim of this section is to describe analogues of general linear projections of $\A_\bk^n$ for flexible varieties. More precisely, we shall prove the following fact.

\bthm\label{gp1.t1}  
Let $X$ and $P$ be smooth algebraic varieties and $Q$ be a normal algebraic variety.  Let $\rho : X \to Q$ and $\tau: Q\to P$ be dominant morphisms such that $\kappa : X \to P$ is smooth for $\kappa=\tau \circ \rho$. 
Suppose that  $Q_0$ is a non-empty Zariski open smooth subset of $Q$ so that for  $X_0=\rho^{-1} (Q_0)$ the morphism $\rho|_{X_0} : X_0 \to Q_0$ is smooth.
Let $G \subset \Aut (X/P)$ be an algebraically generated group acting 2-transitively on each fiber of $\kappa : X \to P$ and
$Z$ be a locally closed reduced subvariety in $X$.

{\rm (i)} Let $\dim Z\times_P Z\leq 2\dim Z- \dim P$ and $\dim Q \geq \dim Z+m$ where $m \geq 1$. Then  there exists an algebraic family $\cA\subset G $ of automorphisms such that
for a general element $\alpha \in \cA$ one can find a subvariety $R$ of $\alpha (Z)\cap X_0$ of dimension $\dim R \leq \dim Z-m$ for which
$\rho (R) \cap \rho (\alpha (Z)\setminus R)=\emptyset$ and 
the morphism $\rho|_{(\alpha (Z)\cap X_0)\setminus R}:  (\alpha (Z) \cap X_0) \setminus R\to Z_\alpha'\cap Q_0$ 
is injective where 
$Z_\alpha'=\rho \circ \alpha (Z)$. In particular, if $\dim Q\geq \dim Z +1$ (resp. $\dim Q \geq 2\dim Z +1$)
for a general element $\alpha \in \cA$ 
the morphism $\rho|_{\alpha (Z)\cap X_0}:  \alpha (Z) \cap X_0 \to Z_\alpha'\cap Q_0$ 
is birational (resp. a bijection).

{\rm (ii)} Let $\dim Z\times_P Z\leq 2\dim Z- \dim P$, $\dim Q \geq  \dim Z+1$
and $F$ be a closed subvariety of $Z$ such 
$\dim F\times_P Z < \dim Q - \dim P$ (which is the case when $F$ is a finite set and $P$ is a singleton).
Then  for a general element $\alpha$ in the family $\cA$ from (i)  and every $z \in \alpha (F) \cap X_0$
one has $\rho^{-1}(\rho (z)) \cap \alpha (Z) =z$.

{\rm (iii)} Suppose that $G$ is generated by a saturated set $\cN$ of locally nilpotent vector fields on $X$
(in particular, every fiber of $\kappa$ is $G$-flexible).
Let $P_0=\tau (Q_0)$, $Z_0=Z\cap \kappa^{-1}(P_0)$ and $\dim T(Z_0/P_0)  \leq \dim Q -\dim P$. \footnote{
In the case of an irreducible $Z$ the above inequality follows from this one: $\dim TZ  \leq \dim Q -\dim P+ \dim \kappa (Z)$.
Indeed, in the presence of the latter  the exact sequence  $0\to  T_z(Z_0/P_0) \to T_zZ_0 \to T_{\kappa (z)} \kappa (Z_0)$  implies the former.
In particular, if $\kappa (Z)$ is dense in $P$ one has to require that $\dim TZ \leq \dim Q$.}
Then there exists an algebraic family $\cA\subset G $ of automorphisms such that
for a general element $\alpha \in \cA$, 
every $z \in \alpha (Z_0) \cap X_0$ the induced map $\rho_* : T_z\alpha (Z_0) \to T_{\rho (z)}Q$ of the tangent spaces is injective.

{\rm (iv)} Suppose that 
$G$ is again generated by a saturated set $\cN$ of locally nilpotent vector fields on $X$. 
Let $\dim Z\times_P Z\leq 2\dim Z- \dim P, 2\dim Z +1 \leq \dim Q \text{ and  }\dim T(Z_0/P_0)  \leq \dim Q -\dim P$. 
Then the family $\cA$ from (i) can be chosen so that 
for a general element $\alpha \in \cA$   the morphism  $\rho|_{\alpha (Z)\cap X_0}:  \alpha (Z) \cap X_0 \to Z_\alpha'\cap Q_0$ is 
injective and it induces an injective map of the Zariski tangent bundle of $\alpha (Z) \cap X_0$ into the Zariski tangent bundle of  $Q_0$.

\ethm

\bproof  For every variety $\cX$ over $P$ denote by $S_\cX $ the variety $S_\cX=(\cX\times_P \cX) \setminus \Delta_\cX$ where $\Delta_\cX$ is the diagonal in $\cX \times \cX$. 
Then every automorphism in $\Aut (X/P)$ can be lifted to an automorphism of $S_X$. In particular, we have a $G$-action on $S_X$ and by the assumption this action is transitive
on every fiber of the projection $S_X \to P$.
Consider the subvariety  $Y\subset S_X$ that is the intersection of $X_0\times_{Q_0} X_0$ and $S_X$ in $X\times_PX$. The codimension of $Y$ in $S_X$ is $\dim Q-\dim P$ and, because of smoothness, all fibers of the natural morphism $Y \to P$
are of the same dimension. Hence, by Remark \ref{agga.r1} (1) $\dim Y \times_P S_Z= \dim Y + \dim S_Z -\dim P$. By 
Theorem \ref{agga.t2} (ii) we can choose algebraic subgroups $H_1, \ldots, H_m$ of $G$ such that for a general element $(h_1, \ldots , h_m) \in H_1 \times \cdots \times H_m$
one has $$\dim W \leq \dim Y+  \dim S_Z -\dim S_X= \dim S_Z-\dim Q + \dim P\leq 2 \dim Z-\dim Q$$ where 
$W= Y \cap \alpha (S_Z)$ for $\alpha =h_1 \cdot \ldots \cdot h_m$.
Hence, in case (i) the dimension of $W$ is at most $\dim Z-m$. Let $R$ be the image of $W$ under one of the two natural projections $X\times_QX \to X$ 
(in particular, $R\subset \alpha (Z)\cap X_0$ and  $\dim R \leq \dim Z-m$).
Note that for $z \in \alpha (Z)\cap X_0$ one has $\rho^{-1}(\rho (z))\cap \alpha (Z)=z$ iff $z \notin R$.
Hence the restriction of $\rho$ to $(\alpha (Z)\cap X_0)\setminus R$ is injective.
Therefore, letting $\cA = H_1 \times \cdots \times H_m$, we get (i).

 In (ii) we let $S_{F,Z} =(F\times_P Z) \cap S_{X}$.
By the assumption of (ii) we have $\dim S_{F,Z} < \dim Q-\dim P$.
By Theorem \ref{agga.t2} (i)  for a general element $(h_1, \ldots , h_m) \in H_1 \times \cdots \times H_m$ the dimension of  the intersection of 
$\alpha (S_{F,Z})$  (where $\alpha= h_1 \cdot \ldots \cdot h_m$) with $Y$
is at most $ \dim S_{F,Z} +\dim Y -\dim S_X=  \dim S_{F,Z} +\dim P-\dim Q<0$, 
i.e this intersection is empty. It remains to note that the fact that $ \alpha (S_{F,Z}) \cap Y=\emptyset$
is exactly the statement that  for every $z \in \alpha (F) \cap X_0$
one has $\rho^{-1}(\rho (z)) \cap \alpha (Z) =z$.

In (iii)  for every variety $\cX$ and a subvariety $\cY$ of the tangent bundle $T\cX$ let $\cY^*=\cY \setminus \cS $ where $\cS$
is the zero section of the natural morphism $T\cX \to \cX$.
Every automorphism $\alpha \in \Aut (X/P)$ generates an automorphism of $T(X/P)$. In particular, $G$ acts on $T(X/P)^*$
and by  Theorem \ref{fm.t2} this action is transitive on every fiber of $T(X/P)^* \to X \overset{\kappa}{\to} P$.
 Note that the codimension of  $Y^* = T(X_0/Q_0)^*$ in $T(X/P)^*$ is equal to $\dim Q-\dim P$ 
 and all fibers of the natural projection $Y^*\to P_0$ are of the same dimension. 
Hence, by Remark \ref{agga.r1} $\dim Y^* \times_P T(Z_0/P_0)^*= \dim Y^* + \dim T(Z_0/P_0)^* -\dim P$.
 By Theorem \ref{agga.t2} and the inequality $\dim T(Z_0/P_0)  \leq \dim Q -\dim P$
  we can choose one-parameter unipotent algebraic subgroups $\tH_1, \ldots, \tH_{\tilde m}$ of $G$ such that 
for a general element $(\tilde h_1, \ldots , \tilde h_{\tilde m}) \in \tH_1 \times \cdots \times \tH_{\tilde m}$ and 
$Z''=(\tilde h_1 \cdot \ldots \cdot \tilde h_{\tilde m}) (Z_0)$ one has $\dim Y^* \cap T(Z''/P_0)^* \leq 0$.
Note that if $Y^* \cap T(Z''/P_0)^*$ contains a point then $\dim Y^* \cap T(Z''/P_0)$ must be at least 1
 (since this point is a vector in $TZ''$ and then $Y^* \cap T(Z''/P_0)^*$ 
contains all nonzero vectors proportional to that one).
That is, $Y^*\cap T(Z''/P_0)^*=\emptyset$.
This implies that for every $z\in Z''\cap X_0$ the restriction of $\rho_*$ to 
$T_z(Z''/P_0)$ is injective. Consequently, the restriction of $\rho_*$ to 
$T_zZ''$ is injective i.e., we have (iii).

In the last statement we note that (i) and (iii) in combination with Proposition \ref{agga.p2} imply that for a general element 
$\alpha$ of $H_1 \times \ldots \times H_m \times \tH_1 \times \ldots \times \tH_{\tilde m}$ 
the restriction of $\rho$ to $\alpha (Z)\cap X_0$ is injective and it induces an injective map of the Zariski tangent bundle of 
$\alpha (Z) \cap X_0$ into the Zariski tangent bundle
of  $ Q_0$. Thus we have (iv).
 \eproof

\brem\label{gp1.r1}  
(1) Since the family $\cA$ in Theorem \ref{gp1.t1}  is the product of connected groups we see that $\cA$ is irreducible and contains the identity map.

(2) Recall that by Proposition \ref{agga.p2} one can choose 
the family $H_1 \times \ldots \times H_m$  in Theorem \ref{agga.t2} 
independent of $Y$ and $Z$. Hence the argument in Theorem \ref{gp1.t1}
implies that in each of the claims (i)-(iv) one can suppose that 
the family $\cA$ does not depend  on the choice of $Z$.
 Furthermore, Theorem \ref{gp1.t1} remains valid (by the same Proposition \ref{agga.p2})
if one replaces $\cA$ with a family $H\times \cA$ (or $\cA \times H$) where $H$ is any connected algebraic subgroup
of $G$. In particular, if $\cA_k\subset G$ ($k=1, \ldots , s$) is an algebraic family of automorphisms of 
the form $H_1 \times \ldots \times H_m$ such that some property $\cP_k$ is satisfied for a general element of 
$\cA_k$ then there exists an algebraic family $\cA\subset G$ of automorphisms such that for a general
element of $\cA$ all these properties are true simultaneously.

(3) Let $P$ be a singleton and $X$ be $G$-flexible. Under the assumption of (iv) we can find $\alpha \in \cA$ such that for a given $z_0 \in Z$ one has $z:=\alpha (z_0) \in X_0$.
In particular,  for a general $\alpha$ the morphism $\rho|_{\alpha (Z)} : \alpha (Z) \to\bar  Z_\alpha'$ induces the injective map
$\rho_* : T_z\alpha (Z) \to T_{\rho (z)}Q$.
Note also that under the assumption of
 (iv) we can suppose that  for some neighborhood $U$ of $z$ in $\alpha (Z)$ the morphism $\rho|_U : U \to Q$ induces an
 injective map of the Zariski tangent bundles.
However, even if  $\rho|_{\alpha (Z)}: \alpha (Z) \to Q$ induces an
injective map of the Zariski tangent bundles it may not be
proper (and, in particular, $ Z_\alpha'$ may not be closed in $Q$). 
As a counterexample one can consider a bijective
morphism of $\C^*$ onto a polynomial curve which has only one singular point and this singularity is a node.


(4) Note that the  natural morphism $T(X/P) \to P$ (resp. $Y^* \to P_0$) in the proof of Theorem \ref{gp1.t1} (iii) has 
equidimensional fibers. Hence by Proposition \ref{agga.p3} we conclude that the assumption that
$\rho$ and $\tau$ are dominant and $\kappa$ is smooth can be omitted in the formulation of Theorem \ref{gp1.t1} (iii).

(5) Similarly, if $\kappa$ is not smooth or even dominant we can suppose as in Remark \ref{agga.r2}
that $\kappa (X)$ as a disjoint union $\bigcup_{k=1}^nP_k$ of
smooth varieties such that for $X_k=\kappa^{-1}(P_k)$ the morphism $\kappa|_{X_k} : X_k \to P_k$ is smooth.
Letting $Y_k=X_k \cap Y$ and $Z_k=X_k\cap Z$ and suing further stratification we can assume that
for every $k$ the morphism $\kappa|_{Z_k}: Z_k\to P_k$ has all fibers of the same dimension
(and, in particular $\dim Z_k\times_{P_k}Z_k = 2 \dim Z_k -\dim P_k$ \footnote{Of course, this does not imply 
that $\dim Z\times_PZ\leq 2 \dim Z -\dim P$.}).
Suppose that  for every $k$ and $Q_k=\tau^{-1}(P_k)$ one has
 \be\label{gp1.eq1}  \dim Q_k\geq \dim Z_k +m \ee 
By Theorem \ref{gp1.t1} (i)  (with $X,Q,P$ and $Z$ replaced by $X_k,Q_k,P_k$ and $Z_k$) there exists an algebraic family $\cA_k\subset G $ of automorphisms such that
for a general element $\alpha \in \cA_k$ one can find a subvariety $R_k$ of $\alpha (Z_k)\cap X_0$ of dimension $\dim R_k \leq \dim Z_k-m$ for which
$\rho (R_k) \cap \rho (\alpha (Z_k)\setminus R_k)=\emptyset$ and 
the morphism $\rho|_{(\alpha (Z_k)\cap X_0)\setminus R_k}:  (\alpha (Z_k) \cap X_0) \setminus R_k\to \rho \circ \alpha (Z_k)\cap Q_0$ 
is injective. By (2) one can suppose that $\cA_k=\cA$ is independent of $k$.  Thus we see that 
the conclusion of Theorem \ref{gp1.t1} (i) (and, therefore, (iv))  remains  valid even in the case when $\kappa$ is not dominant provided that
one replaces the assumption that  $\dim Z\times_P Z\leq 2\dim Z- \dim P$ and $\dim Q \geq \dim Z+m$
by the inequalities  in Formula (\ref{gp1.eq1}).

\erem

\bcor\label{gp1.c1} Let the assumptions of Theorem \ref{gp1.t1} (iv) hold with $Q=Q_0$ being quasi-affine.
 Suppose that  $Z$ is a once-punctured curve \footnote{That is, $Z$ is the complement to a point in a complete curve.}.
Then there exists an algebraic family $\cA \subset G$ of automorphisms of $X$ such that for a general $\alpha \in \cA$
the set $\rho (\alpha (Z))$ is a closed curve in $Q$ and the restriction of $\rho$ yields an isomorphism between
$\alpha (Z)$ and $\rho (\alpha (Z))$.
\ecor

\bproof By Theorem \ref{gp1.t1} (iv) there is a family $\cA$ such that for a general $\alpha \in \cA$ the restriction of $\rho$  yields an injective map of the Zariski tangent bundle of 
$\alpha (Z)$ into the Zariski tangent bundle of $\overline{\rho (\alpha (Z))}$. Note also that,
 being an injective image of $Z$,
$\rho (\alpha (Z))$ is closed in $Q$
since  any once-punctured curve in a quasi-affine algebraic variety is automatically closed.
Therefore, the bijective morphism $\rho|_{\alpha (Z)} : \alpha (Z) \to \rho (\alpha (Z))$ is proper. 
This yields the finiteness of  $\rho|_{\alpha (Z)}$.
Furthermore, since the map $T_z \alpha (Z) \to T_qQ$ is injective,
\cite[Proposition 7]{Ka91}) implies that  $\rho|_{\alpha (Z)} : \alpha (Z) \to \rho (\alpha (Z))$ is an isomorphism
which is the desired conclusion.
\eproof

We shall need later the following fact.

\bprop\label{gp1.p1} Let  $X$ be a $G$-flexible variety for a group $G \subset \SAut (X)$ generated by a saturated set of locally nilpotent vector fields,
$\rho : X \to Q$ be a dominant morphism into another algebraic variety $Q$ and
$Z$ be a locally closed reduced subvariety in $X$. Suppose that $F$ is a finite subset of $Z$ such that
$\dim T_{z_0}Z \leq \dim Q$ for every $z_0\in F$.
Then  one can find an automorphism $\alpha_0\in G$ with the following property:
for any irreducible algebraic family $\cA\subset G$ of automorphisms of $X$ containing $\alpha_0$
and any general element $\alpha \in \cA$ there is a neighborhood $V_\alpha$ of $\alpha (F)$ in $\alpha (Z)$ such that for every $z \in V_\alpha$ and $q=\rho (z)$ the induced map $\rho_* : T_z\alpha (Z) \to T_qQ$
of the tangent spaces is injective.
\eprop

\bproof  
By flexibility and Theorem \ref{fm.t2}  there is an automorphism $\alpha_0$ of $X$ such that
for every $z_0 \in F$ and 
$z_0' =\alpha_0 (z_0)$ the variety $Q$ (resp. the morphism $\rho$) is smooth at (resp. over) $q_0'=\rho (z_0')$ and
the induced map $T_{z_0'}\alpha_0 (Z) \to T_{q_0'}Q$
of the tangent spaces is injective. In the case when $\alpha_0$ is contained in an irreducible family $\cA$ of automorphisms
put $z_1=\alpha (z_0)$, $Z_1=\alpha (Z)$ and $q_1=\rho (z_1)$ for a general $\alpha \in \cA$. Then by continuity
the variety $Q$  (resp. the morphism $\rho$) is smooth at (resp. over) $q_1$
and  the induced map  $\rho_*|_{T_{z_1}Z_1}: T_{z_1}Z_1 \to T_{q_1}Q$ of the tangent spaces is also injective.
That is, the kernel  of $ \rho_*|_{T_{z_1}Z_1}$ is zero. Hence for every $z$ in some neighborhood of $z_1$ in $Z_1$ the kernel of
the induced map $T_{z}Z_1 \to T_{\rho(z)}Q$ is also  zero. 
This yields the desired conclusion.
\eproof

We shall also need a parametric version of Theorem \ref{gp1.t1} (iii).

\bprop\label{gp1.p3} Let the assumptions and notations of Theorem \ref{gp1.t1} (iii) hold with $Q_0=Q$ and let $U$ be a smooth algebraic variety.
Let $\cZ =U \times Z$, $\cP=U\times P$, $\cX = U \times X$, $\cQ = U \times Q$, $\check \kappa =({\rm id}, \kappa): \cX \to \cP$
and $\check \tau =({\rm id}, \tau) : \cQ \to \cP$. Suppose that 
$\check \rho : \cX \to \cQ$ is a smooth morphism  such that $\check \kappa = \check \tau \circ \check \rho$.  Consider the $G$-action 
on $\cX$ such that the natural projection $\cX \to X$ is $G$-equivariant,

{\rm (1)} Then there exists an algebraic family $\cA\subset G $ of automorphisms of $\cX$ such that
for a general element $\alpha \in \cA$, a general $u\in U$ and
every $z \in \alpha (u\times Z)$ the induced map $\check \rho_* : T_z\alpha (u\times Z) \to T_{\check \rho (z)} (u\times Q)$ is injective.

{\rm (2)}  Suppose that the assumptions of Theorem \ref{gp1.t1} (iv) are satisfied as well. Then the morphism 
$\check \rho|_{\alpha (u\times Z)} : \alpha (u\times Z) \to u\times Q$ in (1) is also injective for general $\alpha \in \cA$ and $u \in U$.

{\rm (3)} Furthermore, suppose that $W$ is a subvariety of $X$ for which $\dim W + \dim Z< \dim X$ (resp. $\dim W + \dim Z\leq  \dim X$)
and $\cW =U \times W$. Then one can suppose additionally that $\cW \cap \cX_u$ does not meet $\alpha (\cZ_u)$ (resp. 
$\cW \cap \cX_u\cap \alpha (\cZ_u)$ is finite) for a general $u \in U$.
\eprop

\bproof  Let  $T(X/P)^*$ , $T(Z/P)^*$ and $Y^*$ be as in the proof of Theorem \ref{gp1.t1} (iii). 
Let  $\cT^*= U \times T(X/P)^*$, $\cY^*=U \times Y^*$ and $\cE=U\times T (Z/P)^*$.
Then $G$ acts on $\cT^*$
and by Theorem \ref{fm.t2} this action is transitive on every fiber of $\cT^* \to \cX \overset{\check \kappa}{\to} \cP$.
Note that $\dim \cE^*=\dim T(Z/P) + \dim U \leq \dim Q- \dim P+\dim U$ (because 
$ \dim T(Z/P)\leq \dim Q- \dim P$ by the assumption). 
Since all fibers of the natural projection $\cY^*\to \cP$ are of the same dimension,
 by Remark \ref{agga.r1} (1) one has  $\dim \cY^* \times_\cP \cE^*= \dim \cY^* + \dim \cE^* -\dim \cP$.
 Since the codimension of $\cY^*$ in $\cT^*$ is equal to $\dim Q-\dim P$, by Theorem \ref{agga.t2} 
 we can choose one-parameter unipotent algebraic subgroups $\tH_1, \ldots, \tH_{\tilde m}$ of $G$ such that 
for a general element $\alpha =(\tilde h_1, \ldots , \tilde h_{\tilde m}) \in \tH_1 \times \cdots \times \tH_{\tilde m}$ and 
$\cE''=\alpha_* (\cE^*)$ one has $\dim \cY^* \cap \cE'' \leq \dim U$. Hence for a general $u \in U$ the fiber  $E_u$ of  $\cE''\cap \cY^*$ over $u$ is at most finite.
Note that  up to the zero vector the set $(T\alpha (u\times Z)) \cap \cY^*$ is 
the kernel of the natural projection $T\alpha (u\times Z) \to T(u \times Q)$. This implies that $E_u$ is empty and, therefore,  the map
 $(T\alpha (u\times Z))\cap \cE^* \to T(u \times Q)$ is injective. Consequently, this yields (1).

The similar argument works in the case of Theorem \ref{gp1.t1} (i) with $\dim Q \geq 2\dim Z +1$. That is,
$\check \rho|_{\alpha' (u\times Z)} : \alpha' (u\times Z) \to u\times Q$ is injective for  a general $u \in U$ and
a general automorphism $\alpha'$ in some algebraic family $\cA'\subset G$.
By Remark \ref{gp1.r1} (2) one can suppose that the induced map of the tangent bundles is still injective and we have (2).

For (3) note that Theorem \ref{agga.t2} implies that there is a family $\cA'\subset G$ of algebraic automorphisms of $\cX$ such that for a general $\alpha' \in \cA'$ one
has $\dim \cW \cap \alpha' (\cZ) <\dim U$ (resp. $\dim \cW \cap \alpha' (\cZ) \leq \dim U$). 
Hence $\cW \cap \cX_u$ does not meet $\alpha (\cZ_u)$ (resp. 
$\cW \cap \cX_u\cap \alpha (\cZ_u)$ is finite) for a general $u \in U$. Remark \ref{gp1.r1} (2) implies that one can suppose that $\cA =\cA'$
which yields the desired conclusion.
\eproof. 

\brem\label{gp1.r3} By Remark \ref{gp1.r1}  (4) one can drop the assumption that $\kappa$ is dominant in the first statement of 
Propostion \ref{gp1.p3}. Similarly, this assumption can be dropped in the second statement if one replaces
the inequality  $\dim Z\times_P Z\leq 2\dim Z- \dim P$ and $\dim Q \geq \dim Z+m$
by the ones in Formula (\ref{gp1.eq1}).

\erem

\section{General projections for flexible varieties. II}

\bnota\label{gp2.n1} In this section $X$ and $P$ are smooth algebraic varieties,
$P$ is a closed affine subvariety of $\A_\bk^m$, $\kappa : X\to P$ is a  surjective morphism, and $Z$ is a closed subvariety of $X$.
We suppose that $\hat P$ is a completion of $P$, $\hX$ is a completion of $X$, $\hD=\hX \setminus X$, 
and $\hat \kappa : \hX \dashrightarrow  \hat P$ is a rational extension of $\kappa$. 
Let $\hat f$ be the rational extension of a regular function $f\in \bk [X]$ to $\hX$.
Denote by $\cR (\hat f)$ the subvariety of  $\hD$ that consists of points which are either indeterminacy points of $\hat f$ 
or points at which $\hat f$ is regular and takes finite values.
Given any morphism  of the form $\lambda = (f_1, \ldots , f_N) : X \to \A_\bk^N $
we let $\cR(\hat \lambda) = \bigcap_{i=1}^N \cR(\hat f_i)$ where 
$\hat \lambda$ is the extension of $\lambda$ to $\hX$. In particular, since
$\kappa$ is given by $m$ coordinate functions we can define $\cR (\hat \kappa )$. 
\enota

The aim of this section is to describe some conditions under which  
the morphism $\kappa|_Z : Z\to P$ is proper and, in particular, $\kappa (Z)$ is closed in $P$. 

\bprop\label{gp2.p1} 
Let  $\brZ$ be
the intersection of $\hD$ with the closure of $Z$ in $\hX$.  If
 $\brZ\cap \cR (\hat \kappa )=\emptyset$ then  $\kappa|_Z : Z\to P$ is proper.
\eprop

\bproof Using a resolution $\pi : \bX \to \hX$ of the indeterminacy points of $\hat \kappa$ one gets a morphism
$\bar \kappa = \hat \kappa \circ \pi : \bX \to \hat P$.
Treat $X$ as a subvariety of $\bX$ and denote by $\bZ$ the closure of $Z$ in $\bX$. Note that $\kappa|_Z : Z \to P$ is 
proper if and only if $(\bar \kappa^{-1}(P) \setminus X)\cap \bZ =\emptyset$. Note also that $\pi (\bar \kappa^{-1}(P) \setminus X)$ consists of
all point in $\hD$ where either $\hat \kappa$ has indeterminacy or where $\hat \kappa$ is regular and takes  values in $P$.
Since $P$ is closed in $\A_\bk^m$ the latter means that all coordinate functions of $\hat \kappa$ take finite values and, therefore,
 $\pi (\bar \kappa^{-1}(P) \setminus X)$ is contained in $\cR (\hat \kappa )$. Furthermore, if $\hZ$ is the closure of $Z$ in $\hX$ the
 properness of $\pi$ implies that $(\bar \kappa^{-1}(P) \setminus X)\cap \bZ =\emptyset$ if and only if 
 $\pi (\bar \kappa^{-1}(P) \setminus X)\cap \hZ =\emptyset$.  Since $\pi (\bar \kappa^{-1}(P) \setminus X)\cap \hZ=\brZ\cap \cR (\hat \kappa )$
 we have the desired conclusion. \eproof

 \brem\label{gp2.r1r} One may have properness of  $\kappa|_Z : Z\to P$ even if
$\brZ\cap \cR (\hat \kappa )\ne \emptyset$. Indeed, consider $X=\A_\bk^2 \subset \PP_\bk^2=\hX$, $P=\A_\bk^1$
and the map $X \to P$ given by $(x,y) \to x$ for a coordinates system $(x,y)$ on $\A_\bk^2$.
Let $Z\subset X$ be the parabola given by $y=x^2$. Then $\kappa|_Z : Z\to P$ is proper but $\bZ$ meets $\hD$
at an ideteminacy point of $\hat x$.

\erem

\bcor\label{gp2.c1} Let  $\hD$ be irreducible and $I(\hat f_i)$ be the set of indeterminacy points of $\hat f_i$
for $\kappa = (f_1, \ldots, f_m)$. Suppose that every $f_i$ is non-constant on $Z$ and 
$\brZ \cap \bigcap_{i=1}^m I(\hat f_i) =\emptyset$. Then  the morphism
$\kappa|_Z : Z \to P$ is proper and $\kappa (Z)$ is closed in $P$. 

\ecor

\bproof  
Let $x \in \brZ$.
By the assumption at least for one index $i$ the rational function $\hat f_i$ yields  a regular morphism into $\PP^1$ in a neighborhood of $x$.
Note the value
of $\hat f_i$ at $x$ is $\infty$ (indeed, $\hat f_i^{-1}(\infty) \subset \hD$ and, therefore, $\hat f_i^{-1}(\infty)$ must be equal to $\hD$ since the latter is irreducible). 
Thus $\brZ \cap \cR (\hat \kappa )$ is empty and Proposition \ref{gp2.p1} yields the desired conclusion.
\eproof


\bthm\label{gp2.t1}  Let $\hat \kappa : \hX \to \hat P$ be regular
and  $\tilde \kappa : \tX \to P$ be its restriction to $\tX=\hat \kappa^{-1} (P)$.
Suppose that $H$  is an algebraically generated group acting on $\tX$ so that
 it transforms every fiber of $\kappa$ into itself.
Let $\tD=\hD \cap \tX$ be a finite disjoint union of 
irreducible smooth subvarieties $\{ T_i \}$ such that every set $ \tilde \kappa (T_i)$ is a smooth
variety and
$\tilde \kappa|_{T_i}: T_i \to \tilde \kappa (T_i)$ is a smooth morphism with each fiber being 
an orbit of $H$. 
Let $\theta : X \to P\times \A_\bk^N$ be a morphism of the form $\theta = (\kappa, \lambda)$
and let $S_i=\cR (\hat \lambda ) \cap T_i$ and $Z_i=\brZ \cap T_i$ where $\brZ$ from Proposition \ref{gp2.p1}.
Suppose that for  
every $i$ either $S_i$ is empty or
$\dim Z_i\times_P S_i + \dim \tilde \kappa (T_i) < \dim T_i$. 
Then there exists an algebraic fammily $\cA \subset H$ of automorphisms  of $X$ such that
for a general element $h \in \cA$ the morphism $\theta|_{h(Z)} : h(Z)\to P\times \A_\bk^N$ is proper.

\ethm
\bproof  
Since $P$ is closed in $\A_\bk^n$ the description of $\tX$ implies that for every point  $x \in \hD \setminus \tD$
 there is a coordinate function $f$ of $\kappa$ for which $\hat f$
has a regular value $\infty$ at $x$.
Thus $(\hD \setminus \tD) \cap \cR (\hat \kappa) = \emptyset$ and $\brZ \cap \cR (\hat \theta)
\subset \brZ \cap \cR (\hat \lambda)\cap \tD$. 
Applying Theorem \ref{agga.t2} to the morphism 
$\tilde \kappa|_{T_i}: T_i \to \tilde \kappa (T_i)$ in the case of a non-empty $S_i$
we see that 
for a general element $h$ of some algebraic family $\cA \subset H$ the variety $h(Z_i)\cap S_i$ is  empty
when
$\dim Z_i\times_P S_i + \dim \tilde \kappa (T_i) < \dim T_i$.  
Hence, by Proposition \ref{agga.p2}  $ h(\brZ) \cap \cR (\hat \lambda)\cap \tD= \emptyset$  and we are done by Proposition \ref{gp2.p1}.
\eproof

\bcor\label{gp2.c1'} Let $\tilde \kappa : \tX \to P$, 
$H$, $T_i$, $\theta=(\kappa, \lambda)$ and $S_i$ be as in Theorem \ref{gp2.t1}.
Suppose that for every $i$ either $S_i$ is empty or $\tilde \kappa|_{S_i} : S_i \to \tilde \kappa (T_i)$ is a flat morphism
with fibers of dimension  $l_i$. Let $k=\dim Z$ and $m_i$ be the dimension of the fibers
of $\tilde \kappa|_{T_i}: T_i \to \tilde \kappa (T_i)$. Suppose that
for every $i$   one has $l_i+k-1 <m_i$. 
Then for a general element $h \in H$ the morphism $\theta|_{h(Z)} : h(Z)\to P\times \A_\bk^N$ is proper.
\ecor

\bproof  Let $Z_i$ be as in Theorem \ref{gp2.t1}. By the same theorem we can assume that $S_i$ is not empty.
Since  $ \tilde \kappa (S_i) = \tilde \kappa (T_i)$ Remark \ref{agga.r1} (1) implies that
$Z_i \times_P S_i =\dim S_i + \dim Z_i -\dim  \tilde \kappa (T_i)$, i.e., the desired inequality
from Theorem \ref{gp2.t1} can be rewritten as  $\dim Z_i + \dim S_i  < \dim T_i$.
By the assumption we have $l_i+\dim \tilde \kappa (T_i)+k-1 <m_i+\dim \tilde \kappa (T_i)$.
Note that $\dim S_i=l_i+\dim \tilde \kappa (T_i)$, $\dim T_i=m_i+\dim \tilde \kappa (T_i)$,
and $\dim \brZ=k-1$. Hence we have this inequality $\dim S_i + \dim \brZ<\dim T_i$ which concludes the proof.
\eproof

\bnota\label{gp2.n2}
From now on by $\ED (X)$ we denote $\ED (X) =\max (2\dim X +1, \dim TX)$. 
\enota

\bcor\label{gp2.c2} Let the assumptions of Theorem \ref{gp2.t1} or Corollary \ref{gp2.c1'} hold with $\theta$ being a dominant morphism.
Suppose also that $X$ is a $G$-flexible variety where
$G \subset \SAut (X)$ is generated by a saturated set of locally nilpotent derivations
 and that  $\ED (Z)\leq \dim P+N $.
Then  there exists an algebraic family $\cA \subset G$ of automorphisms such that for a general $\alpha \in \cA$ 
the morphism $\theta|_{\alpha(Z)} : \alpha (Z)\to P\times \A_\bk^N$ is a closed embedding.
\ecor

\bproof  Consider Theorem \ref{gp1.t1} (iv) with $P$ in {\em  its formulation} being a singleton unlike $P$ in the formulation above.  
Applying this special case one can see that there exists an algebraic family $\cA_0 \subset G$ of automorphisms such that for a general $\beta \in \cA_0$ the morphism 
 $\theta|_{\beta(Z)} : \beta (Z)\to P\times \A_\bk^N$ is an embedding. By Proposition \ref{agga.p2} this property remains valid
if one replaces $\cA_0$ by $\cA = H \times \cA$.
 Hence by Theorem \ref{gp2.t1} for a general $h \in H$ and $\alpha = h \circ \beta \in \cA$ the morphism  $\theta|_{\alpha(Z)} : \alpha (Z)\to P\times \A_\bk^N$ is
 also proper. Repeating now the argument from Corollary \ref{gp1.c1} we conclude that
  it is a closed embedding which is the desired conclusion.
\eproof

\brem\label{gp2.r2}  Suppose that unlike in Notation \ref{gp2.n1} $P$ is not an affine variety
but only a quasi-affine one and let $R$ be an affine variety containing $P$ as an open subset.
Assume additionally that $\kappa (Z)$ is closed in $R$. Then the conclusion about properness of 
$\theta|_{h(Z)} : h(Z) \to P\times \A_\bk^N$ in Theorem \ref{gp2.t1} remains valid.
For a more general setting we need the following.
\erem

\bnota\label{gp2.n3}
 Let $ \rho : \cX \to \cP$ be an affine morphism of smooth quasi-affine varieties
and $\cP=\bigcup_{j=1}^mP_j$, where each $P_j$ is an affine Zariski dense open subset of $\cP$.
Let $X_j=\rho_i^{-1}(P_j)$, $\rho_j=\rho|_{X_j}: X_j \to P_j$ and 
$\tilde \rho_j : \tX_j \to P_j$ be an extension of $\rho_j$ to a proper morphism, where $\tX_j$ is a smooth variety
and $\tD_j=\tX_j \setminus X_j$.  Suppose that for every $ j=1, \ldots , m$ there is a set $\cN_j$ of locally nilpotent vector fields on $X_j$
tangent to the fibers of $\rho_j$ and extendable to complete vector fields on $\tX_j$. Let $H_j$ be the group of automorphisms of $X_j$ (and $\tX_j$)
generated by the elements of the flows of the vector fields from $\cN_j$.
Let  $\lambda : \cX \to \A^N_\bk$ be a morphism, $\Theta=(\rho, \lambda): \cX \to \cP\times \A_\bk^N$,
$\cZ$ be a closed subvariety of $\cX$ and $\cZ_j=\cZ \cap X_j$. Suppose also that for every $j=1, \ldots, m$
the assumptions (and, therefore, the conclusion) of Theorem \ref{gp2.t1} are true if $\kappa : X \to P$, $\tilde \kappa :\tX \to P$, $Z$, $\tD$, $H$  and $\theta$ in
the formulation of that theorem are replaced with $\rho_j: X_j \to P_j$, $\tilde \rho_j : \tX_j \to P_j$, $\cZ_j$ $\tD_j$, $H_j$ and $\theta_j = \Theta|_{X_j}$.
\enota

\bthm\label{gp2.t2}  Let Notation \ref{gp2.n3} hold. 
Then there exists an algebraic family $ \cA$ of automorphisms of $\cX$ over $\cP$ such that 
for a general element $h \in \cA$ the morphism $\Theta|_{h(\cZ)} : h(\cZ)\to P\times \A_\bk^N$ is proper.
\ethm

\bproof  
Let $I_j\subset \bk [\cP]$ be the defining ideal of the variety $\cP \setminus P_j$.
For every $\delta \in \cN_j$ one can find $f \in I_j$ such that $f\delta$ is a locally nilpotent vector field on $\cX$ tangent to
the fibers of $\rho$. Furthermore, for a given point $p \in P$ such $f$ can be chosen so that $f (p)=1$.
Consider the set  $\cN_j'$ of all locally nilpotent vector fields of this form $f\delta$.
The elements of their flows generate a group $H_j'\subset \SAut (\cX/\cP)$.
Without loss of generality we can suppose that for every $\sigma \in \cN_j$ and each $r \in \bk$ one has $r\sigma \in \cN_j$.
Then for every point $p \in P_j$ the restrictions of $\cN_j$ and $\cN_j'$ to $\rho_j^{-1}(p)$ (resp. $\tilde \rho_j^{-1} (p)$) coincide.
Hence, the restriction of  the $H_j'$-action to $\rho_j^{-1}(p)$ (resp. $\tilde \rho_j^{-1} (p)$) is the same as the $H_j$-action.
Thus, we can replace $H_j$ in Notation \ref{gp2.n3} by $H_j'$.
In particular, from the beginning we can assume that $H_j \subset \SAut(\cX/\cP)$. By Theorem \ref{gp2.t1} there exists an algebraic family
$\cA_j \subset H_j \subset \SAut(\cX/\cP)$ such that for a general element $h \in \cA_j$ the morphism $\Theta|_{h(\cZ_j)} : h(\cZ_j)\to P_j\times \A_\bk^N$ is proper.
Arguing as in Remark \ref{gp1.r1} (2) we can find an algebraic family $\cA \subset \SAut(\cX/\cP)$ such that for a general element $\alpha \in \cA$ 
the morphism $\Theta|_{\alpha (\cZ_j)} : \alpha (\cZ_j)\to P_j\times \A_\bk^N$ is proper for every $j$ which yields the desired conclusion.
\eproof

\section{General projections for partial quotient morphisms of flexible varieties}\label{gp3}

In the case of partial quotients of $\G_a$-actions one can get local properness under much milder assumptions than in Theorem \ref{gp2.t1}. It is reflected in the following fact
which plays an important role in \cite{KaKuTr}.

\bthm\label{gp3.t1} Let $X$ be a smooth quasi-affine algebraic variety, $\cN$ be a saturated set of locally nilpotent vector fields on $X$, 
and $G\subset \SAut (X)$ be the group generated by $\cN$. Suppose also that $X$ is $G$-flexible.  
Let $\rho_0 : X \to Q$ be a partial quotient morphism
associated with a nontrivial $\delta_0 \in \cN$, $Z$ be a
locally closed reduced subvariety of $X$
of codimension at least 2 and $F$ be a finite subset of $Z$ such that $\dim T_{z_0}Z \leq \dim Q$ 
for every $z_0 \in F$.
Then there exists a connected algebraic family  $\cA\subset G$ of automorphisms 
such that for a general element $\alpha \in \cA$
and  the closure $\bZ_\alpha'$ of $Z_\alpha'=\rho_0\circ \alpha (Z)$ in $Q$ 
one can find a neighborhood $V_0'$ of $\rho_0 (\alpha (F))$ in $\bZ_\alpha'$ such that for $V_0=\rho_0^{-1}(V_0')\cap \alpha (Z)$ the morphism $\rho_0|_{V_0} : V_0 \to V_0'$ 
is an isomorphism.

\ethm

The proof consists mostly of  reminding some results from \cite{FKZ-GW}.

\bprop\label{1.14} {\rm (Proposition 2.15 in  \cite{FKZ-GW})} Let $X$, $G$, and $\cN$ be as in  Theorem \ref{gp3.t1}.
Then for any locally nilpotent derivation
$\de_0\in \cN$ one
can find another one
 $\de_1\in \cN$ such that the subgroup $H \subset G$
generated by $\delta_0$, $\delta_1$ and all their replicas acts with an open orbit on $X$.
\eprop

\brem\label{gp3.r1} In fact we have more. It follows from the proof of  \cite[Proposition 2.15]{FKZ-GW}) that $\delta_1$ can be chosen so that the open orbit of $H$ contains a given finite subset of
$X$.

\erem

\bnota\label{44a} (a) 
Let $\delta_0$ and $\delta_1$ be as in Proposition \ref{1.14} and $U^i$  be the one-parameter unipotent subgroup of $\SAut (X)$ associated with  $\delta_i$. 
Any function $f\in \ker\de_0\backslash \ker \de_1$ yields the one-parameter group $U^0_f$ associated with the replica
$f \de_0$, and similarly
$g\in \ker\de_1\backslash \ker \de_0$ yields the one-parameter group  $U^1_g$  associated with the replica
$g\de_1$.

(b) To any sequence of invariant functions
\be\label{seq}
\cF=\{f_1,\ldots,f_s, g_1, \ldots,g_s\},\,\,\,\,\mbox{where}
\,\,\,\, f_i\in\ker\de_1\backslash\ker
\de_0\,\,\,\,\mbox{and}\,\,\,\, g_i\in\ker\de_0\backslash\ker
\de_1\,,
\ee
we associate an algebraic family
of automorphisms defined by the product
\be\label{00121}
U^\cF=U^1_{f_s}\cdot
U^0_{g_s}\cdot\ldots\cdot U^1_{f_1}\cdot
U^0_{g_1}\subseteq H\,.
\ee
More generally, given a tuple
$\kappa=(k_i,l_i)_{i=1,\ldots,s}\in\N^{2s}$ the product
\be\label{001210}
U_\kappa=U_\kappa^\cF=
U^1_{f_s^{k_s}}\cdot
U^0_{g_s^{l_s}}\cdot\ldots\cdot U^1_{f_1^{k_1}}\cdot
U^0_{g_1^{l_1}}
\subseteq H\,
\ee
is as well an algebraic family of automorphisms.
\enota

\bprop\label{cor-0012}   {\rm (Corollary 5.4 in  \cite{FKZ-GW})}There is a finite collection of invariant functions $\cF$ as in (\ref{seq})
such that for any sequence
$\kappa=(k_i,l_i)_{i=1,\ldots,s}\in\N^{2s}$ the algebraic family of automorphisms
$U_\kappa$ as in (\ref{001210})
has an open orbit  $O(U_\kappa)$ that coincides with $O(H)$ and so
does not depend on the choice of $\kappa\in\N^{2s}$.\eprop

\bnota\label{4.5}
We keep the notation and assumptions from \ref{44a}(a).

(a) Let $\rho_0:X\to Q_0$ and $\rho_1:X\to Q_1$ be partial quotient morphisms with respect
to the unipotent subgroups $U^0$ and $U^1$, respectively.
It is proven in \cite[Lemma 3.3]{FKZ-GW} that there are open embeddings $X\hto \bX$, $Q_0\hto\bQ_0$, and $Q_1\hto\bQ_1$
into normal projective varieties such that the following conditions
are satisfied.
\bnum[(i)]
\item  $\rho_0$ and $\rho_1$ extend to morphisms
$\brho_0:\bX\to\bQ_0$ and $\brho_1:\bX\to \bQ_1$.

\item the unique ``horizontal" divisors $D_0\subset \bX \setminus X$  and $D_1\subset \bX \setminus X$, that 
are mapped  birationally  (by $\brho_0$ and $\brho_1$) onto
$\bQ_0$ and $\bQ_1$ respectively,  are smooth.

\item the completion $\bX$ satisfies some other conditions which we assume to be true but omit because they are not needed for the formulation of
Proposition \ref{prop-two} below.

\enum
(b) Given a closed subscheme $Y\subseteq X$ of codimension at least $2$ we call
$$
\p_0Y=\bY\cap D_0\and \p_1Y=\bY\cap D_1
$$
the {\em partial boundaries.}

(c) 
For a one-parameter group $U$ we let $U^*=U\backslash\{\id\}$
and for $U_\kappa=
U^1_{f_s^{k_s}}\cdot
U^0_{g_s^{l_s}}\cdot\ldots\cdot U^1_{f_1^{k_1}}\cdot
U^0_{g_1^{l_1}}$ as in \eqref{001210}  we let
$$
U_\kappa^*=
U^{1*}_{f_s^{k_s}}\cdot
U^{0*}_{g_s^{l_s}}\cdot\ldots\cdot U^{1*}_{f_1^{k_1}}\cdot
U^{0*}_{g_1^{l_1}}\,.
$$

 \enota

\bprop\label{prop-two} {\rm (Proposition 5.11 in  \cite{FKZ-GW})}
Let $(Y_\alpha)_{\alpha\in A }$ be a flat family of proper closed subsets of $X$.
Assume that the partial boundaries
$\p_i Y_\alpha$ are contained in $E_{\alpha,i}$,
where the $(E_{\alpha,i})_{\alpha\in A}$, $i=0,1$,
form flat families of proper closed subsets of $D_i$.
Then one can find an open dense subset $A^o$ of $A$, flat families of proper, closed subsets
$(E^o_{\alpha, i})_{\alpha\in A^o}$ of $ D_i$ ($i=0,1$), and a sequence
$\kappa=(k_1,l_1,\ldots,k_s,l_s)\in\N^{2s}$ such that for any element $h\in U_\kappa^*$ we have
$$
\p_i (h.Y_\alpha)\subseteq E^o_{\alpha ,i}\, ,\qquad i=0,1\,,\; \forall\,\alpha\in A^o\,.
$$
\eprop

\bproof[Proof of Theorem \ref{gp3.t1}]  Let $Q_0$ be a smooth Zariski open dense subset of $Q$ such that
for $X_0=\rho_0^{-1}(Q_0)$ the morphism $\rho_0|_{X_0} : X_0 \to Q_0$ is smooth. 
Since the $G$-action on $X$ is $m$-transitive for every $m$ replacing $Z$ by $g(Z)$ for some $g \in G$
we can suppose that $F \subset X_0 \cap O(H)$ where $O(H)$ is the open orbit
of $H$ from Proposition \ref{1.14}.
Let $\bZ$ be the closure of $Z$ in $X$. 
By Theorem \ref{gp1.t1}  (i), (ii) and Remark \ref{gp1.r1 (4)} there exists an irreducible family $\cA \subset G$ of algebraic
automorphisms of $X$ such that for a general $\alpha \in \cA$ one has 

(a) $\rho_0|_{\alpha (\bZ)} : \alpha (\bZ)\cap X_0 \to \rho_0 \circ \alpha (\bZ)$ is birational and
$\rho_0^{-1}(\rho_0 (\alpha (F))\cap \alpha (\bZ) =\alpha (F)$ (here we use the fact that 
$\alpha (F) \subset X_0\cap O (H)$ since $F \subset X_0\cap O(H)$).

Recall that by Remark \ref{gp1.r1} (1) we can suppose that $\cA$ is irreducible and it contains the identity map $e$. 
By  the definition of $\SAut (X)$ any automorphism $\alpha_0 \in \SAut (X)$ is an element of an algebraic family
$H_1 \times \ldots \times H_m$  of automorphisms on $X$ as in Theorem \ref{agga.t1} where every $H_i$ is a unipotent one-parameter group. Using  the way to enlarge $\cA$ as in Remark \ref{gp1.r1} (2) we can replace $\cA$ by $H_1 \times \ldots \times H_m\times \cA$ where the last family is still irreducible and contains $e$
but it contains also $\alpha_0$ now.
Hence,  since $\dim T_{z_0}Z \leq \dim TQ$ for every $z_0\in F$,
choosing $\alpha_0$ as in Proposition \ref{gp1.p1}
we also have

(b) a neigborhood $\tV_\alpha$ of $\alpha (F)$ in $\alpha (Z)$ such that for every $z \in \tV_\alpha$ the induced map $T_z\alpha (Z) \to T_{\rho_0 (z)} Q$
of the tangent spaces is injective.

Let  $U_\kappa$ be from Proposition \ref{prop-two}
and let $\beta = h \circ \alpha$ be a general element of the family $U_\kappa \cdot \cA$.
By Proposition \ref{agga.p2} and Remark \ref{gp1.r1} (1) we still have conditions (a) and  (b) for this bigger family.
By Proposition \ref{prop-two} the partial boundary $\p_0 ( \beta (\bZ))=\p_0 (h.  \alpha (\bZ))\subseteq  E^o_{\alpha ,0}$ where
$E^o_{\alpha ,0}$ is a proper subvariety of $D_0$ from Notation \ref{4.5}. 
This implies that  the morphism $\rho_0|_{ \beta (\bZ)} : \beta (\bZ) \to Q$ is proper over
$Q\setminus \brho_0 ( E^o_{\alpha ,0})$ where $\brho_0$ is again  from Notation \ref{4.5}.
On the other hand, by Proposition \ref{cor-0012} 
$\beta (z_0) =h. \alpha (z_0)$ runs over the open set $O(H)$ when $h$ runs over $U_\kappa$. Hence 
$\rho_0 \circ \beta (F)$ does not meet $ \brho_0 ( E^o_{\alpha ,0})$ for a general $\beta$. Therefore,
the morphism $\rho_0|_{ \beta (\bZ)} : \beta (\bZ) \to Q$  is proper over an neighborhood of 
$F'=\rho_0 \circ \beta (F)$ in $Q$. 

Assume that there is an irreducible subvariety $Y' \subset \rho_0\circ \beta (Z) \cap Q_0$ such 
that $\dim Y'< \dim Y$ where $Y$ is an irreducible component of $\rho_0^{-1} (Y') \cap \beta (Z)$ with 
some $z_0'=\rho_0(\beta (z_0))$ from $F'$ contained in the closure of $\rho_0 (Y)$. 
Then by the Chevalley's theorem \cite[13.1.3]{EGA} $\dim \rho_0^{-1}(z_0')\cap Y\geq 1$ (this fiber cannot be empty because of properness) contrary to the fact that by (a)
$\rho_0^{-1}(z_0') \cap \beta (Z)=\beta (z_0)$. 
Hence we can suppose that $\rho_0|_{\beta (Z)}$ is not only birational but also quasi-finite over a neighborhood
$V_0' \subset \rho_0\circ \beta (Z)$ of $F'$ which is combination with properness implies finiteness by the Grothendieck theorem \cite[Theorem 8.11.1]{EGA}.

If for a sufficiently small $V_0'$ and $V_0=\rho_0^{-1} (V_0') \cap \beta (Z)$ the morphism $\rho_0|_{V_0}$ is not injective  
then for some curve $C' \in V_0$ through a point $z_0'\in F'$ and $C=\rho_0^{-1}(C') \cap \beta (Z)$
the finite morphism $\rho_0|_C : C \to C'$ must be ramified at $\beta (z_0)$.  However, this is contrary to the fact that the induced map $T_{\beta (z_0)} T( \beta (Z)) \to T_{z_0'} Q$
of the tangent spaces is injective by (b). That is, we can suppose that $\rho_0 |_{V_0}: V_0\to V_0'$ is injective and proper. Since it is also finite condition (b)   implies that
this map is an embedding (e.g., see \cite[Proposition 7]{Ka91}) which yields the desired conclusion.
\eproof

\brem\label{gp3.r2}
(1) Note that by Remark \ref{gp1.r1} (1) and  construction the family $\cA$ from  Theorem \ref{gp3.t1} is a Zariski open subset in a larger family of automorphisms which contains the identity map.

(2) Similar to Remark \ref{gp1.r1} (2) we also note that the family $\cA$  from Theorem \ref{gp3.t1} does not depend on the choice of the subvariety $Z$.

(3) Let $X$, $Z$, $F$, $\cN$ and $G$ be as in Theorem \ref{gp3.t1}, $S =\{ \delta_1, \ldots, \delta_s\}\subset \cN$
and $\rho_i : X \to Q_i$ be a partial quotient morphism associated with $\delta_i$ for $i=1, \ldots s$.  Then  an easy adjustment of
 the proof yields the following generalization of Theorem \ref{gp3.t1}.

{\em There exists a connected algebraic family  $\cA\subset G$ of automorphisms such that for a general element $\alpha \in \cA$, 
every $i=1, \ldots , s$, and  the closure $Z_\alpha^i:=\overline{\rho_i \circ \alpha (Z)}$ of $\rho_i\circ \alpha (Z)$ in $Q_i$ 
one can find a neighborhood $V_i'$ of $\rho_i (\alpha (F))$ in $Z_\alpha^i$ such that for $V_i=\rho_i^{-1}(V_i')\cap \alpha (Z)$ the morphism 
$\rho_i|_{V_i} : V_i \to V_i'$  is an isomorphism.}

\erem

\section{The case of Gromov-Winkelmann flexible varieties}

The aim of this section is the following fact.

\bthm\label{cgw.t1} Let $Z$, $Y_1$, and $Y_2$ be closed subvarieties of $\A_\bk^n$  such that $Y_1 \cap Z=Y_2 \cap Z =\emptyset$,  $\dim Z\leq n-3$
and $\ED (Y_1)\leq n-2$ (where $\ED (Y_1)$ is as in Notation \ref{gp2.n2}). Let $\varphi : Y_1 \to Y_2$ be an isomorphism and $X =\A_\bk^n \setminus Z$.
Suppose also that either

{\rm (a)} $\dim Z + \dim Y_1 \leq n-3$, or

{\rm (b)} $\dim Y_1 =1$ and  $\dim Z = n-3$.

Then there exists an automorphism $\gamma \in \SAut (X)$ for which $\gamma|_{Y_1}=\varphi$.
\ethm

The proof requires some preparations.

\bnota\label{cgw.n1} Further in this section we write $\A^n$ instead of $\A_\bk^n$ and the symbol
$\A_{u_1, \ldots , u_n}^n$ means that $\A^n$ is equipped with a fixed coordinate system
$\bar u= (u_1, \ldots, u_n)$. In particular, this system induces an embedding $\A^n \hookrightarrow \PP^n$
into a projective space. 
We  also suppose  that
$Z$, $Y_1$ and $Y_2$ are closed subvarieties of $\A^n_{u_1, \ldots , u_n}$ such that $\dim Z \leq n-3$, $\ED (Y_1)\leq n-2$,
$Y_1$ and $Y_2$ are disjoint from $Z$ and there exists an isomorphism $\varphi : Y_1 \to Y_2$.

\enota

The following result will be important this section: in the case when $Z$ is of codimension at least 2  Gromov \cite{Gr1}  observed that  $\A^n\setminus Z$ is a flexible variety 
and Winkelmann \cite{Wi}
showed that $\SAut ( \A^n\setminus Z)$ acts transitively on $ \A^n\setminus Z$ which is equivalent by virtue of Theorem \ref{fm.t1}.
(In particular, for a finite $Y_1$ Theorem \ref{cgw.t1} is valid even when $\dim Z=n-2$.)

The next fact is well-known (e.g., see \cite{Ka91}).

\bprop\label{cgw.p1}  Let  $\ED (Z) \leq k\leq n$ and $\Lin (\A^n, \A^k)$ be the affine variety of surjective linear maps  $\A^n \to \A^k$.
Then  for a general element $\rho \in \Lin (\A^n,\A^k)$ the restriction $\rho|_{Z} : Z \to \A^k$ is a closed embedding. 
\eprop

\bprop\label{cgw.p2} Let $\rho : \A^n \to \A^k$ be a general element of $\Lin (\A^n, \A^k)$.

{\rm (1)}  If  $k \leq \dim Z$ then $\rho|_Z : Z \to \A^k$ is surjective
and for every $w \in \A^k$ the fiber $F=\rho^{-1}(w)\cap Z$ is of dimension $\dim Z-k$. 

{\rm (2)} If $k=\dim Z$ then $\rho|_Z : Z \to \A^k$ is finite.

{\rm (3)} Let $T(Z/\A^k)$ be as in Notation \ref{gp1.d1}. Then for every $k$ one has $\dim T(Z/\A^k)\leq \dim TZ -k$.
\eprop

\bproof Let $\rho = (f_1, \ldots , f_k)$ be the coordinate form  of $\rho$ and $\bar \rho = (\bar f_1, \ldots , \bar f_k)$ be the rational extension of $\rho$ to $\PP^n$.
 Let $\cR (\bar \rho) = \bigcap_{i=1}^k \cR (\bar f_i)$ be the same as in Notation \ref{gp2.n1}. Since $\rho$ is surjective
(i.e., $f_1, \ldots , f_k$ are linearly independent) $\cR (\bar \rho)$ is of codimension $k$ in  $D=\PP^n\setminus \A^n$.
 The natural action of $\SL (n,\bk)$ on $\A^n$ extends to an action to $\PP^n$ whose restriction to $D$
is transitive. In particular, replacing $\rho$ by $\rho \circ h$ where $h$ is a general element of $\SL (n,\bk)$, by virtue of Theorem \ref{agga.t1} we can assure that 
 the intersection $\bZ \cap \cR (\bar \rho)$ is of dimension $\dim Z -1 -k$ where $\bZ$ is the closure $\bZ$ of  $Z$ in $\PP^n$
 (and, therefore,  $\bZ \cap \cR (\bar \rho)=\emptyset$ when $\dim Z=k$). 

Assume that there exists a fiber $F$ of $\rho|_Z$ with $\dim F > \dim Z -k$. Then the closure $\bF$ of $F$ meets $D$ along a subvariety of dimension at least $\dim Z -k$.
However, $\bF \cap D \subset \bZ \cap \cR (\bar \rho)$ which yields a contradiction. 
Thus $\dim F =\dim Z -k$ (since it cannot be less than $\dim Z -k$ \cite[Chap. 1, Sec. 6, Theorem 7]{Sh}).
  In particular, $\rho|_{Z} : Z \to \A^k$ is quasi-finite when $\dim Z=k$. Since $\bZ \cap \cR (\bar \rho) =\emptyset$ in the latter case
 the morphism $\rho|_Z : Z \to \A^k$ is proper by Proposition \ref{gp2.p1} and, therefore, it is finite by \cite[Theorem 8.11.1]{EGA}. 
 Hence, we have (1) and, in particular,  this map is surjective.  This implies that when $k < \dim Z$
we have also surjectivity for a general $\rho$ which is (2).
 
In (3) treat every point $x \in \A^n$ as a vector and consider the map $\psi: TZ \to \A^n$ that sends each $v \in T_zZ, \, z \in Z$
to $z+v$.  Denote by $V$ (resp. $V_0$) the closure of $\psi (TZ)$ (resp. $T(Z/\A^k)$) in $\PP^n$. Note that $\dim V$ does not exceed $\dim TZ$.
Hence, for a general $\rho$ and $\cR (\bar \rho)$ as above the dimension of $V\cap \cR (\bar \rho)$ does not exceed $\dim TZ -k-1$.
It remains to note that $V_0\cap D$ is contained in $V\cap \cR (\bar \rho)$ (indeed, if $v \in \Ker \rho_*$ then $f_1(v)=\ldots =f_k(v)=0$
and  the closure
of every line $\{ z +tv| \, t \in \bk \}$ in $\PP^n$ meets $D$ at a point of $\cR (\bar \rho)$). That is,  $\dim V_0\cap \cR (\bar \rho)\leq \dim TZ -k-1$
and, thus,  $\dim V_0\leq  \dim TZ -k$ which concludes the proof.
\eproof




\bnota\label{cgw.n2} For every group $G$ acting on an algebraic variety $X$ and every subscheme $W$ of $X$
we denote by $G_W$ the subgroup of $G$ such that the action of every element of $G_W$ yields the identity map on $W$.

\enota

\blem\label{cgw.l4} Let $X=P\times \A^{N+1}$,  $\rho : X \to P$ be the natural projection and $m>0$.
Suppose that $G=\SAut (X/P)$ and $X_p: =\rho^{-1}(p), \,p \in P$. Let $\Lin (\A^{N+1}, \A^{N})$ be the
variety of surjective linear maps $\A^{N+1}\to \A^{N}$ and $Z$ be a closed subvariety of $X$ such that for every $p \in P$
one has $\dim Z\cap X_p \leq N-1$
and for a general
$\lambda \in \Lin (\A^{N+1}, \A^{N})$ and $\theta = (\rho, \lambda)$ the morphism $\theta|_Z: Z \to P \times \A^{N}$ is  proper. 
Then each variety $X_p \setminus Z$ is $G_{Z_m}$-flexible where $Z_m$ is the $m$-th infinitesimal neighborhood
of $Z$ in $X$ and $G_{Z_m}$ has the same meaning as in Notation \ref{cgw.n2}.

\elem

\bproof     Let $\psi : P \times \A^{N} \to P$ be the natural projection. Since $\theta|_Z$ is proper one has
$\theta (Z\cap X_p)= \theta (Z) \cap \psi^{-1}(p)$. 
Let $x$ be any point in $X_p\setminus Z$. Since $\dim Z\cap X_p\leq N-1$ we see that $\theta (x) \notin \theta (Z\cap X_p)$
for a general $\lambda$ and hence $\theta (x) \notin \theta (Z)$.
Note also that $\theta$ can be treated as a quotient morphism of a locally nilpotent vector $\delta$ field on $\tX$ whose value at $x$ 
is a general vector $\nu$ in $T_x X_p$ (since $\lambda$ is general).
Note also that the kernel of $\delta$ coincides with  $\bk (P \times \A^{N})$ where we treat $\bk (P \times \A^{N})$ as a subring of  $\bk (X)$ under the natural embedding.
Hence for every regular function $f$ on $P\times \A^{N}$ that vanishes on $\theta (Z)$ but not at $\theta (x)$
(which exists since $\theta (Z)$ is closed in $P\times \A^{N}$) the vector field $f^m \delta$ is locally nilpotent
and its value at $x$ coincides with this general vector $\nu$ up to a nonzero factor. Furthermore, the flow of $f^m\delta$ is a one-parameter group in $G_{Z_m}$.
Hence $X_p\setminus Z$ is $G_{Z_m}$ flexible by the definition of flexibility and we are done.
 \eproof

\bnota{\label{cgw.n0}
Further in this section by $X$ we denote  $X=\A^n_{u_1, \ldots , u_n}\setminus Z$.}
\enota

\blem\label{cgw.l2} Let $\kappa : \A^n \to \A^s$ be the natural projection
$(u_1, \ldots , u_n) \to (u_1, \ldots, u_s)$ where $s \leq \ED (Y_1)$. Let $\kappa =\kappa \circ \varphi$, $\dim Z \leq n-s-2$
and $F\subset Y_1$ be a finite set for which $\varphi (y)=y$ for every $y\in F$.
Then  for every $m>0$ there exists an automorphism $\beta$ of $X\subset \A^n$ over $\A^s$ such that  for each $y \in F$
the $m$-jet of $\varphi$ at $y$ coincides with the restriction of the $m$-jet of $\beta$.
\elem

\bproof 
By Lemma \ref{cgw.l4} every fiber  of $\kappa|_X$ is $G_Z$-flexible.
By Theorem \ref{rv.t2} we can suppose now that the tangent spaces of  $Y_1$ and $Y_2$ 
at every $y \in F$ coincide and the isomorphism of these tangent spaces induced by $\varphi$ is the identity map. 
Without loss of generality we can also suppose that $T_yY_1$ is contained in $\{ u_{k+1}= \ldots = u_n=0\}$ where $k\geq s$
and $y$ is the origin $\bar 0$ of $\A^n$.
Then the $m$-jet $\varphi_m$ of $\varphi$ at $y$ can be chosen in the 
following form $(u_1, \ldots , u_n) \to (u_1+h_1, \ldots , u_n+h_n)$
where $h_1, \ldots, h_n$ are polynomials in $u_1, \ldots , u_k$ without free and linear terms. 
Let $\tilde h_{n} =\frac{\p h_1}{\p u_1}+\ldots +\frac{\p h_k}{\p u_k}$.
Then the $m$-jet $(u_1, \ldots , u_{n-1}, u_n) \to (u_1+h_1, \ldots , u_{n-1}+h_{n-1}, u_n-u_n\tilde h_{n})$ is of divergence 1 and, hence,
it satisfies the assumption on
the Jacobian in Notation \ref{rv.n2} (e.g., see the proof of \cite[Lemma 4.13]{AFKKZ}). 
By Theorem \ref{rv.t2} we can find an automorphism $\beta_1\in G_Z$ of $X$ over $\A^s$ with such precise jets at every point of $F$. 
Hence $\varphi_m$  is a composition of $\beta_1$ with an $m$-jet $\psi_m$ of the
form $(u_1, \ldots , u_{n-1}, u_n) \to (u_1, \ldots , u_{n-1}, u_n+ h)$ where a priori $h$ depends on $u_1, \ldots , u_n$.
However, restricting $\psi_m$ to the \'etale germ of $Y_1$ at  $y$ we can replace $h$ with a function of $u_1, \ldots , u_k$, i.e., 
we can view $\psi_m$ as the restriction of the $m$-jet of an automorphism $\beta_2\in G_Z$. Letting $\beta =\beta_2 \circ \beta_1$ we get 
the desired conclusion.
\eproof

\blem\label{cgw.l1} Let $\rho : \A^n \to \A^k$ be the natural projection $(u_1, \ldots , u_n) \to (u_1, \ldots, u_k)$
with $ k\geq \max (\dim Z  +1, \ED (Y_1))$.  Assume that for the  closure $Z'$ of $\rho (Z)$ in $\A^k$
the set $Y_1\cap \rho^{-1}(Z')$ is finite.    
Let  $\rho|_{Y_1} : Y_1 \to \A^k$ be a closed embedding and
$\varphi (y)=y$ for every  $y \in Y_1\cap \rho^{-1}(Z')$.  Then there exists $m>0$ such that whenever
the $m$-th jets of $\varphi$ are equal to the jets of the identity map for all $y\in Y_1\cap \rho^{-1}(Z')$, one has
an automorphism $\alpha$ of $\A^n$ over $\A^k$ for which the following holds

{\rm (a)} $\alpha (z)=z$ for every $z \in Z$ (in particular, $\alpha$ is an automorphism of $X =\A^n \setminus Z$)  and

{\rm (b)}  after replacing $Y_1$ by $\alpha (Y_1)$ and $\varphi$ by $\varphi \circ \alpha^{-1}$ one has $u_i |_{Y_1}=u_i\circ \varphi$
for every $u_i $ with $i \geq k+1$. 

\elem

\bproof Let $I'\subset \bk [\A^k]$ be the vanishing ideal of $Z'$,  $Y_1'=\rho (Y_1)$, $I\subset \bk [Y_1']$ be the image of $I$ under the 
natural morphism $\bk [\A^k]\to \bk [Y_1']$ and $J\subset \bk [Y_1']$ be the vanishing ideal of $Y_1'\cap Z'$. By Nullstellensatz one has
$J^m \subset I$ for some $m>0$. 
Let $g_i=u_i \circ \varphi -u_i|_{Y_1}$ for $i \geq k+1$. Since and $Y_1$ and $Y_1'$ are isomorphic we can consider $g_i$ as a function on $Y_1'$. 
Note that it always belongs to $J$ and, furthermore, under the assumption it belongs to $J^m$.  
Thus $g_i$ is the restriction of a polynomial $\tilde g_i$ on $\A^k$ which vanishes on $Z'$.
Consider the automorphism $\psi_{i}$ of $X$ given by $u_i \to u_i+\tilde g_i(u_1, \ldots , u_k)$
and $u_j \to u_j$ for $j\ne i$.
Then the composition $\psi_{k+1} \circ \ldots \circ \psi_{n}$  yields the desired automorphism $\alpha$.
\eproof

\brem\label{cgw.r1} 
Since $m$ in the proof above can be chosen as large as we wish one can suppose that the extension $\tilde g_i$ vanishes not only on
$Z'$ but also on the $m$-th infinitesimal neigborhood of $Z'$. 
This implies that the  automorphisms $\{ \psi_{i} \}$ in the proof are elements of the flows of locally nilpotent vector fields that vanish on $Z$ with multiplicity $m$.
Since $\alpha$ from Lemma \ref{cgw.l1} is a composition of such automorphisms we can suppose that
the restriction of $\alpha$ to the $m$-th infinitesimal neigborhood $Z_m$ of $Z$ is the identity map. The same is true for $\beta$ from
Lemma \ref{cgw.l2} since  by Lemma \ref{cgw.l4} $X_p\setminus Z$ is not only $G_{Z}$-flexible
but also $G_{Z_m}$-flexible.

\erem

\bproof[Proof of Theorem \ref{cgw.t1}  (a)]
Let $H=\SL(n,\bk)$ act naturally on $\A^n$,
$\rho_{k}  : \A^n \to \A^k$ be  the natural projection given by $(u_1, \ldots , u_n) \to (u_1, \ldots , u_k)$ and
$\rho_{k,h}$ be the composition $ \rho_k \circ h$ where $h\in H$.  
By Propositions \ref{cgw.p1} and \ref{cgw.p2} we have the following for a general $h\in H$.

(1)   For $k= \dim Z$ the morphism $\rho_{k,h}|_{Z}: Z \to \A^{k}$ is finite
while the morphism $\rho_{n-2,h}|_{Z}: Z \to \A^{n-2}$ is a closed embedding.

(2)  
For every $k<\dim Y_2$ the morphism $\rho_{\sigma,k}|_{Y_2}: Y_2 \to \A^{k}$ is surjective
with the dimension of each fiber equal to $\dim Y_2 -k$, for $k=\dim Y_2$
the morphism $\rho_{k,h}|_{Y_2}: Y_2 \to \A^{m}$ is finite and
for every $k\leq \dim Y_2$ one has $\dim T(Y_2/\A^k) \leq \dim TY_2-k$.

Let us show inductively that for every $k \leq \dim Y_1$ and a general $h$ in $H$ one can also assume the following.

(3)  $\rho_{k,h} \circ \varphi =\rho_{k,h}|_{Y_1}$ and, in particular, by (2) 
if $k<\dim Y_1$ then $\rho_{\sigma,k}|_{Y_1}: Y_1\to \A^{k}$ is surjective
with the dimension of each fiber equal to $\dim Y_1 -k$, 
if $k=\dim Y_1$ then $\rho_{k,h}|_{Y_1}: Y_1 \to \A^{m}$ is finite and in any case
$\dim T(Y_1/\A^k) \leq \dim TY_1-k$.

(4)  $\rho_{n-2,h}|_{Y_1}: Y_1 \to \A^{n-2}$ is a closed embedding.

(5) $\rho_{n-2,h} (Y_1)$ does not meet $\rho_{n-2,h} (Z)$.

For $k=0$ the assumption (3) is automatically true  while (4) and (5) follow from Proposition \ref{cgw.p1} applied to the variety  $Z\cup Y_1$. 
Let us assume that (3)-(5) are valid
for $k-1$ and let us establish that this implies their validity for $k$.
Since (4) and (5) hold for $k-1$, by Lemma \ref{cgw.l1} we obtained an automorphism $\beta \in \SAut (X)$ over $\A^{n-2}$ such that after replacing 
$Y_1$ by $\beta (Y_1)$ one has $u_i |_{Y_1}=u_i\circ \varphi$  for $i=n-1$ and $i= n$.  
Consider the element $\sigma$ of $H$ such that $\sigma_* (u_k)=u_n, \sigma_* (u_n)=-u_k$ while the rest coordinates remain intact.
Note that since $h$ is general we can suppose that $\sigma \circ \rho_{k,h}$ is still of the form $\rho_{k,h_1}$ where $h_1$ is
a general element of $H$. Hence, replacing $\rho_{k,h}$ with $\sigma \circ \rho_{k,h}$ we get (3) while keeping (1) and (2).

Treat now $\rho_{k,h} : X \to \A^k$ (resp. $\rho_{n-2,h}: X \to \A^{n-2}$) as $\kappa : X \to P$ (resp. $\rho : X \to Q)$) in Theorem \ref{gp1.t1}.
Note the assumptions of Theorem \ref{gp1.t1} are valid for the group $G=\SAut(X/P)$ and the subvariety $Y_1$ of $X$ (indeed, the
fibers of $\kappa$ is $G$-flexible by Lemma \ref{cgw.l4} while $\dim Y_1 \times_PY_1\leq 2 \dim Y_1 -\dim P$ by Remark \ref{agga.r1}
and the fact that the fibers of $\kappa|_{y_1} : Y_1 \to P$ are equidimensional because of (3)).
By Theorem \ref{gp1.t1} (i) and (iii) we can find an algebraic family $\cA_1\subset G$  (a priori depending on $h$) such that for a general element $\alpha_1 \in \cA_1$
the morphism $\rho|_{\alpha_1 (Y_1)}: \alpha_1 (Y_1) \to Q$ is injective and generates an injective map of the Zariski tangent bundles.
Similarly, since $\rho^{-1} (\rho (Z))$ is of dimension $\dim Z+2 < \dim X -\dim Y_1$ by Theorem \ref{agga.t2} we can find
an algebraic family $\cA_2\subset G$  (also depending on $h$) such that for a general element $\alpha_2 \in \cA_2$ the variety
$\alpha_2 (Y_1)$ does not meet $\rho^{-1} (\rho (Z))$ which is equivalent to (5). By Remark \ref{gp1.r1} (2) we can suppose that
$\cA_1=\cA_2$.

Hence, in order to prove (4) and (5) one needs to show that  the morphism $\rho|_{\alpha (Y_1)}: \alpha (Y_1) \to Q$ can be made proper.

Let $H'$ be the subgroup of $H$ consisting of all elements whose action on the first summand of 
$\A^n=\A_{u_1, \ldots , u_k}^k\oplus \A_{u_{k+1}, \ldots , u_{n-2}}^{n-k-2}\oplus \A^2_{u_{n-1},u_n}$ is the identity map.
Consider Proposition \ref{gp1.p3} with $U=H'$ and $\check \rho : H' \times X =\cX \to  \cQ =H \times Q$ given by $ (h',x) \to (h',\rho_{n-2,h'h}(x))$. 
Then Proposition \ref{gp1.p3} implies that we can find a family $\cA_1$ as above which provides the injectivity condition for general $h'$ in $H'$
simultaneously. Note that the composition of the $H'$-action
with the projection $\A^n \to \A_{u_{k+1}, \ldots , u_{n-2}}$ can be viewed as the variety $\Lin(\A^n, \A^{n-k-2})$ of surjective
linear maps $\rho_{n-k-2,h'h}' : \A^n \to \A^{n-k-2}$ and, hence,  Proposition \ref{cgw.p2} implies
that $\rho'_{n-k-2,h'h}|_{Y_1} : Y_1 \to \A^{n-k-2}$ is proper for a general $h'$ because $\dim Y_1 < n-k-2$
by the assumption of Theorem \ref{cgw.t1}. Since for a general $h'$ the element $h'h$
is still general in $H$, replacing $h$ by $h'h$ we get (4) without ruining (1)-(3) while the second statement of Proposition \ref{gp1.p3}
provides us with (5).

Now we shall establish by induction conditions (3)-(5) for a general $h$ and  $k=k_1+1, \ldots , n$ where $k_1=\dim Y_1$.
Since $\rho_{h,k_1}|_{Y_1}: Y_1 \to \A^{k_1}$ is finite by (3)  we can suppose now every $\rho_{h,k}|_{Y_1}: Y_1 \to \A^k$ is proper.
Assume that (4) and (5) hold for $k-1$.  Applying Lemma \ref{cgw.l1} as above we can achieve (3) for $k$ without ruining this properness
(but may be violating (4) and (5)). 
Using notation $\kappa : X \to P$ and $\rho : X \to Q$ in the same meaning as before we observe that though 
the assumption $\dim Y_1 \times_PY_1\leq 2 \dim Y_1 -\dim P$ does not hold we can still apply Theorem \ref{gp1.t1} (i) and (iv)
by virtue of Remark \ref{gp1.r1} (5).
Hence the same argument as above yields an algebraic family $\cA \subset G$ of automorphisms of $X$ over $\A^{k}$
such that for a general  $\alpha \in \cA$
the morphism $\rho_{n-2,h}|_{\alpha (Y_1)}: \alpha (Y_1) \to Q$ is injective and generates an injective map of the Zariski tangent bundles,
and, furthermore,  $\alpha (Y_1)$ does not meet $\rho_{n-2,h}^{-1} (\rho_{n-2,h} (Z))$. The latter fact is (5) and in combination with the properness
the former one yields (4). 

Note that this process of replacing $Y_1$ by its automorphic image in $X$ in the case of $k=n$ yields the equality $Y_1=Y_2$ and, hence,
the desired conclusion.\eproof

\blem\label{cgw.l3} Let $m>0$ and $\cL$ be a finite collection of $3^m$ distinct lines in $\PP^2$
such that the intersection of any three of them is empty. Then for every finite subset $F\subset \PP^2$ of cardinality $m$ there is a line $L \in \cL$ such that $L \cap F=\emptyset$.

\elem

\bproof Let us use induction by $m$. For $m=1$ the claim follows from the definition of $\cL$. Suppose that 
the statement is true for $m-1$ and $\cL$ is a collection of $3^m$ lines with the desired property. 
Present $\cL$ as a disjoint union $\cL_1 \cup \cL_2\cup \cL_3$ where each collection $\cL_i$ consists of $3^{m-1}$ lines.
Let  $F=F_0\cup \{x \}\subset \PP^2$  be of cardinality $m$, i.e., $F_0$ is of cardinality $m-1$.
By the assumption we can find lines $L_1 \in \cL_1$, $L_2 \in \cL_2$ and $L_3 \in \cL_3$ such that $L_i \cap F_0=\emptyset$ for every $i$.
Since $L_1\cap L_2\cap L_3=\emptyset$ one of these lines does not contain $x$ and we are done.
\eproof

\bproof[Proof of Theorem \ref{cgw.t1} (b)] 
Let $\bar u=(u_1, \ldots, u_n)$ and $\rho_{k,h} : \A^n \to \A^k$ have the same meaning as in the proof of Theorem \ref{cgw.t1} (a) where $h \in H=\SL(n,\bk)$.
We can suppose again that (1) as in the proof of Theorem \ref{cgw.t1} (a) holds.

Let us show inductively that for $k\leq n-3$ and some general $h$ (depending on $k$) 
we can also assume 

(2) $\rho_{k,h} \circ \varphi =\rho_{k,h}|_{Y_1}$ 

(3)  $\rho_{n-2,h}|_{Y_1}: Y_1 \to \A^{n-2}$ is a closed embedding.

(4) $\rho_{n-2,h} (Y_1)$ meets $\rho_{n-2,h} (Z)$ at a finite number of points and, furthermore, these points are the images of smooth points of $Y_1$. 

(5) for every singular point $y$ of $Y_1$ the image of $T_yY_1$ in $T\A^k$ under the map induced by $\rho_{k,h}$ is of
dimension $\min (k, \dim T_yY_1)$ (which implies that $\dim T(Y_1/\A^k) \leq \dim TY_1-k$) and the cardinalities
of ${\rm Sing} \, Y_1$ and $\rho_{n-2,h} ({\rm Sing} \, Y_1)$ are the same.

{\em Step 1} (condition (5) and the case of $k=0$). Let  $k=0$ and $h$ be a general element of $H$. Then the assumption (2) is automatic and (3) follows from Proposition \ref{cgw.p1}.
Denote by $E$ the set of singular points of $Y_1$. Since $X$ is $\SAut (X)$-flexible,  replacing $Y_2$ by its automorphic image in $X$
we can suppose by Theorem \ref{fm.t1a}  that $E$ is also the set of singular points of $Y_2$ and $\varphi (y)=y$ for every $y \in E$.
Furthermore, by Theorem \ref{fm.t2} we can suppose that for every $y\in E$ and every $m\geq 0$ the image of $T_yY_2$ in $T\A^m$ under the map induced by 
$\rho_{m,h}$ is of dimension $\min (m, \dim T_yY_2)$. Applying now an automorphism of $X$ to $Y_1$ we can suppose
also that $\varphi_* : T_yY_1\to T_yY_2$ is the identity map. Let $G$ be the subgroup of $\SAut (X)$ whose restriction to
the second infinitesimal neighborhood $E_2$ of $E$ is the identity map. From now on we are going to use replacements of 
$Y_1$ and $Y_2$ by their images under the actions of some elements of $G$. Note that this quarantees (5) for all $m\geq 0$.
Returning to the case of $k=0$ we observe that since $h$ is general by Proposition \ref{cgw.p1} the morphism $\rho_{n-2, h}|_{E \cup Z}$ is injective (and, in particular,
the sets $\rho_{n-2, h}(E)$ and $\rho_{n-2,h} (Z)$ are disjoint) and  also $\dim Y_1 + \dim \rho_{n-2,h}^{-1}(\rho_{n-2,h}(Z))=\dim X$. 
Hence, since $X'=X\setminus E$ is $G$-flexible, applying to $Y_1$ an element of $G$ we get (4)  by Theorem \ref{agga.t2} which concludes Step 1.

{\em Step 2.} Now let us assume that (1)-(4) are valid for $k-1$ and proceed with
the case of $k\leq n-3$. 

Let $F_1=Y_1 \cap \rho_{n-2,h}^{-1}(\rho_{n-2,h}(Z))$ and $F_2=\varphi (F_1)$.
Since $\rho_{k-1,h}\circ \varphi = \rho_{k-1, h}|_{Y_1}$ we see that the sets $\rho_{k-1,h}^{-1} (w) \cap F_1$ and $\rho_{k-1,h}^{-1} (w) \cap F_2$
are of the same cardinality for every point $w \in \A^{k-1}$. Since by Lemma \ref{cgw.l4} every fiber of $\rho_{k-1,h}|_{X'}$ is
$G$-flexible by Theorems \ref{fm.t2} and \ref{rv.t1} we can find an automorphism $\alpha_1 \in G$
such that replacing $Y_2$ by $\alpha_1 (Y_2)$ one has $F_1=F_2$ and $\varphi|_{F_1}={\rm id}|_{F_1}$.
Suppose that $m$ is as in Lemma \ref{cgw.l1} (with $\rho$ replaced by $\rho_{n-2,h}$). Choose $\beta$
as in Lemma \ref{cgw.l2} and replace $Y_2$ by $\beta (Y_2)$, i.e., the $m$-jet of $\varphi$ at any point of $F_1$ is now
the $m$-jet of the identity map\footnote{By flexibility such $\beta$ can be chosen as an element of $G$
since the finite sets $E$ and $F_1$ are disjoint.}. By Lemma \ref{cgw.l1} we can find $\alpha_2 \in G$ such that after replacement of $Y_1$ by
$\alpha_2(Y_1)$ one has $u_n|_{Y_1}=u_n\circ \varphi$.
Consider $\sigma \in H$ such that $\sigma_* (u_k)=u_n, \sigma_* (u_n)=-u_k$ while the rest coordinates remain intact.
Since we started with a general $h$ one has $\sigma \circ \rho_{k,h}=\rho_{k,h_1}$ where $h_1$ is
a general element of $H$. Thus, replacement of $h$ by $h_1$ leave (1) valid while providing us with (2).

If we want to apply Theorem \ref{gp1.t1} further  we need to guarantee the condition (\ref{gp1.eq1}) in Remark \ref{gp1.r1} (5)
with $\rho_{k,h}$ as $\kappa$. For $k=1$ we can get this condition by
choosing a couple of points $y_1$ and $y_2$ in each irreducible component of $Y_2$ and requiring that for $\alpha_1$
as above one has additionally $\rho_{1,h} (\alpha_1 (y_1))\ne \rho_{1,h}(\alpha_1 (y_2))$. Then $\rho_{1,h}|_{Y_2}: Y_2 \to \A^1$ is a quasi-finite morphism
(after the replacement of $Y_2$ by $\alpha_1(Y_2)$)
and we have (\ref{gp1.eq1}) for $Y_2$. By (2) we have it also for $Y_1$.
Furthermore,  this condition will survive for $k>1$ provided that in further replacements of $h$ by another general $\tilde h\in H$ one has $\rho_{2,h}=\rho_{2,\tilde h}$.

As in the previous proof we treat  $\rho_{k,h} : X' \to \A^k$ (resp. $\rho_{n-2,h}: X' \to \A^{n-2}$) as $\kappa : X \to P$ (resp. $\rho : X \to Q)$) in Theorem \ref{gp1.t1}.
By Theorem \ref{gp1.t1} (i) and (iii) we can find an algebraic family $\cA_3\subset G$ such that for a general element $\alpha_3 \in \cA_3$
the morphism $\rho|_{\alpha_3 (Y_1\setminus E)}: \alpha_3 (Y_1\setminus E) \to Q$ is injective and generates an injective map of the Zariski tangent bundles.
By (5) this is also true for the morphism $\rho|_{\alpha_3 (Y_1)}: \alpha_3 (Y_1) \to Q$.
Similarly, since $\rho^{-1} (\rho (Z))$ is of dimension $\dim Z+2 \leq \dim X -\dim Y_1$ by Theorem \ref{agga.t2} we can find
an algebraic family $\cA_4\subset G$  such that for a general element $\alpha_4 \in \cA_4$ the variety
$\alpha_4 (Y_1)$ meet $\rho^{-1} (\rho (Z))$ at a finite number of points which is equivalent to (4). By Remark \ref{gp1.r1} (2) we can suppose that
$\cA_3=\cA_4$. We need to show that  the morphism $\rho|_{\alpha (Y_1)}: \alpha (Y_1) \to Q$ can be made proper for which
we need to change our general element $h \in H$.

Let $H'$ be the subgroup of $H$ as in the proof of Theorem \ref{cgw.t1} (a). Applying the same argument as in that proof (i.e., using Proposition \ref{gp1.p3} (1)) 
we see that for a general $h' \in H'$ the replacement of $h$ by $\tilde h = h'h$ provides us
with properness of the morphism
$\rho'_{n-k-2,h'h}|_{Y_1} : Y_1 \to \A^{n-k-2}$  as soon as $n-k-2 \geq \dim Y_1= 1$ (i.e., $k\leq n-3$) while preserving
the injectivity  of these morphism and injectivity of the induced morphism of the Zariski tangent bundles
(as well as (1), (2), and (5)).  This yields (3) and, similarly, by Proposition \ref{gp1.p3} (3) we get (4).
Thus we have $\rho_{n-3,h} \circ \varphi =\rho_{n-3,h}|_{Y_1}$ which concludes Step 1. 

{\em Step 3}.  Let us show that replacing $Y_1$ and $Y_2$ with their automorphic images in $X$ 
and replacing $h$ by another general element of $H$ we can suppose that 
$\rho_{n-2,h} \circ \varphi =\rho_{n-2,h}|_{Y_1}$ while keeping the morphism $\rho_{n-2,h}|_{Y_1}: Y_1 \to \A^{n-2}$ finite.

Consider the embedding $\A^n=P\times \A^3_{u_{n-2}, u_{n-1},u_n} \subset P\times \PP^3=:W$
where $P=\A^{n-3}_{u_1, \ldots , u_{n-3}}$.
Let $\tD =W\setminus \A^{n}=P\times \PP^2$ and for every subset $Y$ of $\A^n$ denote by $\breve Y$ the intersection of 
$\tD$ with the closure of $Y$ in $W$. Note that if $Y$ is a curve then $\breve Y$ consists of no more than $l$ points where $l$
is the number of punctures of $Y$ (that are the points in the complement to a normalization of $Y$ in a smooth completion of this normalization).
Furthermore, let $V\subset \A^n$ be given by $u_n=0$. Then $\breve V=P \times L$ where $L$ is a line in $\PP^2$
and the natural morphism $Y\to \A^1_{u_n}$ is finite if and only if $\breve Y$ does not  meet $P\times L$. Consider the natural action
of $H''=\SL(3, \bk)$ on $\A^3_{u_{n-2}, u_{n-1},u_n}$ and  extend this action to the action on $\A^n=\A^{n-3}\oplus \A^3$ such that
it is trivial on the first summand (this enables us to treat this $H''$ as a subgroup of $H=SL(n,\bk)$ acting naturally on $\A^n$).
Consider general elements $h_1'', \ldots , h_{3^l}''$ of $H''$ where $l$ is the same as above but for $Y=Y_1$.
Let $V_i$ be the zero locus of $u_n\circ \rho_{h_i''h}$. Then $\breve V_i=P\times L_i$ where $L_1, \ldots , L_{3^l}$ are general lines 
in $\PP^2$ and, hence, the intersection of any three of them is empty. 

Consider $F_1^i=Y_1 \cap \rho_{n-2,h_i''h}^{-1}(\rho_{n-2,h_i''h}(Z))$ and $F_2^i=\varphi (F_1^i)$. As before
replacing $Y_2$ by its automorphic image one has $F_1^i=F_2^i$ and $\varphi|_{F_1^i}={\rm id}|_{F_1^i}$ for every $i$.
Furthermore, for every point in $\bigcup F_1^i$ we can suppose that the $m$-jet of $\varphi$ coincides with the $m$-jet of
the identity map where $m$ is as in Lemma \ref{cgw.l1}.  By Lemma \ref{cgw.l3} for one of the lines, say $L_1$, the intersection
$\breve Y_2$ with $P\times L_1$ is empty and thus for $u_n''=u_n\circ (h_1''h)$ the morphism $Y_2 \to A^1_{u_n''}$ is finite.
Applying Lemma \ref{cgw.l1} we can replace $Y_1$ by its automorphic image such that $u_n''|_{Y_1}=u_n''\circ \varphi$.
Exchanging the role of $u_{n-2}$ and $u_n$ we can suppose that $\rho_{n-2,h_1''h} \circ \varphi =\rho_{n-2,h_1''h}|_{Y_1}$.
Since $h_1''h$ is a general element of $H$, after the replacement of $h$ by $h_1''h$ we still have (1) and (5). The finiteness of the morphism
$Y_2 \to A^1_{u_{n-2}''}$ implies
also that the morphism $\rho_{n-2, h_1''h}|_{Y_1}: Y_1 \to \A^{n-2}$ is finite which concludes Step 3. 

{\em Step 4.} To finish the proof let us treat
 $\rho_{n-2,h} : X' \to \A^{n-2}$ (resp. $\rho_{n-1,h}: X' \to \A^{n-1}$) as $\kappa : X \to P$ (resp. $\rho : X \to Q)$) in Theorem \ref{gp1.t1}.
By Theorem \ref{gp1.t1} (i) and (iii) (in combination with (5)) we can find an algebraic family $\cA_5\subset G$ such that for a general element $\alpha_5 \in \cA_5$
the morphism $\rho|_{\alpha_5 (Y_1)}: \alpha_5 (Y_1) \to Q$ is injective and generates an injective map of the Zariski tangent bundles.
Being also finite (since $\kappa|_{Y_1}: Y_1\to P$ is finite and $\alpha_5$ is in $\Aut (X/P)$) 
it is a closed embedding. Furthermore, by Theorem \ref{agga.t2} we can suppose that 
$\alpha_5 (Y_1)$ does not meet $\rho^{-1}(\rho (Z))$ since the sum of the dimensions of these varieties
is less than $\dim X$. Consider the natural action
of $H'''=\SL(2, \bk)$ on $\A^2_{ u_{n-1},u_n}$ and  extend this action to the action on $\A^n=\A^{n-2}\oplus \A^2$ such that
it is trivial on the first summand (i.e., $H'''$ is again a subgroup of $H$ acting naturally on $\A^n$).
By Proposition \ref{gp1.p3} for a general  $h'''\in H'''$ the morphism
$\rho_{n-1, h'''h}|_{\alpha_5 (Y_1)}: \alpha_5 (Y_1) \to Q$ is still a closed embedding with $\alpha_5 (Y_1)$ not meeting
$\rho_{n-1, h'''h}^{-1}(\rho_{n-1, h'''h} (Z))$.  Since $h'''$ is general the same is true for the replacement of $h$ by $\sigma'''h'''h$ where
$\sigma'''\in H'''$ is up to a sign the transposition of coordinates $u_{n-1}$ and $u_n$.
By Lemma \ref{cgw.l1} (with $k=n-1$) replacing $Y_1$ by its automorphic image we can suppose that for $u_n'''=u_n \circ (h'''h)$
one has $u_n'''|_{Y_1}=u_n'''\circ \varphi$. Exchanging the role of $u_{n-1}$ and $u_n$ 
we can suppose that $\rho_{n-1,\sigma'''h'''h} \circ \varphi =\rho_{n-1,\sigma'''h'''h}|_{Y_1}$.
The assumptions of Lemma \ref{cgw.l1} are true with $\rho$ in that lemma equal to $\rho_{n-1,\sigma'''h'''h}$.
Thus, after additional replacement of $Y_1$ by its automorphic image we can make $Y_1=Y_2$ which is the desired conclusion.

\eproof

\brem\label{cgw.r2} (1) The assumption that $Y_1$ and $Y_2$ are closed in $\A^n$ (and not in $\A^n \setminus Z$) cannot be dropped.
Indeed, consider $Z$ that does not admit nontrivial automorphisms, and let $L_1\simeq \A$ and $L_2\simeq \A$ be disjoint curves in $\A^n$
each of which meets $Z$ at one point only. Then $Y_1=L_1 \setminus Z$ and  $Y_2=L_2 \setminus Z$ are isomorphic but there is no way to extend
this isomorphism to an automorphism of $\A^n \setminus Z$.

 (2)   We constructed an automorphism $\gamma$ of $X$ such that $\gamma|_{Y_1}=\varphi$ as a composition of elements of $G_Z$.
 However, Remark \ref{cgw.r1} implies that we can choose $\gamma$ as a composition of elements of $G_{Z_m}$.
 More precisely, the proof of Theorem \ref{cgw.t1} yields the following fact.
\erem

\bprop\label{cgw.p3} For every $m>0$ the automorphism $\gamma$ from Theorem \ref{cgw.t1}
can be chosen as a composition of elements of the flows of locally nilpotent vector fields that vanish on $Z$ with multiplicity $m$.
In particular, the restriction of $\gamma$ to the $m$-th infinitesimal neighborhood $Z_m$ of $Z$ is the identity map.
\eprop

It is interesting to understand how sharp is Theorem \ref{cgw.t1}. Therefore, we would like to pose the following.

{\bf Question.} {\em Let $Z$ be a closed subvariety in $\A^n$ of codimension 2, and let $Y_1$ and  $Y_2$ be two smooth closed
isomorphic curves in $\A^n$ disjoint from $Z$. Can one always find an automorphism of $\A^n \setminus Z$
transforming $Y_1$ onto $Y_2$?}

\section{The case of quadrics}

\bnota\label{cqu.n1} In this section $m\geq 6$ and
$X$ is a hypersurface in $\A_\bk^m$ that is a nonzero fiber of a non-degenerate quadratic form.
That is, when $m=2n$ (resp. $m=2n+1$) we can suppose that $X$ is given by $u_1v_1+ \ldots + u_nv_n=1$
(resp. $u_0^2+ u_1v_1+ \ldots + u_nv_n=1$) in a suitable coordinate system $(\bar u, \bar v)=(u_1, \ldots, u_n, v_1, \ldots ,v_n)$(resp. $(u_0, \bar u, \bar v)$).
Recall that such an $X$ is a homogeneous space of a special orthogonal group $G: ={\rm SO} (m, \bk )$ acting linearly on $\A_\bk^m$. 
In particular, it is flexible since homogeneous spaces of any semi-simple group are flexible \cite{AFKKZ}. Furthermore,
the coordinate system determines an embedding $\A_\bk^m \hookrightarrow \PP_\bk^m$ such that the action of $G$ 
extends to $\PP_\bk^m$. The closure $\bX$ of $X$ in $\PP_\bk^m$ yields a completion of $X$ for which $\bX \cap H$ is a quadric
in $H =\PP_\bk^m\setminus \A_\bk^m\simeq \PP_\bk^{m-1}$. 
\enota

\blem\label{cqu.l1}  Let Notation \ref{cqu.n1} hold. Then $G$ acts transitively on $\bX \cap H$.

\elem

\bproof Making a linear coordinate change we can suppose that the coordinate system $(w_1, \ldots , w_m)$ on $\A_\bk^m$
is such that $X$ is given by $w_1^2+ \ldots + w_m^2=1$. Then $\PP_\bk^m$ is equipped with the coordinate
system $(\tilde w_0: \tilde w_1: \ldots : \tilde w_m)$ so that $w_i=\tilde w_i/\tilde w_0$ for $i=1, \ldots ,m$
and $H$ is given by $\tilde w_0=0$. Hence the equation of $\bX \cap H$ in $H$ is
$\tilde w_1^2+ \ldots + \tilde w_m^2=0$. 

Let $Z_i=\bX \cap H \cap \{ \tilde w_i=0\}$ and $Z_i'=(\bX \cap H)\setminus Z_i$.
Note that $Z_i'$ is isomorphic to a nonzero fiber of  a non-degenerate quadratic form on $\A_\bk^{m-1}$ and, therefore,
the subgroup $G_i\simeq {\rm SO} (m-1, \bk )$ of $G$ that preserves the coordinate $w_i$ acts transitively on $Z_i'$. 
On the other hand
for a point $z \in Z_m$ at least one of the coordinates $\tilde w_1 :\ldots : \tilde w_{m-1}$ is nonzero (say $\tilde w_1$).
Furthermore, applying an element of $G_m$ we can make $\tilde w_2$ also different from zero.
Hence the subgroup $G_2$ can transform $z$ into
a point of $Z_m'$ since the action of $G_2$ on $Z_2'$ can switch coordinates $\tilde w_1$ and $\tilde w_m$
(and change the sign of one of them). This yields the desired statement about transitivity because
$\bX \cap H=Z_m \cup Z_m'$.
\eproof

The aim of this section is the following theorem.

\bthm\label{cqu.t1} Let Notation \ref{cqu.n1} hold and
$\varphi : Y_1\to Y_2$ be an isomorphism of two closed subvarieties of $X$. Let  $\ED (Y_i) +\dim Y_i \leq m-2$.
Then $\varphi$ extends to an automorphism of $X$ which belongs to $\SAut (X)$.

\ethm

\bproof Consider the case of $m=2n$, i.e., $X \subset \A_\bk^m$ is given by  $u_1v_1+ \ldots + u_nv_n=1$.
Note that for $m\geq 6$ the assumption that $\ED (Y_i) +\dim Y_i \leq m-2$ implies that 
$k:=\dim Y_i +2 \leq n$. 
Let $(u_1, \ldots, u_n, v_{k+1}, \ldots , v_n)$ be a coordinate system on $\A_\bk^{m-k}$ and 
$\rho : \A_\bk^m \to \A_\bk^{m-k}$ be the natural projection. Suppose that 
$\bar \rho : \bX \dashrightarrow \A_\bk^{m-k}$ is the rational
extension of $\rho|_X$. Then $\cR (\bar \rho)$ as in Notation \ref{gp2.n1} has dimension $k-1$.
Note also that if $\bY_j$ is the closure of $Y_j$ in $\bX$ then $\dim \bY_i \cap H <m-k-2$ since
$\ED (Y_i)\leq m-k$. Hence $\dim \bY_j\cap H + \dim \cR (\bar \rho)< \dim \bX \cap H$.
Lemma \ref{cqu.l1} implies that Corollary \ref{gp2.c2}  is applicable and
for a general element $\alpha$ of some algebraic family $\cA\subset \SAut (X)$
of automorphisms $\rho|_{\alpha (Y_j)} : \alpha (Y_j) \to \A_\bk^{m-k}$ is a closed embedding. Furthermore,
by Theorem \ref{agga.t1} we can suppose that $\alpha (Y_j)$ does not meet the subvariety $F$ of $X$
given by $u_2=\ldots = u_k=0$ since $\dim Y_j < k-1=\codim_X F$. Hence replacing each $Y_j$
by its automorphic image we suppose that $\rho|_{Y_j} : Y_j \to \A_\bk^{m-k}$ is a closed embedding
and $Y_j \cap F =\emptyset$. This implies that $Y_j':=\rho (Y_j)$ is isomorphic to $Y_j$ and it
does not meet the subspace of $\A_\bk^{m-k}$
given by the same system of equations $u_2=\ldots = u_k=0$. In particular, we can treat $v_1|_{Y_j}$
as the lift of a function $f_j$ on $Y_j'$ and by the Nullstellensatz one has $1-f_j=\sum_{i=2}^k u_ig_{j,i}$ for some regular
functions $g_{j,2}, \ldots , g_{j,k}$ on $Y_j'$.

{\em Claim}. For every $Y_j$ as in the statement of the Theorem there exists an automorphism of $X$
that sends $Y_j$ onto a subvariety of  the hypersurface $S$ in $X$ given by $v_1=1$ (in particular, $S$ is isomorphic to $\A_\bk^{m-2}$).

Indeed, let $\delta_i, \, i=2, \ldots , k$ be the locally nilpotent 
vector fields on $X$ given by $\delta_i = u_i \frac{\p}{\p v_1} - u_1\frac{\p}{\p v_i}$.
Let $\psi_{j,i}$ be the flow of $g_{j,i}\delta_i$ at time $t=1$. Note that the automorphism $\psi_{j,2} \circ \ldots \circ \psi_{j,k}$
transforms $Y_j$ into a subvariety of $X$ on which the restriction of $v_1$ is identically 1. This concludes the proof
of the Claim.

Thus we can suppose from the beginning that $Y_1$ and $Y_2$ are contained in $S$. 
By Theorem \ref{in.t0} there is an automorphism $\beta$ of $S$ that transforms $Y_1$ into $Y_2$. 
Since $X\setminus \{ v_1=0\} \simeq \A_\bk^*\times S$ we can suppose that $\beta$ is the restriction of an
automorphism of $X\setminus \{v_1=0\}$. Thus by Theorem \ref{rv.t1}  $\beta$ can be extended to an automorphism
of $X$ which concludes the case of $m=2n$.

If $m=2n+1$ then we treat $(u_0, u_1, \ldots, u_n, v_{k+1}, \ldots , v_n)$ as a coordinate system on $\A_\bk^{m-k}$
and consider the natural projection $\rho : \A_\bk^m \to \A_\bk^{m-k}$.
The rest of the proof works without change and we are done for the case of $m=2n+1$ as well.
\eproof

The assumption on $\ED (Y_i)$
 implies that in the smooth case (i.e., in the case when $\ED (Y_i)=2\dim Y_i+1$) we have the following.

\bcor\label{cqu.c1} Let $X$ be a nonzero fiber of a non-degenerate quadratic form in $\A_\bk^m$.
Suppose that $\varphi : Y_1\to Y_2$ is an isomorphism of two closed smooth subvarieties of $X$
such that $\dim Y_i $ does not exceed $\frac{m}{3}-1$. Then $\varphi$ extends to an automorphism of $X$.

\ecor

\section{Comparable morphisms}

\bdefi\label{cmo.d1} (1) Let $X$ (resp. $X'$) be a smooth algebraic variety and $\rho : X \to X'$ be
a morphism.  Consider a family $\cF$ of closed subvarieties of $X$.
Suppose that, given $Y_1$ and $Y_2 \in \cF$ 
isomorphic over $X'$ and such that $\rho|_{Y_i} : Y_i \to X'$ is a closed embedding,
there is an automorphism $\alpha \in \Aut (X/X')$ of $X$ over $X'$ that transforms $Y_1$
onto $Y_2$. In this case we say that $\rho$ is comparable on $\cF$. 

When the ground field $k$ is $\C$ we also say that $\rho$
is holomorphically comparable on $\cF$ if one can find a holomorphic automorphism of $X$ over $X'$
that transforms $Y_1$ onto $Y_2$.

(2) Let $X$ (resp. $X'$) be a smooth algebraic variety with a group $G\subset \Aut (X)$ (resp. $G'\subset \Aut (X')$) acting on it.
We say that a morphism $\rho : X \to X'$ is $(G,G')$-comparable if for every $g' \in G'$ there exists $g \in G$ such that $\rho \circ g = g' \circ \rho$. 

\edefi

The next fact follows from  Definition \ref{cmo.d1}.

\bprop\label{cmo.p1} Let $\rho : X \to X'$ be $(G,G')$-comparable. Let $Y_1$ and $Y_2$ be closed subvarieties of $X$ and $\varphi : Y_1 \to Y_2$
be an isomorphism. Suppose that each morphism $\rho|_{Y_i}: Y_i \to X'$ is a closed embedding, $Y_i'=\rho (Y_i)$, and $\varphi' : Y_1' \to Y_2'$ is the isomorphism for
which $\varphi' \circ \rho|_{Y_1} = \rho|_{Y_2} \circ \varphi$. Let  $\varphi'$ extend to an automorphism $g'$ of $X'$ which is an element of $G'$. Then there exists an element $g \in G$ such that $g(Y_1)$ is naturally isomorphic
to $Y_2$ over $X'$.
Furthermore, if $\rho$ is also comparable (resp. holomorphically comparable) on a family $\cF$ containing $Y_1$ and $Y_2$
then there exists an algebraic (resp. holomorphic) automorphism $\alpha$ of $X$
for which $\alpha |_{Y_1} =\varphi$.

\eprop

\bnota\label{cmo.n0} Given a smooth algebraic variety $X$ and an algebraic group $H$  consider the sheaf
generated by the presheaf whose elements are morphisms from Zariski open subsets of $X$ to $H$.
Denote by $\cH^1(X,H)$
the first $\check{\rm C}$ech cohomology of $X$ with coefficient in this sheaf. That is, for a Zariski open cover $\{ U_i \}$ of $X$ a one-cocycle
associated with this cover is given by morphisms $\psi_{ij} : U_i \cap U_j \to H$.
\enota

\bprop\label{cmo.p2} Let $H$ be an algebraic group and $H'$ be a closed algebraic subgroup of $H$.
Suppose that $H$ acts freely on a smooth algebraic variety $\tX$ so that the quotient map
$\tilde \rho: \tilde X \to \tX/H=: X'$ (resp. $\tau :\tX \to \tX/H'=:X$)  is a  principal  $H$-bundle (resp. $H'$-bundle)
over a smooth algebraic variety (and, hence,
for a Zariski locally trivial fiber bundle 
$\rho : X\to X'$  one has $\tilde \rho=\rho \circ \tau$). 
Let $h.x$ be the natural action of $h \in H$ on $x \in X$ and $\varphi : Y_1 \to Y_2$ be an isomorphism of closed subvarieties of $X$ over $X'$ such that $\rho|_{Y_k} : Y_k\to X'$ is
a closed embedding.  Suppose that $Y'=\rho (Y_k)$  and the group $\cH^1(Y', H')$ is trivial (which is 
the case when $H'$ itself is trivial). Then there exists a morphism $\theta : Y' \to H$ such that

{\rm (a)} $\varphi (y_1) = \theta (\rho (y_1)).y_1$ for every $y_1 \in Y_1$ and

{\rm (b)} $\varphi$ extends to an automorphism of $X$ over $X'$ if $\theta$ extends
to a morphism $X' \to H$.

Furthermore, if the ground field $\bk$ is $\C$ then

{\rm (c)} $\varphi$ extends to a holomorphic automorphism of $X$ over $X'$ if $\theta$ extends
to a holomorphic map $X' \to H$.

\eprop

\bproof  Let $y_2=\varphi (y_1)$, $y'=\rho (y_k)$ and $\tilde y_k \in \tau^{-1} (y_k)$. Since $y_1$ and $y_2$ are in $\rho^{-1}(y')\simeq H/H'$
we have an element $h \in H$ for which $y_2=h. y_1$  where $h$ is determined up to an element of $H'$. However, $h$ becomes unique
if we require additionally that $h$ transfer $\tilde y_1$ to $\tilde y_2$.  Since $\tau : \tX \to X$ is  a principal $H'$-bundle 
we can choose local sections of $Y_1$ (resp. $Y_2$) in $\tX$ which determines a local choice of $\tilde y_k$ as above.
Hence, for a Zariski open cover $\{ U_i\}$ of $Y'$ we can choose morphisms $\psi_i : U_i \to H$ for which
one has $\psi_i (y').y_1=y_2$ as soon as $y' \in U_i$.  If $y' \in U_i\cap U_j$ then letting $\psi_{ij}(y')=\psi_j(y')\circ \psi_i(y')^{-1}$
we get an element of $\cH^1(Y',H')$. Since the latter group is trivial we can construct $\theta : Y' \to H$ as in (a).
If $\theta$ extends to a morphism (reps. holomorphic map) $\Theta : X' \to H$ then the desired extension of $\varphi$
in (b) (resp. (c)) is given by $ \Theta (\rho (x)).x$ for $x \in X$.

\eproof

\bcor\label{cmo.c1} Let $\rho : X \to X'$, $H$ and $H'$ be as in Proposition \ref{cmo.p2} and 
$X'$ be quasi-affine.

{\rm (i)} Suppose that either $H'$ is trivial or $H'\simeq \SL(m,\bk)$ for some $m>0$.
Let $\cF_{\rm aff}$ be the family of closed subvarieties of $X$ isomorphic to $\A_\bk^k$
with $k=1,2, \ldots $.
Then $\rho$ is comparable on $\cF_{\rm aff}$.

{\rm (ii)} Suppose that $H'$ is trivial and  $H$ is isomorphic  as an algebraic variety to $\A_\bk^n$.
Consider the family $\cF$ of closed subvarieties $Z$ of $X$ such that
$\rho (Z)$ is closed in an affine variety $X''$ containing $X'$ as an open subvariety. Then $\rho$ is comparable on $\cF$.

\ecor

\bproof Let $Y_1$ and $Y_2\in \cF_{\rm aff} $ be as in Definition \ref{cmo.d1} (1), i.e., $Y'=\rho (Y_i)\simeq \A_\bk^k$ is closed in $X'$. 
 Note that by the Quillen-Suslin theorem \cite{Qui}, \cite{Sus} the group $\cH^1(Y',H')$ is trivial even in the case of $H'\simeq \SL(m,\bk)$\footnote
 {Formally, the  Quillen-Suslin theorem implies the group $\cH^1(\A^n_\bk ,{\rm GL}(m,\bk))$ is trivial. However,  each \v Cech
 cocycle on $\A^n_\bk$ with coefficients in ${\rm GL}(m,\bk)$ differs from a cocycle with coefficients in $\SL(m,\bk)$ by multiplication of, say, the first
coordinate of $\A^n_\bk$ by invertible regular functions. Since the first cohomology of $\A^n_\bk$ with coefficients in invertible regular functions is trivial
by the Serre Theorem B we see that $\cH^1(\A^n_\bk,{\rm SL}(m,\bk))$ is also trivial.},
 i.e., we are under the assumption of Proposition \ref{cmo.p2}.
Consider a morphism $g : X' \to \A_\bk^k$ whose restriction to $Y'$ is the identity map (such a $g$ exists since $X'$
is quasi-affine). Let $\theta$ be as in Proposition \ref{cmo.p2}. Then $\theta \circ g$ yields an extension of $\theta$ to $X'$ 
which by Proposition \ref{cmo.p2} concludes (i).

In (ii) note that since $Z':=\rho (Z), \, Z \in \cF$ is closed in the affine variety $X''$ this morphism $\theta$ extends to 
a morphism $\Theta : X'' \to \A_\bk^n \simeq H$ and again Proposition \ref{cmo.p2} yields (ii).

\eproof

\brem\label{cmo.r1}   Let  $\rho : X \to X'$ be an affine morphism of smooth varieties
and  $G\subset \Aut (X)$ (resp. $G'\subset \Aut (X')$).

(1) Note that if $\rho$ is $(G,G')$-comparable 
and $G'$ acts transitively on $X'$ then all fibers of $\rho$ are isomorphic and, furthermore, $\rho : X \to X'$
is a locally trivial fiber bundle
in \'etale topology (it follows from the fact that any affine morphism $W\to V$ of algebraic varieties with pairwise
isomorphic general fibers is locally trivial over an \'etale neighborhood of a general point of $V$ \cite{KrRu}).

(2) Let $G$ (resp. $G'$) be generated by a set $\cN$ (resp. $\cN'$) of complete algebraic vector fields on $X$ (resp. $X'$).
In order for $\rho$ to be $(G,G')$-comparable it suffices to require that for every $\delta' \in \cN'$ there exists $\delta \in \cN$ such that for every $x\in X$ and $x'=\rho (x)$
one has 
\be\label{cmo.eq1} \rho_* (\delta_x)=\delta_{x'}'\ee
where $\delta_x$ (resp. $\delta_{x'}'$) is the value of the field $\delta$ at $x$ (resp.  $\delta_{x'}'$ at $x'$).
\erem

\bdefi\label{cmo.d2} We call a pair $\delta, \delta'$ of vector fields (on $X$ and $X'$ respectively) comparable if they satisfy Formula (\ref{cmo.eq1}). 
Similarly, we call the pair  $(\cN, \cN')$  from Remark \ref{cmo.r1} (2) comparable
if for every $\delta' \in \cN'$ there exists $\delta \in \cN$ so that the pair $\delta, \delta'$ is comparable.

\edefi

\bexa\label{cmo.e1}
Let $X=\SL (n,\bk)$, thus $\dim X=n^2-1$. 
 Let $A=[a_{i,j}]_{i,j=1}^n$ be a matrix from $\SL(n,\bk)$
and let $k\leq n-1$ and $m\leq n$. Consider the $(m\times k)$-matrix $A'$ obtained from $A$ by removing all rows starting with $(m+1)$-st and all columns 
starting with $(k+1)$-st.
Then one has the natural morphism of $\rho :X\to X', \, A \to A'$ into the space $X'\simeq \A_\bk^{km}$ of $(k\times m)$-matrices.
Let  $1\leq i \ne j\leq n$ and  $$\delta_{ij}= \sum_{l=1}^n a_{l,i} \frac{\p}{a_{l,j}},$$
i.e., $\delta_{ij}$ is a locally nilpotent vector field on $X$ whose  flow is the addition of multiples of the $i$-th column in $A$ to the $j$-th column.
Note that $$\delta_{ij}'=\sum_{l=1}^m a_{l,i} \frac{\p}{a_{l,j}}, \,\, i,j \leq k$$
is a locally nilpotent vector field on $X'$ for which Formula  (\ref{cmo.eq1}) is valid.  
Similarly, the locally nilpotent vector fields
$\sigma_{ij}= \sum_{l=1}^n a_{i,l} \frac{\p}{a_{j,l}}$ and $\sigma_{ij}'= \sum_{l=1}^k a_{i,l} \frac{\p}{a_{j,l}}$ on $X$ and $X'$ respectively also satisfy
Formula  (\ref{cmo.eq1}).

\eexa

\bnota\label{cmo.n1} Suppose that $\rho : X \to X'$ is a smooth morphism of smooth algebraic varieties. 
Let $\cN$ (resp. $\cN'$) be a set of locally nilpotent vector fields on $X$ (resp. $X'$) such that the pair  $(\cN, \cN')$ is comparable.
Furthermore, we suppose that $X$ is a closed subvariety of $X' \times \A_\bk^m$ with $\rho$ being the restriction of the natural projection $\hat \rho : X' \times \A_\bk^m \to X'$.
That is, for every $\delta$ and $\delta'$ satisfying Formula  (\ref{cmo.eq1}) we can write $\delta = \delta' + \delta''$ where $\delta''$ is tangent to each fiber  $\hat \rho^{-1}(x')\simeq \A_\bk^m$ of $\hat \rho$
 (where by abuse of notation we identify $\delta'$  with its natural lift to $X' \times \A_\bk^m$).
\enota

Recall that a vector field on $\A_\bk^m$ (with a fixed coordinate system  $\bar u = (u_1, \ldots ,u_m)$)
is a linear  vector field if it is of the form $\dot {\bf x} =A {\bf x} + {\bf b}$ where {\bf x} and {\bf b} are vectors in $\A_\bk^m$ and $A$ is a square $(m \times m)$-matrix (i.e., it is a non-homogeneous system of linear differential equations).

\bprop\label{cmo.p3} 
Let Notation \ref{cmo.n1} hold. Suppose that for every comparable pair $(\delta, \delta') \in \cN \times \cN'$  
the following condition is true:

{\rm (A)} the restriction of the field $\delta''=\delta - \delta'$ to every fiber $\hat \rho^{-1}(x')\simeq \A_\bk^m$ is a linear vector field  (depending on $ x' \in X'$). 

Let $\tilde \cN'$ be the smallest saturated set of locally nilpotent vector fields on $X'$ 
such that it contains $\cN'$. 
Suppose that $\cL'$ is the Lie algebra generated by the fields of the form
$a \delta'$ where $\delta' \in \tilde \cN'$ and $a \in \bk [X']$. 
Then there exists a Lie algebra $\cL$ of vector fields on $X$ such that for every complete 
(resp. locally nilpotent) vector field $\sigma'\in \cL'$ there exists a complete (resp. locally nilpotent) 
vector field $\sigma \in \cL$ such that the pair $(\sigma, \sigma')$ is comparable. 

\eprop

\bproof 
Note that for every $a \in \bk [X']$ and every comparable
pair $(a \delta, a\delta')$ is comparable and satisfies condition (A).
Furthermore, let $(\delta_i, \delta_i'), \,  i=1,2, \ldots s$ be a collection of comparable pairs satisfying condition (A) and let $\ell$ be a linear form in $s$ variables. 
Let  $\kappa' = \ell (\delta_1', \ldots, \delta_s')$ (resp. $\kappa = \ell (\delta_1, \ldots, \delta_s)$). Then the pair $(\kappa, \kappa')$ is also comparable and satisfies 
condition (A). 

For the Lie brackets $\delta_0=[\delta_1, \delta_2]$ and  $\delta_0'=[\delta_1', \delta_2']$ the pair $(\delta_0, \delta_0')$ is, of course, comparable. Let us check condition (A) for this pair.
Note that $$[\delta_1, \delta_2] =[\delta_1', \delta_2'] + [\delta_1'', \delta_2'']+ [\delta_1', \delta_2'']+[\delta_1'', \delta_2']$$
where the first term belongs to $\cL'$ and the restriction of the second one to  every fiber  of  $\hat \rho$ is a linear vector field. Consider, say, the third term $ [\delta_1', \delta_2'']$.
Since the restriction of $\delta_2''$ to $\hat \rho^{-1}(x')\simeq \A_\bk^m$ is of the form $\dot {\bf x} =A {\bf x} + {\bf b}$ we see that
$ [\delta_1', \delta_2'']$ is of the form $\delta_1'(A) {\bf x} +\delta_1' ( {\bf b})$ and, in particular, it
is again a vector field whose restriction to  $\hat \rho^{-1}(x')$ is linear. Thus  the pair $(\delta_0, \delta_0')$ satisfies condition (A).

Now, a comparable pair $(\delta, \delta')$ of complete vector fields  induces the
flows $\varphi_t : X \to X$ and $\varphi_t': X' \to X'$ (where $t \in \bk$) so that $\rho\circ \varphi_t =\varphi_t'\circ \rho$.
Hence $\varphi_t$ tranforms a fiber $\rho^{-1}(x_0')\simeq \A_\bk^m$ to the fiber $\rho^{-1}(\varphi_t'(x_0'))\simeq \A_\bk^m$. 
When condition (A) is satisfied this map of fibers is an element of the flow of a linear  non-autonomous vector field
which is automatically an affine map. Hence for every comparable pair 
$(\sigma, \sigma')$ of complete vector fields  the conjugation by $\varphi_t$ and $\varphi_t'$ yields a comparable pair
$(\tilde \sigma, \tilde \sigma')$ of complete vector fields such that condition (A) is satisfied.   

Applying these operations of taking Lie brackets, conjugations, and linear combinations we construct 
the desired saturated set $\tilde \cN'$  
so that the pair $(\tilde \cN, \tilde \cN')$ is comparable for some set $\tilde \cN$ of locally nilpotent vector fields on $X$.
Let $\cL'$ (resp. $\cL$) be the Lie algebra generated by the fields of the form
$a \delta'$ where $\delta' \in \tilde \cN'$ and $a \in \bk [X']$ (resp. $a \delta$ where $\delta\in \tilde \cN$ and $a \in 
\bk [X']\subset \bk [X]$).
By construction for every
element $\sigma' \in \cL'$ there is an element of $\sigma \in \cL$ so that the pair $(\sigma, \sigma')$ is also comparable and satisfies condition (A).

Suppose now that $\sigma'$ is complete and $O'$ is an integral curve of this field  (in particular, $O'$ is isomorphic to either $\A_\bk$ or $\A_\bk^*$).
To show that $\sigma$ is complete 
it suffices to prove that $O'$ admits a lift to an integral curve $O\subset \hat \rho^{-1}(O')\simeq O'\times \A_\bk^m$
of the field $\sigma =\sigma' + \sigma''$. The restriction of $\sigma''$ to $O'\times \A_\bk^m$ is of the form $\dot {\bf x} =A {\bf x} + {\bf b}$ where the matrix $A$ and the vector ${\bf b}$
depend on the parameter $t \in O'$. That is, we are dealing with a non-autonomous system of linear equations.
Such a system has a solution for all values of $t$ which yields the desired lift of $O'$ to an integral curve of $\sigma$.

If $\sigma'$ is locally nilpotent then in order to show that $\sigma=\sigma'+\sigma''$ is locally nilpotent one needs to prove
that for every $b \in \bk [X]$ there exists $n$ for which $\sigma^n (b)=0$. This is true when $b \in \bk [X']\subset \bk [X]$
since  in this case $\sigma (b)=\sigma' (b)$. Since $\bk [X]$ is generated over $\bk [X']$ by the coordinate functions on $\A^m_\bk$
it suffices to prove that if  $b$ is a coordinate on $\A_\bk^m$ then $\sigma^n (b)=0$ for $n>>0$.
Note that in this case $\sigma (b) =\sigma'' (b)$ and by condition (A) one has $\sigma''(b) \in \bk [X']$. This yields the desired conclusion.
\eproof

\section{Holomorphic extension of $\theta$}

In this section the ground field $\bk$ is $\C$ and its aim is to describe conditions under which $\theta: Y' \to H$ from Proposition \ref{cmo.p2}
admits a holomorphic extension $\Theta : X' \to H$.

\bprop\label{hex.p1} Let the assumption of Proposition \ref{cmo.p2} hold and $X'$ be Stein.
Suppose that the map $\theta$ is homotopy equivalent to a constant map
from $Y'$ to $H$. Then $\theta$ admits  a holomorphic extension $\Theta : X' \to H$.

\eprop

\bproof By the Oka-Grauert principle \cite[Theorem 5.4.4]{For} it suffices to construct a continuous extension of $\theta$.
Let $\tilde \theta (t) : Y \to H, \, t \in [0,1]$ be a homotopy of $\theta$ to a constant map, i.e., $\tilde \theta (*,0)=\theta$ and 
$\theta (*,1)$
sends $Y'$ to a point $h_0 \in H$.

The argument is rather simple in the case of a smooth $Y'$ since 
one can consider its closed tubular neighborhood $U$. Note that $U$ admits a continuous map $g: U\to [0,1]$
such that $g^{-1}(0)=Y'$ and $g^{-1}(1)$ is the boundary $\p U$ of $U$. Define a continuous extension $\theta' : U \to H$
of $\theta$ via the formula $\theta' (u)= \tilde \theta  ( p(u), g(u)) $ where $u \in U$ and $p: U \to Y'$ is
the natural projection. Then we can extend $\theta'$ further to $\theta'' : X' \to H$
by putting $\Theta(x') =h_0$ for $x' \in X' \setminus U$.

In the general case we recall that complex algebraic varieties admit triangulation \cite{Lo}.
That is, one can view $Y'$ as a subcomplex in a complex $X'$. Then the second baricentric derived neighborhood
$U$ of $Y'$ is regular (e.g., see \cite{Hir}). This means, that $Y'$ is a deformation retract of $U$ and
the map $U \to Y'$ (which is identity on $Y'$)
can be constructed via a sequence of collapses $U_{i-1} \searrow U_{i}, \, i=1, \ldots , k $ where each $U_i$ is a subcomplex of
the second baricentric subdivision, $U_0=U$, $U_k$ coincides with the second baricentric subdivision of $Y'$, 
and $U_{i-1} \setminus U_{i} =\{ \sigma_i, \tau_i)  $ where
$\sigma_i$ is an $l$-dimensional simplex 
 and $\tau_i$ is the only $(l-1)$-dimensional face of $\sigma_i$ that is not contained
in $U_{i}$ (viewed as a complex). In particular,  the union $\bigcup_{i=1}^k\tau_i$
contains the boundary of $U_0$. Let us establish the following.

{\em Claim.} A continuous map
$\theta_i: U_i \to H$, admitting a homotopy $\tilde \theta_i$ to the constant map $h_0$ and such that $\theta (\p \tau_i)=h_0$, can be extended to a similar continuous map
$\theta_{i-1}: U_{i-1} \to H$ with a homotopy $\tilde \theta_{i-1}$ to  the constant map $h_0$ such that  $\tilde \theta_{i-1} (\tau_i)=h_0$.

Let $\p \sigma_i$ (resp. $\p \tau_i$) be the boundary of $\sigma_i$ (resp. $\tau_i$). 
Then there is a homotopy $g: \sigma_i \times [0,1]\to \sigma_i$ such
that $g(*,1)=\id|_{\sigma_i}$, ${\rm Im} \, g(*, 0)= \p \sigma_i \cap U_{i-1}$, and
for every $t \in [0,1]$ the restriction $g(*,t)|_{\p \sigma_i \cap U_{i-1}}$ is the identity
map on $\p \sigma_i \cap U_{i-1}$. 
In particular,
each collapse $U_{i-1} \searrow U_{i}$
induces a strong deformation retract $f_i :U_{i-1} \to U_i$.
Furthermore, $g$ can be chosen so that it yields 
a natural homeomorphism between $\sigma_i \setminus \p \tau_i$ and $\dot{\tau}_i\times [0,1]$
where $\dot {\tau_i}=\tau_i \setminus \p \tau_i$ is the interior of $\tau_i$. 
\footnote{In order to see this, treat $\sigma_i$ as 
a closed $l$-dimensional ball, $\tau_i$ as the upper $(l-1)$-dimensional semi-sphere in its boundary,
$\p \sigma_i \cap U_{i-1}$ as the lower semi-sphere and $\p \tau_i$ as its equator.}
Hence we can now define the extension $\theta_{i-1}$ of $\theta_i$ so that every $u \in \dot \tau_i$ and $t \in [0,1]$ 
one has $\theta_{i-1} (g(u,t))=\tilde \theta_i (f_i(u),t)$. By construction $\theta_{i-1} (\tau_i)=h_0$ and the extension is
homotopic to a constant map (indeed, define  the homotopy $\tilde  \theta_{i-1}(g(u,s))$  on $U_{i-1} \setminus \p \tau_i$ via
$\tilde \theta_i (f_i(u),t+s (1-t))$). This yields the the Claim.

It remains to show that for such extensions one has $\tilde \theta_i|_{\p U}= h_0$ as soon as $U_i \cap \p U\ne \emptyset$. First
note that if $\sigma_i$ is one-dimensional and meets $Y'$ then $\tau_i$ is a singleton (i.e., $\p \tau_i =\emptyset$) and, therefore, $\tilde \theta_i (\tau_i)=h_0$.
Let  $S$ be the collection of $\sigma_i$ such that each of them has a face contained in $U_i$. For every $\sigma_i \in S$
there is the only one vertex not in $U_i$. Suppose that $T$ is the set of such vertices. Since each of $\sigma_i \in S$ collapses before
the  one-dimensional simplexes mentioned above we can suppose that 
$\tilde \theta_i|_T=h_0$. Let $S_m$ be the collection of $m$-dimensional simplexes $\sigma_j\notin S$ that have faces in $\bigcup_{\sigma_i \in S}\sigma_i$.
For every $\sigma_j \in S_2$ we have a one-dimensional $\tau_j$ such that $\p \tau_j\in T$. Thus by the  Claim
we can suppose that the restriction of $\tilde \theta_j$ to $\bigcup_{\sigma_j \in S_2}\tau_j$ is the constant map to $h_0$.
Similarly, for every $\sigma_k \in S_3$ we have $\p \tau_k$ in $\bigcup_{\sigma_j \in S_2}\tau_j$. Hence, proceeding by induction
we see that the restriction  of every $\tilde \theta_i$ to $\bigcup_{j\geq 2} \bigcup_{\sigma_k \in S_j}\tau_k$ is $h_0$.
Consequently, for every $\sigma_i \notin S \cup\bigcup_{j\geq 2} S_j$ we have  $\tilde \theta_i ( {\tau_i})= h_0$. 
Since we deal with the second baricentric subdivision  we have $\bigcup_{\sigma_i \in S} \tau_i \cap \p U=\emptyset$ which
yields the desired conlusion.
\eproof

\blem\label{hex.l2} Let $m\geq 2$ and $Z$ be a complex affine algebraic variety such that $H_i(Z)=0$ for
$i \geq m+1$. Then $Z$ can be embedded into a contractible topological space $\hZ$
such that $\hZ \setminus Z$ is a union of components each of which is homeomorphic to a ball $B$
whose real dimension is between 2 and $m+2$.

\elem

\bproof Since any algebraic variety is a finite CW-complex, $Z$ has a finitely generated fundamental group.
Gluing $Z$ with discs along generators of this group we obtain a topological space $Z_1 \supset Z$ such that
$Z_1$ is simply connected. Furthermore, the Mayer-Vietoris sequence implies that the embedding
$Z \hookrightarrow Z_1$ induces an isomorphism $H_i(Z) \simeq H_i(Z_1)$ for $i\geq 3$.
By the Hurewicz theorem $\pi_2(Z_1)$ is naturally isomorphic
to $H_2 (Z_1)$.  Thus, we can choose a nontrivial element $\omega$ of $H_2(Z_1)$ presented by
a two-dimensional cell in the CW-complex which is homeomorphic to a two-sphere.  
 Glue $Z_1$ with a three-dimensional
ball along this sphere and let $V_1 =Z_1\cup B_1$ be the resulting simply connected topological space.  
Consider the Mayer-Vietoris sequence 
$$\to H_{k+1}(V_1)\overset{d_{k+1}}\to H_k(Z_1\cap B_1) \to H_k (Z_1)\oplus H_k (B_1) \to H_k (V_1) \overset{d_k}  \to H_{k-1}(Z_1\cap B_1)
\to  .$$
Note that if $\omega$ generates a free group then the map $H_2(Z_1\cap B_1) \to H_2 (Z_1)$ is a monomorpism and,
therefore, $d_3$ sends $H_3(V_1)$ to zero. That is, after such a gluing $H_i (V_1)=H_i(Z_1)$ for $i \geq 3$
while the rank of $H_2(V_1)$ is less than the rank of $H_2(Z_1)$ since $\omega$ induces the zero element in $H_2(V_1)$.
Continuing this procedure we can embed $Z_1$ into a simply connected $V_2$ such that $H_i (V_2)=H_i(Z_1)$ for $i \geq 3$
while $H_2(V_2)$ is finite. Let $V_3$ be the result of gluing a three-ball  to $V_2$ along $\omega$ as before and 
$V_3=V_2 \cup B_2$ be the resulting simply connected topological space.
 In this case $\omega$ is of finite order $l$ and the kernel of $H_2(V_2\cap B_2) \to H_2 (V_2)$ in the associated
 Mayer-Vietoris sequence is isomorphic
 to $l \Z \subset \Z \simeq H_2(V_2\cap B_2)$. In particular, the image of the
 map $H_{k+1}(V_3)\to H_k(V_2\cap B_2)\simeq  \Z$ is $l\Z$. Since the latter is a free $\Z$-module the map $H_{k+1}(V_3)\to l\Z$
 has a right inverse. Thus it follows from the Mayer-Vietoris sequence
 that $H_3(V_3)$ is naturally isomorphic to $H_3(V_2) \oplus \Z$ while $H_i (V_3)=H_i(V_2)$ for $i \geq 4$.
 Since the number of elements in $H_2(V_3)$ is less than the number of elements in $H_2(V_2)$ continuing this
 procedure we embed $Z_1$ into a simply connected $Z_2$ such that $H_2(Z_2)=0$, 
 $H_i (Z_2)=H_i(Z_1)=H_i (Z)$ for $i \geq 4$, and $H_3(Z_2)$ is naturally isomorphic 
 to the direct sum of $H_3(Z)$
 and a finitely generated free $\Z$-module.

In the same manner we can embed $Z$ into a simply connected $Z_m$ such that $H_i(Z_m)=0$ for $i \leq m$, 
 $H_i (Z_m)=H_i(Z)=0$ for $i \geq m+2$, and $H_{m+1}(Z_m)$ is naturally isomorphic to the direct sum of $H_{m+1}(Z)$
 and a finitely generated free $\Z$-module. That is, $H_{m+1}(Z_m)$ is a finitely generated free $\Z$-module
 since $H_{m+1}(Z)=0$. Gluing $Z_m$ with balls of real dimension $m+2$ we can obtain
 a simply connected $\hZ$ with $H_i(\hZ)=0$ for $i \leq m+1$. Furthermore, as we saw before, 
 since $H_{m+1}(Z_m)$ is free,
 this procedure does not affect the $i$-homology groups with $i\geq m+2$.
Hence $H_i (\hZ) =0$ for $i\geq 1$ and
by the Hurewicz theorem $\pi_i(\hZ)=0$ for $i \geq 1$. By the Whitehead theorem $\hZ$ is contractible which
yields the desired conclusion. 
\eproof

\brem\label{hex.r1}
It follows from the proof that if  in Lemma \ref{hex.l2} the group $H_m (Z)$ is free then in order  to construct $\hZ$
one needs only balls of real dimension at most $m+1$.
\erem

\bprop\label{hex.p2} Let $\rho : X \to X'$ be a principal bundle for some algebraic group $H$ and
$Y'$ be a closed subvariety of $X'$.
Let $X'$ be Stein and  $H$ 
simply connected with $H_i (H)=0$ for $i\leq m$ where $m\geq 2$.
Suppose also that $H_i(Y')=0$ for $i \geq m+1$ and
$H_m (Y')$ is free. Then any morphism $\theta : Y' \to H$ admits a holomorphic extension $\Theta : X' \to H$.

\eprop

\bproof By Lemma \ref{hex.l2} we can embed $Y'$ into a contractible topological space $Y$.
Furthermore, by Remark \ref{hex.r1} we can suppose that each component of $Y \setminus Y'$
is a ball $B$ whose real dimension $k$ is at most $m+1$
and whose boundary is an image of a $(k-1)$-dimensional sphere
in $Y'$. Composition with $\theta$ yields a continuous map of this sphere to $H$.
Since $\pi_{k-1}(H)=0$ we can extend any continuous map of 
a $(k-1)$-dimensional sphere in $H$ to a map from a closed $k$-ball. Hence 
we can extend the morphism $\theta : Y' \to H$ to a continuous 
map $\hat \theta : Y \to H$. Let $\psi : Y \times [0,1] \to Y$ be a contraction of $Y$, i.e., $\psi (*, 0) =\id|_Y$
and $\psi (*,1)=y_o \in Y$. Then $\tilde \theta = \hat \theta \circ \psi$ is a homotopy of $\theta$ to 
a constant map. Hence we are done by Proposition \ref{hex.p1}.

\eproof

\section{The case of $\SL(n,\bk)$}

\bnota\label{qub.n1}
Let us fix notation for the rest of the paper. From now on $X$ will be always an affine algebraic variety isomorphic to $\SL (n,\bk)$ with $n\geq 3$, i.e.,  $\dim X=n^2-1$. 
We treat points in $X$ as matrices $A=[a_{i,j}]_{i,j=1}^n$ from $\SL(n,\bk)$ and denote by $\{ A_{ij} \}$ the cofactors of this matrix.
For $m\leq n$ let $P_m$ be the space of $(m\times n)$-matrices,
$P_m^0$ be its subvariety consisting of matrices of rank $m$ and
$\rho_{m} : X \to P_{m}$ be the natural projection that sends $A$ to the matrix consisting of the first $m$ rows of $A$. 
Consider also the natural projection $\rho: X\to Q$ where $Q$ is the quadric given by the equation
 \be\label{csl.eq1} a_{11}A_{11}+a_{12}A_{12} + \ldots + a_{1n}A_{1n}=1\ee
in the affine space $\A_\bk^{2n}$ equipped with coordinates $(a_{1,1}, \ldots , a_{1,n}, A_{1,1}, \ldots , A_{1,n})$.
Suppose that
the group $H=\SL (n,\bk)\times \SL(n,\bk)$ acts on  $X$ so that $(B, C) \in \SL (n,\bk)\times \SL(n,\bk)$
sends $A \in X$ to $BAC$. Note that $H$ contains the subgroup $H_m=\SL (m,\bk) \times \SL (n,\bk)$
that acts naturally on $P_m$, while $H_1$ acts naturally on $Q$. 
Under these action the morphism $\rho_{m}$ (resp. $\rho$) is $H_m$-equivariant (resp. $H_1$-equivariant).
\enota

\blem\label{qub.l1} Let Notation \ref{qub.n1} hold. Then the  $H_m$-action (resp. $H_1$-action) on $P_{m}^0$ (resp. $Q$)
is transitive.
\elem

\bproof  Since any finite sequence of column and row operations on a matrix $D\in P_m$ coincides with
the action of some element of $H_m$ on this matrix we see that for  $D\in P_m^0$
there is an element $h \in H_m$ for which $h. D =[d_{i,j}]$
where $d_{i,i}=1$ for every $i=1, \ldots, m$, while $d_{i,j}=0$ for $i \ne j$. This yields the transitivity
of the $H_m$-action on $P_m^0$. 

Similarly, an elementary column operation induced by an element of $H_1$ yields an automorphism
of $Q$ such that
\be\label{csl.eq2q} a_{1, j}\to a_{1,j} + ta_{1,i} \text{ and } A_{1,i}\to A_{1,i}-t A_{1,j} \ee
 (where $t\in \bk$ is a multiple),
while it keeps the rest of coordinates the same.
Since $a_{1,1}, \ldots, a_{1,n}$ cannot vanish simultaneously we can send an arbitrary point $q$ in $Q$ via such automorphism to a point $q_1$ with
\be\label{csl.eq3q}  a_{1,1}=1 \text{ and } a_{1,2}=\ldots =a_{1,n}=0. \ee 
In this case $A_{1,1}=1$ because of Formula (\ref{csl.eq1}). 
Let us use now the similar automorphisms which send $a_{1,1} \to a_{1,1} + ta_{1,j}$
and $A_{1,j} \to A_{1,j} -t A_{1,1}$, while preserving the rest of coordinates on $Q$. 
Then we can send $q_1$ to a point $q_0$ with $A_{1,2}= \ldots = A_{1,n}=0$, 
while keeping Formula (\ref{csl.eq3q}) valid. 
That is, we can send an arbitrary point $q$ to this given point $q_0$ which yields the transitivity of the $H_1$-action on $Q$
and the desired conclusion. \eproof

\bnota\label{qub.n2} Let $I_k$ be the identity $(k\times k)$-matrix for $k\geq 1$. Consider  the following subgroups of block matrices
in $\SL(n,\bk)$.
\[
F_{n-m}'=\left[ {\begin{array}{cc}
   I_m & \bar 0 \\      
   M& \SL (n-m,\bk)  \end{array} }    \right], 
   \text{ resp. } 
F''=\left[ {\begin{array}{cc}
   I_1 & \bar 0 \\      
 \bar 0 & \SL (n-1,\bk)  \end{array} } \right]  \]
 where $M$ is the set of $(n-m)\times m$ matrices.
Then we have the subgroup $H_{n-m}'=F_{n-m}'\times I_n$ of $H$ acting on $X$ and the subgroup 
$H''=F''\times I_n$ acting on $X$.

\enota

\bprop\label{qub.p1} Let Notations \ref{qub.n1} and \ref{qub.n2} hold and $n-m\geq 2$. 
Then 

{\rm (i)} $\rho_{m} : X \to P_m^0$ is a principal $H_{n-m}'$-bundle over $P_{m}^0$;

{\rm (ii)} $\rho : X \to Q$ is a principal $H''$-bundle over $Q$.
\eprop

\bproof 
Note that the action of $H_{n-m}'$ on $X$ is free and it preserves every fiber of the morphism $\rho_m$.
Therefore,   the fibers of $\rho_{m}$ are of dimension at least 
$\dim H_{n-m}'=m(n-m)+(n-m)^2-1=n(n-m)-1$. Furthermore, the action of $H_m$ transforms each orbit of $H_{n-m}'$ into 
another orbit.
Hence, the fibers  of $\rho_m$ are of the same dimension because $\rho_{m}$ 
is $H_m$-equivariant and the $H_m$-action on  $P_m^0$  is transitive.
Since $\dim P_m^0=nm$ and $\dim X =n^2-1$,
observing  the equality   and $ n^2-1=nm +n(n-m)-1$ we see that the fibers of
$\rho_m$ are nothing but the orbits of $H_{n-m}''$ which shows that
$\rho_{m} : X \to P_m^0$ is a principal $H_{n-m}'$-bundle. Thus we have (i).

Similarly, the $H''$-action on $Q$ is free, it preserves every fiber of the morphism $\rho$ and it commutes with the $H_1$-action from Notation \ref{qub.n1}.
Hence, the fibers of $\rho$ are of dimension at least 
$\dim \SL (n-1,\bk)=(n-1)^2-1$. All these fibers are again of the same dimension because $\rho$
is $H_1$-equivariant and the $H_1$-action on $Q$ is transitive by Lemma \ref{qub.l1}.
On the other hand, the dimension of $Q$ is  $2n- 1$.
Observing  the equality $\dim \SL (n,\bk) = n^2-1=(n-1)^2-1 +(2n-1)$ we conclude that the dimension of each fiber is $(n-1)^2-1$.
That is, the fibers of $\rho$
are nothing but the orbits of the $H''$-action which yields  (ii) and the desired conclusion.
\eproof

\bnota\label{qub.n3} (1) Let $\cN=\{ \delta_{i,j} | 1\leq i \ne j\leq n \} \bigcup \{ \sigma_{i,j} | 1\leq i \ne j\leq n \} $ be the set of locally nilpotent vector fields on $X$ given by
$$\delta_{i,j}= \sum_{l=1}^n a_{l,i} \frac{\p}{a_{l,j}} \, \text{  and   } \,  \sigma_{i,j}= \sum_{l=1}^n a_{i,l} \frac{\p}{a_{j,l}},$$ 
$\cN_m'=\{ \delta_{i,j}' | 1\leq i \ne j\leq n \} \bigcup \{ \sigma_{i,j}' | 1\leq i \ne j\leq m \} $ be the set of locally nilpotent vector fields on $P_m$ given by
$$\delta_{i,j}'= \sum_{l=1}^m a_{l,i} \frac{\p}{a_{l,j}}  \, \text{  and   } \,  \sigma_{i,j}'= \sum_{l=1}^n a_{il} \frac{\p}{a_{j,l}},$$
and $\cN''=\{ \delta_{i,j}'' | 1\leq i \ne j\leq n \} $ be the set of locally nilpotent vector fields on $Q$ given by
$$\delta_{i,j}''= a_{1,i} \frac{\p}{a_{1,j}} -A_{1,j} \frac{\p}{A_{1,i}}.$$

(2) Let $\tilde \cN_m'$ (resp. $\tilde \cN''$) be the smallest saturated set of locally nilpotent vector fields on $P_m$ containing $\cN'$
(resp. on $Q$ containing $\cN''$).
In particular,  $\tilde \cN_m'$ generates a group  $G_m'$ of automorphisms of of $P_m$ and $\tilde \cN''$ generates a
group $G''$ of automorphisms of $Q$.
\enota

\blem\label{qub.l2}  Let Notation \ref{qub.n3} hold,  $n\geq 3$ and $1\leq m\leq n-2$.

{\rm (1)}  The pair $(\cN, \cN_m')$  (resp. $(\cN, \cN'')$) is comparable 
for the morphism $\rho_m$ (resp. $\rho$).

{\rm (2)} Every element of $(\cN, \cN_m')$ (resp. $(\cN, \cN'')$) satisfies Condition (A) in 
Proposition \ref{cmo.p3}.

{\rm (3)} The variety $P_m^0$ is $G_m'$-flexible and the variety $Q$ is $G''$-flexible.

{\rm (4)} Let $\cL_m'$ be the Lie algebra of vector fields on $P_m$ generated
by all fields of the form $a'\delta'$ where $\delta'\in \tilde \cN_m'$ and $\deg_{\delta'}a' \leq 1$.
Then $\cL_m'$ contains the space of all algebraic vector fields on $P_m$ that vanish at $Z=P_m\setminus P_m^0$
with some multiplicity $s>0$.

{\rm (5)} Let $\cL''$ be the Lie algebra of vector fields on $Q$ generated
by all fields of the form $a''\delta''$ where $\delta''\in \tilde \cN''$ and $\deg_{\delta''}a'' \leq 1$.
Then $\cL''$ coincides with the space of all algebraic vector fields on $Q$.
\elem

\bproof The first and second statements  follow from Notation \ref{qub.n3} and Definition \ref{cmo.d2}.
For (3) note that all elementary row and column operations on $P_m$ can be viewed as elements of flows
of fields from $\cN_m'$. Since such operations generate the action of the group $H_m$ we see that $G_m'$ contains
$H_m$. By the similar reason $G''$ contains $H''$. Since by Lemma \ref{qub.l1} $H_m$ acts transitively on $P_m^0$ (resp. 
$H_1$ acts transitively on $Q$) we have (3) by Theorem \ref{fm.t1a}.

The vector fields $\delta_{1,2}'$ and $\delta_{3,2}'$ commute,
i.e., we have a compatible pair of locally nilpotent vector fields in $\tilde \cN_m'$. Since $P_m^0=P_m \setminus Z$
is $G_m'$-flexible Theorem \ref{adp.t1} implies that
the Lie algebra $\cL_m'$ contains all algebraic vector fields $\AVF (P_m)$ on $P_m$
vanishing on $Z$ with some multiplicity $s>0$, i.e., we have (4).

The same argument and the commutativity of $\delta_{1,2}''$ and $\delta_{3,2}''$ imply (5).
\eproof

\bprop\label{qub.p2} Let $\tG=\SAut (X)$, $\tG''=\SAut (Q)$ and let
$\tilde G_m'= \SAut_{Z_s} (P_m)$ be the subgroup of $\SAut (P_m)$ generated  by the elements of the flows of locally nilpotent vector fields
on $P_m$ whose restriction to the $s$-infinitesimal neighborhood $Z_s$ is zero, where $s$ is as in Lemma \ref{qub.l2}. 
Then $\rho_m : X \to P_m^0$ is  a $(\tG,\tG_m')$-comparable  morphism and $\rho : X \to Q$
is a $(\tG,\tG'')$-comparable morphism.
\eprop

\bproof By Lemma \ref{qub.l2} (4) the group $\tG_m'$ is generated by elements of the flows of locally nilpotent
vector fields from $\cL_m'$. Since Conditions (A) is satisfied by Lemma \ref{qub.l2} (2) we see that
$ \rho_m$ is $(\tG,\tG_m')$-comparable by Proposition \ref{cmo.p3}. Similarly, by Lemma \ref{qub.l2}
the group $\tG''$ is generated by elements of the flows of locally nilpotent vector fields from $\cL''$ and 
and the same reasoning yields the second statement.
\eproof

\blem\label{csl.l1a}  Let Notation \ref{qub.n1} hold. Then the codimension of the subvariety $Z=P_m\setminus P_m^0$ in $P_m$ is $n-m+1$.
\elem

\bproof A matrix $A\in P_m$ does not belong to $P_m^0$ if and only if its rank $k$ is at most $m-1$.
Applying some element of $H_m$ to $A=[a_{i,j}]$ with $k\leq m-1$ we can suppose that $A=[a_{ij}]$ has
$a_{i,j}=1$ where $i=j$ and $i \leq k$, while $a_{i,j}=0$ in all other cases. Note that the isotropy subgroup $F$ of $A$ in $H_m$
consists of all elements $(B,C)$ such that $B$ is a block $(m\times m)$-matrix of the form

\[
  \left[ {\begin{array}{cc}
   A' &  B' \\      
   \bar 0 &D'    \end{array} } \right]
\]

and  $C$ is a block $(n\times n)$-matrix of the form

\[
  \left[ {\begin{array}{cc}
   (A')^{-1} & \bar 0 \\      
    C' & E'    \end{array} } \right]
\]
where $A'\in  \GL (k,\bk) $, $D' \in \GL (m-k,\bk)$, $E' \in \GL (n-k,\bk)$  with $\det D'=\frac {1}{\det A'}$ and $\det E'=\det A'$, while $B'$ (resp. $C'$) is an aribrary
matrix of size $k\times (m-k)$ (resp. $(n-k)\times k$). Hence, the dimension of $F$
is $k(m-k)+ (m-k)^2 +(n-k)^2 +k(n-k)+k^2-1$.   Since the subvariety of matrices of rank $k$ in $P_m$ is isomorphic to $H_m/F$
we see that its dimension is $\dim H_m-\dim F=m^2-1+ n^2-1-(k(m-k)+ (m-k)^2 +(n-k)^2 +k(n-k)+k^2-1)=m(m-k)+ k(n-k)-k^2-1$. For $k$ running from 1 to $ m-1$ 
the maximum of latter expression is achieved for $k=m-1$ and it is equal to
$m+(m-1)(n-m+1)-(m-1)^2-1=(n+1)(m-1)=nm-n+m-1$. Since $\dim P_m=nm$ we see that the codimension of $Z$ in $P_m$ is $n-m+1$ which is
the desired conclusion.
 \eproof
 
 \blem\label{csl.l2a}  Let Notation \ref{qub.n1} hold, $0<m<n$ and
$Y$ be a closed subvariety of $X$ of dimension at most $m$.
Suppose that $\brH \simeq \SL (n, \bk)$ acts on $X$ via left multiplications. 
Then for a general element $h \in \brH$ the morphism $\rho_m|_{h(Y)}: h(Y) \to P_m$ is proper.

\elem

\bproof  Suppose that $\cX=\A^{n^2}_\bk$ is equipped with the coordinates system $(a_{1,1}, \ldots , a_{n,n})$ and, thus, $\cX$ contains $X$
as a closed subvariety.
Present $\cX$ as $\cX =\prod_{i=1}^n \cX_{i0}$ where $\cX_{i0}=\A^n_\bk$ has the coordinate system $(a_{1,i}, \ldots , a_{n,i})$.
Let $\bar \cX_i \simeq \PP^n$ be a completion of $\cX_{i0}$, i.e., $\bar \cX_i=\cX_{i0}\sqcup\cX_{i1}$ (where $\cX_{i1}\simeq \PP^{n-1}$)
and $\bar \cX =\prod_{i=1}^n\bar \cX_i$ is a completion of $\cX$.
Note that the $\brH$-action on $X$ extends naturally to an $\brH$-action on $\cX$ (resp. $\bar \cX$) and we have also 
the transitive $\brH$-action on each $\cX_{i0}$ with an extension to an action on $\bar \cX_i$ for which $\cX_{i1}$ is an $\brH$-orbit and
the natural projection $\cX \to \cX_i$ (resp. $\psi_i: \bar \cX \to \bar \cX_i$) is  $\brH$-equivariant. Let $J$ be the set $\{0,1\}^n$ without
the element $(0, \ldots, 0)$. Note that the boundary $\bar \cX =\bar \cX \setminus \cX$ can be presented as
$\bigcup_{\bar j \in J} C_{\bar j}$ where $C_{\bar j} =\prod_{i=1}^n \cX_{ij_i}$ and $\bar j =(j_1, \ldots , j_n)\in J$.
In particular, for the closure $\bY$ of $Y$ in $\bar \cX$ one has $\bar Y\setminus Y=\bigcup \brY_{\bar j}$ where $\brY_{\bar j}
=(\bar Y\setminus Y)\cap C_{\bar j}$. Let, say, $j_1=1$ for some $\bar j' \in J$ and, hence,  $\psi_1( C_{\bar j'})\subset  \bX_{11}$.
Since $\bY\setminus Y$ is of dimension $m-1$ we see that $\psi_1 (\brY_{\bar j'})$ is of dimension at most $m-1$. 
Note that the set  $\cR (a_{1,1}, \ldots, a_{m,1})$ of common indeterminacy points of the functions $a_{1,1}, \ldots, a_{m,1}$ on $\bar \cX_1$
is a subset of $\cX_{11}$ of codimension $m$. By Theorem \ref{agga.t1} for a general $h \in \brH$ the intersection
of $h(\psi_1 (\brY_{\bar j'}))$ with $\cR (a_{1,1}, \ldots, a_{m,1})$ is empty. Since $\psi_1$ is $H$-equivariant we have
$h(\brY_{\bar j'})\cap \psi_1^{-1}(\cR (a_{11}, \ldots, a_{m,1}))=\emptyset$ where $\psi_1^{-1}(\cR (a_{1,1}, \ldots, a_{m,1}))$
is the intersection of $C_{\bar j'}$ with  the set of common indeterminacy points of $a_{11}, \ldots, a_{m,1}$ in $\bar \cX$.
Note that the similar
claims are true for all $\bar j \in J$. Hence, for a general $h \in \brH$ the variety $h(\bY)$ does not meet 
the set of common indeterminacy points of the functions $a_{1,1}, \ldots, a_{m1}, a_{1,2}, \ldots, a_{m,n}$.
Furthermore,  the boundary $\tilde \cX \setminus \cX$ does not contain points where these functions are regular
and take finite values.
Therefore, the desired conclusion follows now from Proposition \ref{gp2.p1}.
\eproof

\bprop\label{csl.p1a} Let the assumptions of Lemma \ref{csl.l2a} hold and $\ED (Y)\leq m$ where $m \leq n-2$. Then
for some element $\beta \in \SAut (X)$ the morphism $\rho_m|_{\beta(Y)}: \beta (Y) \to P_m$ is a closed embedding.
\eprop

\bproof  By Proposition \ref{qub.p1} we have $X/H_{n-m}'\simeq P_m^0$. 
Note that the $H_{n-m}'$-action on $X$ is generated by elements of the flows of some set $\cS$ of locally nilpotent vector fields
on $X$.  Let $\tilde \cS$ be the smallest saturated set of locally nilpotent 
vector fields on $X$ containing $\cS$ and tangent to the fibers of $\rho_m$. 
Then $\tilde \cS$ generates a group $F$ of automorphisms of $X$ over $P_m$ such that it contains $H_{n-m}'$
and every fiber of $\rho_m$ is $F$-flexible. By Theorem \ref{gp1.t1} (iv)
there exists an algebraic family $\cA \subset F$ such that for a general element $\alpha \in \cA$
the morphism $\rho_m\circ \alpha : Y\to P_m$ is injective and it induces an injective map
of tangent bundles. By Proposition \ref{gp1.p3} the same is true for a morphism $ (\rho_m\circ h) \circ \alpha : Y\to P_m$ 
where $h$ is a general element of $\brH$. By Lemma \ref{csl.l2a} the latter morphism is proper.  Thus, letting $\beta =h \circ \alpha$ we get
the desired conclusion. 
\eproof

\bthm\label{csl.t2} Let $X=\SL (n,\C)$ and
$\varphi : Y_1\to Y_2$ be an isomorphism of two closed subvarieties of $X$  such that either

{\rm (i)} $\ED (Y_i)+\dim Y_i \leq n-2$, $H_i (Y_1)=0$ for $i \geq 3$ 
and $H_2 (Y_1)$ is a free abelian group; or

{\rm (ii)} $\dim Y_1$ is a curve and $\ED (Y_i)\leq n-2$, or;

{\rm (iii)} $Y_1$ is a once-punctured curve and $\ED (Y_1) \leq 2n-3$. 

Then there exists a holomorphic automorphism $\beta$ of $X$ such that  $\beta|_{Y_1} =\varphi$.
\ethm

\bproof  Let $m=\ED (Y_i)$. By
Proposition \ref{csl.p1a} we can suppose that $\rho_m|_{Y_i}: Y_i \to P_m$ is a closed embedding.
Let  $\varphi' : Y_1' \to Y_2'$ be the isomorphism
induced by $\varphi$, where $Y_i'= \rho_m (Y_i)$. 
By Lemma \ref{csl.l1a} one has $\dim Z=nm-(n-m+1)$, where $Z=P_m \setminus P_m^0$. Therefore,
the assumption $m+\dim Y_i \leq n-2$ implies that
$\dim Y_i' + \dim Z\leq \dim P_m -3$ in case (i). 
Hence, by Theorem \ref{cgw.t1} (a) and Proposition \ref{cgw.p3} $\varphi'$ can be extended
to an automorphism $\alpha' \in \SAut_{Z_s} (P_m)$ for any $s>1$. By Proposition \ref{qub.p2}  $s$ can be chosen so that
$\rho_m$ is $(\SAut (X), \SAut_{Z_s} (P_m))$-comparable and, hence, by Proposition \ref{cmo.p1}
we can suppose that $Y_1'=Y_2'$ and $\varphi'$ is the identity map.

By Proposition \ref{cmo.p1} it suffices to establish now that $\rho_m$ is holomorphically comparable
on a family of algebraic varieties containing $Y_i$. 
By Proposition \ref{qub.p1} $ \rho_m$ is a principal $H_{n-m}'$-bundle. 
Consider the morphism $\theta : Y_1' \to H_{n-m}'$ as in Proposition \ref{cmo.p2} (with $H$ and $Y'$ replaced by
$H_{n-m}'$ and $Y_1'$).
In order  to 
prove holomorphic comparability it suffices to show that $\theta$ extends to a holomorphic map $\Theta : P_m \to H_{n-m}'$.
Note that as an affine variety $H_{n-m}'$ is isomorphic to the direct product of $\C^{m(n-m)}$ and $\SL(n-m,\C)$.
Hence $\pi_1 (H_{n-m}')=\pi_2(H_{n-m}')=0$ and the existence of an extension $\Theta$ is provided by Proposition \ref{hex.p2} 
which yields the desired conclusion in (i).  

For (ii) exactly the same argument works
with Theorem \ref{cgw.t1} (a) replaced by Theorem \ref{cgw.t1} (b).

For (iii) we recall that $\rho : X \to Q$ is a principal $H''$-bundle by Proposition \ref{qub.p1}.
Arguing as in the proof of Proposition \ref{csl.p1a} we can find a subgroup $F\subset \SAut (X/Q)$
containing $H''$ such that every fiber of $\rho$ is $F$-flexible. Hence, by Corollary \ref{gp1.c1}
we can suppose that each $\rho|_{Y_i}: Y_i \to Q$ is a closed embedding. Let $Y_i''=\rho (Y_i)$
and $\varphi'' : Y_1'' \to Y_2''$ be the isomorphism induced by $\varphi$. Since $\ED (Y_1) \leq 2n-3$
we see that $\varphi''$ extends to an automorphism $\alpha'' \in \SAut (Q)$ by Theorem \ref{cqu.t1}.
By Proposition \ref{qub.p2}  $\rho$ is $(\SAut (X), \SAut (Q))$-comparable.  Hence, by Proposition \ref{cmo.p1}
we can suppose that $Y_1''=Y_2''$ and $\varphi''$ is the identity map.
As before we have a morphism $\theta'' : Y_1'' \to H''$ as in Proposition \ref{cmo.p2} which extends
to a holomorphic map $\Theta'' : Q \to H''$. Thus we have holomorphic comparability of $\rho$
which yields (iii) and concludes the proof.
\eproof

\bcor\label{csl.c2a}
Let $X=\SL (n,\C)$ and
$\varphi : Y_1\to Y_2$ be an isomorphism of two smooth closed subvarieties of $X$  such that 
$\dim Y_i \leq \frac{n}{3}-1$, $H_i (Y_1)=0$ for $i \geq 3$ 
and $H_2 (Y_1)$ is a free abelian group.
Then there exists a holomorphic automorphism $\beta$ of $X$ such that  $\beta|_{Y_1} =\varphi$.

\ecor

\bproof Since $\dim Y_i \leq \frac{n}{3}-1$ the smoothness assumption implies that
$\ED (Y_i)+\dim Y_i \leq n-2$. Hence, Theorem \ref{csl.t2} implies the desired conclusion.
\eproof

\bthm\label{csl.t1} Let $\varphi : Y_1\to Y_2$ be an isomorphism of two closed
 subvarieties of $X\simeq \SL (n,\bk )$ with $n \geq 3$ such that $Y_i$ is isomorphic to $\A_\bk^k$. Suppose that 
 either $k\leq \frac{n}{3}-1$ or $k=1$.
Then there exists $ \alpha \in \SAut (X)$ such that $\alpha|_{Y_1} =\varphi$.\footnote{The case of $k=1$ is, of course,
the theorem of Van Santen \cite{St} which we prove by other means.}

\ethm

\bproof  Let $k\leq \frac{n}{3}-1$, $m=\ED (Y_1)$ and let $Y_i'= \rho_m (Y_i), \, i=1,2$. Repeating the argument in the proof of Theorem \ref{csl.t2} we can suppose
that $\rho_m|_{Y_i} : Y_i \to Y_i'$ is a closed embedding, $Y_1'=Y_2'$ and the induced isomorphism $\varphi' : Y_1' \to Y_2'$
is the identity map. 
Since $ \rho_m$ is a principal $H_{n-m}'$-bundle,
$\rho_m$ is comparable on the family $\cF_{\rm aff}$ as in Corollary \ref{cmo.c1}.
The desired conclusion follows now from Proposition \ref{cmo.p1}.

Similarly, if $k=1$ then as in Theorem \ref{csl.t2} (iii) we can suppose that $\rho|_{Y_i}: Y_i \to Q$ is a closed embedding
and for $Y_i''= \rho(Y_i)$ the induced isomorphism $\varphi'' : Y_1''\to Y_2''$ is the identity map.
Since $ \rho$ is a principal $H''$-bundle,
$\rho$ is comparable on the family $\cF_{\rm aff}$.
The desired conclusion follows again from Proposition \ref{cmo.p1}.
\eproof

The following necessary condition for the positive solution of the extension problem is straightforward.

\bprop\label{tob.p1} Let a group $G\subset \Aut (W)$ act on an algebraic variety $W$ 
and $\varphi: Y_1\to Y_2$ be an isomorphism of
closed subvarieties of $W$. Suppose that $\iota_k : Y_k \hookrightarrow W$ is the natural embedding.
Suppose also that for every $\alpha \in  G$ there exists $k\geq 1$
such that the homomorphism $\pi_k (Y_1) \to \pi_k (W)$ induced by $\alpha \circ \iota_1$
is different from the similar homomorphism induced by $\iota_2$.
Then $\varphi$ cannot be extended to an automorphism from  $G$.
\eprop

\brem\label{tob.r1} In the complex case note that if  $\iota_1 : Y_1 \hookrightarrow W$  induces 
a trivial homomorphism $\pi_k (Y_1) \to \pi_k(X)$ for some $k$
while the similar homomorphism induced by  $\iota_2 : Y_2 \hookrightarrow W$ is nontrivial
then $\varphi$ cannot be extended even to a homeomorphism of $W$. 
Furthermore,  if $G$ is contained in the connected component of identity in $\Aut (X)$
then the extension problem does not have a positive solution for the group $G$ if there is no  homotopy
of $\varphi$ to the identity map via closed embeddings of $Y_1$ into $W$. 
\erem

In the rest of this section we present a concrete (and more or less obvious) example illustrating Proposition \ref{tob.p1} in the case of $X \simeq \SL(n,\C)$.

\blem\label{tob.l2}  There is a closed embedding of $Y \simeq \SL (2,\C)$ into $X\simeq \SL (n,\C)$
such that it generates an isomorphism $\pi_3 (Y)\simeq \pi_3 (X)$ of the homotopy groups.
\elem

\bproof Let $\rho : X\to Q$ be as in Notation \ref{qub.n1}, i.e., by Proposition \ref{qub.p1}
 it is a principal $H''$-bundle with  fiber $F\simeq \SL (n-1,\C)$. Since $Q$ is a complexification of a real $(2n-1)$-dimensional sphere
 it has a homotopy type of this sphere\footnote{Actually, it is well-known that $Q$ is diffeomorphic to the tangent bundle of
 the sphere.}.
Hence $\pi_k (Q)=0$ for $1\leq k \leq 2n-2$ and $\pi_{2n-1}(Q)=\Z$.
Now the exact homotopy sequence for
the fiber bundle $\rho : X \to Q$ implies that the natural embedding of
 $\SL (n-1,\C)\simeq F \hookrightarrow X\simeq \SL (n,\C)$
induces an isomorphism  $\pi_k (X)\simeq \pi_k (F)$
for $k< 2n-2$. Choosing a natural sequence 
$\SL (2,\C) \subset \SL (3,\C) \subset \ldots \subset \SL (n-1,\C) \subset \SL (n,\C)$
of closed embeddings we get the desired conclusion.
\eproof

\bthm\label{tob.t1} Let $X$ be an affine algebraic variety isomorphic to $\SL( n,\C)$ with $n \geq 3$.
There are two closed subvarieties $Y_1$ and $Y_2$ in $X$ isomorphic to $\SL (2,\C)$ and such that
there is no automorphism $\alpha$ of $X$ for which $\alpha (Y_1)=Y_2$.
\ethm

\bproof By Lemma \ref{tob.l2} we can suppose that the natural embedding $Y_1 \hookrightarrow X$
generates an isomorphism $\pi_3 (Y_1) \simeq \pi_3 (X)$. By Proposition \ref{tob.p1} and Remark \ref{tob.r1}
in order to prove Theorem \ref{tob.t1} it suffices to present an embedding $Y_2 \hookrightarrow X$ such that the induced homomorphism
 sends $\pi_3 (Y_2)$ into the zero element of $\pi_3 (X)$. 
Treat a point in $X$ as a matrix $A=[a_{i,j}]$. Let $A'$ be the $(n-1)\times (n-1)$ matrix
obtained from $A$ by removing the first row and the first column. 
Consider the subvariety $X_0$ of $X$ that consists of matrices $A$ such that $A'$ is the identity matrix.
Note that $X_0$ is naturally isomorphic to $\C^{2n-1}$ with coordinates $a_{1,2}, \ldots a_{1n}, a_{21}, \ldots , a_{n1}$
since $a_{1,1}$ can be expressed as function of these coordinates because of the equation $\det A=1$.
Choose in $X_0$ the quadric $Y_2\simeq \SL (2,\C)$ given by $a_{1,2}a_{1,3}-a_{2,1}a_{3,1}=1$ and $a_{1,k}=a_{k,1}=0$ for all $k\geq 4$. 
Then the embedding $Y_2 \hookrightarrow X$ induces the zero map  $\pi_3 (Y_2) \to \pi_3 (X)$
since it factors through $Y_2 \hookrightarrow X_0$ and $X_0$ is contractible. This yields the desired conclusion.
\eproof


\begin{thebibliography}{KaMi} 

\bibitem{AbMo} S. Abhyankar, T.-T. Moh, {\em Embeddings of the line in the plane}, J. Reine Angew. Math. {\bf 276} (1975), 148-166.

\bibitem{AFKKZ} I.~V.~Arzhantsev, H.~Flenner, S.~Kaliman,
F.~Kutzschebauch, M.~Zaidenberg:
{\em Flexible varieties and automorphism groups}.  Duke Math.\ J.\ 162:4 (2013), 767--823.

\bibitem{AM} M.~F.~ Atiyah, I.~G.~MacDonald, {\em Introduction to Commutative Algebra}, University of Oxford, Addison-Wesley Publishing Company, 1969. 
 
\bibitem{BS}  J.~Blanc, Immanuel van Santen (n\'e Stampfli), {\em Embeddings of affine spaces into quadrics}, preprint,  mathAG. 

\bibitem{BuHu} G.T. Buzzard, G.T., J.H. Hubbard, {\em A Fatou?Bieberbach domain avoiding a neighborhood of a variety of codimension 2}
Math. Ann. {\bf 316}(4), 699-702 (2000).

\bibitem{Cr}  P. C. Craighero, {\em A result on m-flats in $\A_\bk^n$}, Rend. Sem. Mat. Univ. Padova 75 (1986), 39-46.


\bibitem{Ei} D.~Eisenbud: {\em  Commutative algebra. With a view toward algebraic geometry.}
Graduate Texts in Mathematics, 150. Springer-Verlag, New York, 1995. xvi+785 pp.
%


\bibitem{FS} P.~Feller, I.~ van S. Stampfli, {\em Uniqueness of Embeddings of the Affine Line into Algebraic Groups},  preprint, math.AG.arXiv:1609.02113.

\bibitem{FKZ-GW} H.~Flenner, S.~Kaliman, and M.~Zaidenberg, {\em A Gromov-Winkelmann type theorem for flexible varieties},  J. Eur. Math. Soc. (JEMS) 18 (2016), no. 11, 2483-2510. 
%
\bibitem{For01}  F.~Forstneri\^c,  {\em On complete intersections}, Ann. Inst. Fourier {\bf 51}(2), 497-512 (2001).

\bibitem{For} F.~Forstneri\^c, {\em Stein Manifolds and Holomorphic Mappings.
The Homotopy Principle in Complex Analysis}, Springer-Verlag, Berlin-Heidelberg,  2011.



\bibitem{FR}  F. Forstneri\^c, J.-P. Rosay, {\em Approximation
of biholomorphic mappings by automorphisms of $\C^n$}, Invent. Math. 112:2 (1993), 323-349.

\bibitem{Fre} G. ~Freudenburg, {\em Algebraic Theory of Locally Nilpotent Derivations}  (Encyclopaedia of Mathematical Sciences)  Springer Berlin-Heidelberg-New York, 2006.

\bibitem{Gr1} M.~Gromov: {\em Partial differential relations.}
Ergebnisse der Mathematik und ihrer Grenzgebiete (3) 9.
Springer-Verlag, Berlin, 1986.

\bibitem{EGA} A.  Grothendieck, J.  Dieudonn\'e,  {\em \'El\'ements de g\'eom\'etrie algŽbrique. IV. \'Etude locale des sch\'emas et des morphismes de sch\'emas. III.} 
(French) Inst. Hautes \'Etudes Sci. Publ. Math. {\bf 28}  (1966), 5-255.

\bibitem{Har} R.~Hartshorne, {\em Algebraic Geometry},
Springer-Verlag, New York-Heidelberg, 1977.

\bibitem{Hir} M. Hirsch, {\em Smooth regular neighborhoods}, Ann. of Math. (2) {\bf 76} (1962) 524-530.

\bibitem{Je}  Z. Jelonek, {\em The extension of regular and rational embeddings}, Math. Ann. 277 (1987), no. 1, 113Ð120.

\bibitem{Ka91} S.~ Kaliman, {\em Extensions of isomorphisms between affine algebraic subvarieties of $k^n$ to automorphisms of $k^n$}, Proc. Amer. Math. Soc. 113 (1991), no. 2, 325-334.

 
\bibitem{KaKu08} S.~ Kaliman, F.~ Kutzschebauch, {\em Criteria for the density property of complex manifolds}, Invent. Math. 172 (2008), no. 1, 71--87. 

\bibitem{KaKuTr} S.~ Kaliman, F.~ Kutzschebauch, T.~ T.~ Truong, {\em On subelliptic manifolds}, preprint,  19 p., arXiv:1611.01311, to appear in Israel J. of Math. 

\bibitem{Kl} S.~L.~ Kleiman, {\em The transversality of a general translate},  Compositio Math. {\bf 28} (1974), 287-297. 

\bibitem{KrRu}  H.~Kraft, P.~ Russell, {\em Families of group actions, generic isotriviality, and linearization},
 Transform. Groups 1{\bf 9} (2014), no. 3, 779-792.

\bibitem{Lo}  S. Lojasiewicz, {\em Triangulation of semi-analytic sets}, Ann. Scuola Norm. Sup. Pisa (3) {\bf 18} (1964), 449-474. 

\bibitem{Ma}  H.~Matsumura,  {\em Commutative Algebra}, Reading, Massachusetts The Benjamin/Cummings Publishing Company, 1980.

 \bibitem{Oni}     A.~L.~ Onishchik, {\em Methods of sheaf theory and Stein spaces} (Russian).
Current problems in mathematics. Fundamental directions, Vol. 10 (Russian), 5-73, 283, Itogi Nauki i Tekhniki, Akad. Nauk SSSR, Vsesoyuz. Inst. Nauchn. i Tekhn. Inform., Moscow, 1986. 


\bibitem{PV}
V.~L.~Popov,  E.~B.~Vinberg: {\em Invariant Theory.} In: Algebraic
geometry IV,  A.\ N.\ Parshin, I.\ R.\ Shafarevich (eds.), Berlin,
Heidelberg, New York: Springer-Verlag, 1994.
%
\bibitem{Pr} C. ~Procesi: Lie groups. An approach through invariants
and representations.
Universitext. Springer, New York, 2007.
%
\bibitem{Qui} D.~ Quillen, {\em Projective modules over polynomial rings},  Invent. Math. {\bf 36} (1976), 167-171. 
%
\bibitem{Ra1} C.~P.~Ramanujam, {\em A note on automorphism groups of algebraic varieties}, Math.\ Ann.\
156 (1964), 25--33.

\bibitem{RoRu}  J.-P. Rosay, W. Rudin, {\em Holomorphic maps from $\C^n$ to $\C^n$}, Trans. Amer. Math. Soc. {\bf 310} (1988), no. 1, 47-86.

\bibitem{Sh} I. Shafarevich, {\em Basic algebraic geometry. 1. Varieties in projective space.} Second edition. Translated from the 1988 Russian edition and with notes by Miles 
Reid. Springer-Verlag, Berlin, 1994.

\bibitem{St}  I. ~Stampfli, {\em Algebraic embeddings of $\C$ into $\SL_n(\C)$}, Transform. Groups 22 (2017), no. 2, 525-535.


\bibitem{Sus} A.~Suslin, {\em Projective modules over polynomial rings are free}, (Russian) Dokl. Akad. Nauk SSSR {\bf 229} (1976), no. 5, 1063-1066.

\bibitem{Su} M. Suzuki, {\em Propi\'et\'es topologiques des polynomes de deux variables complexes, et automorphismes alg\'earigue de l'espace $\C^2$}, J. Math. Soc. Japan,{\bf  26} (1974), 241-257.

\bibitem{VS} V. Srinivas, {\em  On the embedding dimension of an affine variety},  Math. Ann., 289(1):125-132, 1991.


\bibitem{Wi} J.~Winkelmann:
{\em On automorphisms of complements of analytic subsets in $\CC^n$}.
Math.\ Z.\  204  (1990),  117--127.


\end{thebibliography}
\end{document}